\def\real{{\tt I\kern-.2em{R}}}
\def\nat{{\tt I\kern-.2em{N}}}

\def\realp#1{{\tt I\kern-.2em{R}}^#1}
\def\natp#1{{\tt I\kern-.2em{N}}^#1}
\def\hyper#1{\ ^*\kern-.2em{#1}}
\def\monad#1{\mu (#1)}

\def\st#1{{\tt st}(#1)}
\def\hyperreal{{^*{\real}}}
\def\hyperrealp#1{{\tt ^*{I\kern-.2em{R}}}^#1} 
\def\hypernat{{^*{\nat }}}
\def\hypernatp#1{{{^*{{\tt I\kern-.2em{N}}}}}^#1} 
\def\eskip{\hskip.25em\relax}

\def\Hyper#1{\hyper {\eskip #1}}
\def\leaderfill{\leaders\hbox to 1em{\hss.\hss}\hfill}
\def\srealp#1{{\rm I\kern-.2em{R}}^#1}
\def\sp{\vert\vert\vert} 

\def\power#1{{{\cal P}(#1)}}
\def\iff{\leftrightarrow}
\def\qed{{\vrule height6pt width3pt depth2pt}\par\medskip}
\def\pars{\par\smallskip}
\def\parm{\par\medskip}

\def\b#1{{\bf #1}}
\def\r#1{{\rm #1}}
\def\sig{{^\sigma}}
\magnification=\magstep0
\tolerance 10000
\hoffset=.64in
\hsize 5.00 true in
\vsize 8.50 true in
\baselineskip=14pt
\nopagenumbers
 at 18truept
 at 14truept
\pageno=63
\headline={\ifnum\pageno=59\hfil\quad \hfil \else\ifodd\pageno\rightheadline 
\else\leftheadline\fi\fi}
\def\rightheadline{\hfil \bf The Theory Ultralogics\hfil\tenbf\folio}
\def\leftheadline{{\tenbf\folio} \hfil\bf Robert A. Herrmann \hfil}
\voffset=1\baselineskip
\medskip
\centerline{{\bf 7. DEVELOPMENTAL PARADIGMS} \footnote*{Chapters 7-11 last revised 15 OCT 2012}}\par
\bigskip
\leftline{\bf 7.1 Introduction.}
\medskip
Consider the real line. If you believe that time is the ordinary continuum$,$ then the 
entire real line can be your {time line}. Otherwise$,$ you may consider only a 
subset of the real line as a time line. In the original version of this section$,$ the time concept for the
MA-model was presented in a unnecessarily complex form. As shown in [3]$,$ one can assume an
absolute substratum time within the NSP-world. It is the infinitesimal light-clock time measures
that may be altered by physical processes. In my view$,$ the theory of quantum electrodynamics
would not exist without such a NSP-world time concept.\pars 
Consider a small interval $[a,b),\ a<b$ as our basic time interval where as the real numbers
increase the time is intuitively considered to be increasing. In the following approach$,$ one may
apply the concept of the persistence of mental version relative to descriptions for the behavior of a
Natural (i.e. physical) system at a moment of time within this interval. An exceptionally small subinterval can be
chosen within $[a,b)$ as a maximum subinterval length $=M.$ 
``Time'' and the size of a ``time'' interval as they are used in this and the following sections refer to an intuitive concept used to aid in comprehending the notation of an event sequence. [See below.]
First$,$ let $a = t_0.$ Then choose $t_1$ such
that $a< t_1 <b.$
There is a partition $t_1, \ldots, t_m$ of $[a,b)$ such that $t_0<t_1< \cdots < t_m < b$ and
$t_{j+1}-t_j \leq M.$ The final subinterval $[t_m, b)$ is now separated$,$ by induction$,$ say be
taking midpoints$,$ into an increasing sequence of times $\{t_q\}$ such that $t_m <t_q <b$
for each $q$ and 
$\lim_{q \to \infty}\ t_q = b.$\pars  
Assume the prototype $[a,b)$ with the time subintervals as defined above. Let $[t_j, t_{j+1})$ be
any of the time subintervals in 
$[a,b).$  For each such subinterval$,$ let ${\r W}_i$ denote the readable 
sentence\parm
\line{\hfil
This$\sp$frozen$\sp$segment$\sp$gives$\sp$a$\sp$description$\sp$for$\sp$the$\sp$\hfil}\pars
\line{\hfil
time$\sp$interval$\sp$that$\sp$has$\sp$as$\sp$its$\sp$leftmost$\sp$endpoint$\sp$the \hfil}\pars
\line{\hfil
$\sp$time$\sp\lceil t_i\rceil\sp$that$\sp$corresponds$\sp$to$\sp$the$\sp$natural$\sp$number$\sp$i.\hfil}\parm 
\noindent Let ${\r T}_i = \{ \r x\r W_i\mid  \r x \in {\cal W}\}.$ The set ${\r T}_i$ is called 
a {{\it totality}} and each member of any such ${\r T}_i$ is called a
{{\it frozen segment.}} Notice that since
 the empty word is not a member of ${\cal W},$ then the cardinality of each 
member of ${\r T}_i$ is greater than that of ${\r W}_i.$ Each ${\rm T}_i$ is 
a (Dedekind) denumerable set$,$ and if $i \not=j,$ then ${\rm T}_i \cap {\rm T}_j 
= \emptyset.$\pars
\hrule\smallskip\hrule \smallskip
It is obvious that the concept of  ``time''  need not be the 
underlying interpretation for these intervals. Time simply refers to an 
external event ordering concept. 
For other purposes$,$ simply 
call these intervals  ``event intervals.'' In the above descriptions for  
$\r W_i,$ simply replace ``time$\sp$interval'' with ``event$\sp$interval'' 
and replace the second instance of the word ``time'' with the word  
``event.'' If this interpretation is made$,$ then other compatible  
interpretations would be necessary when applying a few of the following 
results. \pars \hrule\smallskip\hrule\smallskip
I point out two minor aspects of the above constructions. First$,$5\noindent within 
certain descriptions there 
are often  ``symbols'' used for real$,$ complex$,$ natural numbers etc. These 
objects also exist as abstract objects within the structure ${\cal M}.$ No 
inconsistent interpretations should occur when these objects are 
specifically modeled within ${\cal M}$ since to my knowledge all of 
the usual mathematical objects used within physical analysis are disjoint 
from ${\cal E}$ as well as disjoint from any finite Cartesian product of
${\cal E}$ with itself. If for future research within physical applications 
finite partial sequences of natural numbers and the finite equivalence classes 
that
\noindent appear in ${\cal E}$ are needed and are combined into one model for 
different purposes than the study of 
descriptions$,$ then certain modifications would need to be made so that 
interpretations would remain consistent. Secondly$,$ I have tried whenever 
intuitive strings are used or sets of such strings are 
defined to 
use Roman letter notation for such objects. This  
only applies for the intuitive model. Also $\r W_i$ is only an identifier and may be altered.\par
\bigskip
\leftline{\bf 7.2 Developmental Paradigms}
\medskip
It is clear that if one considers a time interval of the type $(-
\infty,+\infty),$ $(-\infty,b)$ or $[a, +\infty),$ 
then  each of these may be considered as the union of a denumerable collection of time 
intervals of the type $[a,b)$ with common endpoint names displayed. Further$,$ although 
$[a,b)$ is to be considered as subdivided into denumerably 
many subintervals$,$ it is not necessary that each of the time intervals $[t_j,t_{j+1})\subset [a,b)$ 
be accorded a corresponding description for the appearance of a specific 
Natural system that is distinct from all others that occur throughout the time subinterval.
Repeated descriptions only containing a different last natural number i  in the next to last position
will suffice. Each basic developmental paradigm will be restricted$,$ at present$,$ to such a time
interval $[a,b).$\pars
 Where human 
perception and descriptive ability is concerned$,$ the least controversial 
approach would be to consider only finitely many descriptive choices as 
appropriate. A finite set is recursive and such a choice$,$ since the result is 
such a set$,$ would be considered to be the simplest type of algorithm. You ``simply'' 
check to see if an expression is a member of such a finite set.   
If we limited ourselves to finitely many human choices for Natural 
system descriptions 
from the set of all totalities and did not allow a denumerable or a continuum 
set to be chosen$,$ then the next result establishes that within the 
{Nonstandard Physical world (i.e. NSP-world)} 
such a finite-type of choice can be applied and a continuum of descriptions 
obtained.\pars
The following theorem is not insignificant even if we are willing to accept a 
denumerable set of distinct descriptions --- descriptions that are not only 
distinct in the next to the last symbol$,$ but are also distinctly different in 
other aspects as well. For$,$ if this is the case$,$ the results of Theorem 
7.2.1 still apply. The same finite-type of process in the NSP-world yields 
such a denumerable set as well. \pars
The  term  ``NSP-world'' will signify a certain second type of interpretation 
for nonstandard entities. In particular$,$ the subtle logics$,$ unreadable 
sentences$,$ etc. This 
interpretation will be developed throughout the remainder of this book. One 
important aspect of how descriptions are to be interpreted is that a 
description 
correlates directly to an assumed or observed real Natural phenomenon$,$ and 
conversely. In these investigations$,$ the phenomenon is called an {{\it 
event.}}\par
In order to simplify matters a bit$,$ the following notation is employed. 
Let ${\cal T}= \{{\r T}_i\mid i \in
\nat \}.$ Let $F({\cal T})$ be the set of all {\bf nonempty} and {\it finite}
subsets of ${\cal T}.$ This symbol has been used previously to include the 
empty set$,$ this set is now excluded. Now let ${\r A} \in F({\cal T}).$ 
Then there exists a finite choice set \ s \ such that $\r x \in \r s$ iff
there exists a unique ${\r T}_i \in {\r A}$ and $ \r x \in {\r T}_i.$ Now let the 
set $\cal C$ denote the set of all such finite choice sets. As to interpreting these results within 
the NSP-world$,$ the following is essential. Within nonstandard analysis the 
term  ``hyper'' is often used for the result of the * map. For example$,$ you 
have $\hyperreal$ as the hyperreals since $\real$ is termed the real numbers. 
For certain$,$ but not all concepts$,$ the term ``hyper'' or the corresponding 
* notation will be universally replaced by the term ``ultra.'' Thus$,$ certain 
purely  subtle 
words or *-words become ``ultrawords'' within the developmental paradigm 
interpretation. [Note: such a word was previously called a superword.] Of 
course$,$ for other scientific or philosophical systems$,$ such abstract 
mathematical objects can be reinterpreted by an appropriate technical term 
taken from those disciplines. \pars
As usual$,$ we are working within any enlargement and all of the 
above intuitive objects are embedded into the G-structure. Recall$,$ that to 
simplify 
expressions$,$ we often suppress within our first-order statements a specific 
superstructure element that bounds a specific quantifier. 
The alphabet ${\cal A}$ is now assumed to be countable.\pars

[Note 2 MAY 1998: The material between the [[ and the ]] has been altered from the original that appeared in the 1993 revision.] [[[Although theorem 7.2.1 may be significant$,$ it is no longer used for the other portions of this research. The set of all developmental paradigms corresponds to the set of all choice functions define on $\cal T$. Also see http://arxiv.org/abs/math/0605120  \parm
 {\bf Theorem 7.2.1} 
{\sl Let $\emptyset \not= \zeta \subset \nat$ and $\widetilde{\cal T} = \{\b T_i 
\mid i \in \zeta \}.$ There exists a set of 
sets $\cal S$ determined by hyperfinite set $Q$ and hyper
finite choice defined on $Q$ such that:\pars {\rm 
(i )} $s^\prime  \in \cal S$ iff for each $ \b T \in \widetilde{\cal T}$ there 
is one and only one $[g] \in \Hyper {\b T}$ such that $[g] \in s^\prime,$  
 and if $x \in s',$ then there is some $ \b T \in \widetilde{\cal T}$ and some $[g] \in \Hyper {\b T}$ such that $x = [g].$ 
{\rm (}If 
$\Hyper [g] \in {^\sigma {\b T}},$ then $[g] = [f] \in \b T.${\rm)}} \pars 
Proof. (i) Let $A \in F(\widetilde{\cal T}).$ Then from the definition of 
$\widetilde{\cal T},$ there exists some $ n \in \nat$ such that $A = \{
\b T_{j_i} \mid i = 0, \ldots, n \land j_i \in \nat \}.$ From the definition of 
$\b T_k,$ each $\b T_k$ is denumerable. Notice that any $[f] \in \b T_k$ is 
associated with a unique member of $A_1.$. Simply consider the unique $f_0 \in 
[f].$ The unique member of $A_1$ is by definition $f_0(0).$ Thus each member of 
$\b T_k$ can be specifically identified. Hence$,$ for each $\b T_i$ there is a 
denumerable $M_i \subset \nat$ and a bijection $h_i\colon M_i \to \b T_i$ such 
that $a_i \in \b T_i$ iff there is a $k_i \in M_i$ and $h_i(k_i) = a_i.$ 
Consequently$,$ for each $i = 0,\ldots,n$ and $a_{j_i} \in \b T_{j_i},$ we have 
that $h_{j_i}(k_{j_i}) = a_{j_i},$ and conversely for each $i =0,\dots,n$ and 
$k_{j_i} \in M_{j_i},\ h_{j_i}(k_{j_i}) \in \b T_{j_i}.$ Obviously$,$ 
$\{h_{j_i}(k_{j_i}) \mid i = 0, \ldots, n \}$ is a finite choice set.
All of the above may be translated into the following sentence that holds in 
${\cal M}.$ (Note: Choice sets are usually considered as the range of choice functions. Further$,$ ``bounded formula simplification'' has been used.)\pars
\line{\hfil $ \forall y(y \in F(\widetilde{\cal T}) \to \exists s((s \in \power {{\cal E}})\land \forall 
x((x \in y) \to \exists z((z \in x) \land (z \in s) \land $\hfil}\pars
\line{\hfil $\forall w(w \in {\cal E} \to ((w  \in s) \land(w \in x) \iff (w = z)))))\land$\hfil}\pars
\line{(7.2.1)\hfil$\forall u(u \in {\cal E} \to ((u \in s) \iff \exists x_1((x_1 \in y)\land (u \in x_1))))))$\hfil }\parm
\noindent For each $A \in F(\widetilde{\cal T}),$ let $S_A$ be the set of all 
such choice sets generated by the predicate that follows the first $\to$ formed from (7.2.1) by deleting 
the $\exists s$ and letting $y = A.$ Of course$,$ this set exists within our set 
theory. Now let ${\cal C} = \{S_A \mid A \in F(\widetilde{\cal T})\}.$ \pars
Consider $\Hyper {\cal C}$ and $\Hyper (S_A).$ Then $s \in \Hyper (S_A)$ iff 
$s$ satisfies (7.2.1) as interpreted in $\Hyper {\cal M}.$ Since we are 
working in an enlargement$,$ there exists an internal $Q \in \Hyper 
(F(\widetilde{\cal T}))$ such that $^\sigma{\widetilde{\cal T}} \subset Q 
\subset \Hyper {\widetilde{\cal T}}.$ Recall that $^\sigma{\widetilde{\cal T}} 
=\{\Hyper {\b T} \mid \b T \in \widetilde{\cal T}\}.$ Also $^\sigma{\b T} 
\subset \Hyper {\b T}$ for each $\b T \in \widetilde{\cal T}.$  From the 
definition of $\Hyper {\cal C},$ there is an internal set $S_Q$ and  $s \in  
S_Q$ iff $s$ satisfies the internal defining predicate for members of $S_Q$
and this set is the set of all such $s.$ $(\Rightarrow)$ Consequently$,$ 
since for each $\b T\in \widetilde{\cal T},\ \Hyper {\b T} \in Q$$,$ then the 
generally external $s^\prime = \{ s\cap \Hyper {\b T}\mid \b T \in 
\widetilde{\cal T}\}$ satisfies the $\Rightarrow$ for (i). Note$,$ however$,$ that for 
$\Hyper {\b T},\ \b T \in \widetilde{\cal T},$ it is possible that 
$s \cap \Hyper {\b T} = \{\Hyper {[f]}\}$ and $\Hyper {[f]}\in {^\sigma{\b 
T}}.$ In this case$,$ by the finiteness of $[f]$ it follows that $[f] = \Hyper 
{[f]}$ implies that $s \cap\Hyper {\b T} = \{[f]\}.$ Now let
${\cal S} = \{s^\prime \mid s \in S_Q\}.$ In general$,$ $\cal S$ is an external 
object. \pars
($\Leftarrow$) Consider the internal set $S_Q.$  Let $s^\prime$ be the set as defined by the right-hand side of (i). For each internal $x \in s^\prime$
and applying$,$ if necessary$,$ the *-axiom of choice for *-finite sets$,$ we have the internal set $A_x = \{y\mid (y \in S_Q)\land (x \in y)\}$ is nonempty. The set $\{A_x\mid x \in s^\prime\}$ has the finite intersection property. For$,$ let nonempty internal $B=\{x_1,\ldots, x_n\}.$ Then the set $A_B =\{y\mid (y \in S_Q)\land (x_1 \in y)\cdots \land (x_n \in y)\}$ is internal and nonempty by the *-axiom of choice for *-finite sets. Since we are in an enlargement and $s^\prime$ is countable$,$ then $D=\bigcap\{A_x\mid x \in s^\prime\} \not= \emptyset.$ Now take any $s\in D.$ Then $s \in S_Q$ and from the definition of ${\cal S},\ s^\prime \in {\cal S}$. This completes the proof.\qed 
[Note: Theorem 7.2.1 may be used to model physical developmental paradigms associated with event sequences.]
\par
Although it is not necessary$,$ for this particular investigation$,$ the set $\cal S$ may  be considered               
{{\it a set of all developmental paradigms.}}  Apparently$,$  $\cal S$ 
contains every possible developmental paradigm for all possible frozen segments
and $\cal S$ contains paradigms for any *-totality $\Hyper {\b T}.$ There are 
*-frozen segments contained in various $s^\prime$ that can be assumed to be 
unreadable sentences since $^\sigma{\b T} \not= \Hyper {\b T}.$]]\pars
Let $A \in F(\widetilde{\cal T})$ and $M(A)$ be a subset of $S_A$ for which 
there exists a written set of rules that selects some specific member of 
$S_A.$ Obviously$,$ this may be modeled by means of functional relations. 
First$,$ $M(A) \subset S_A$ and it follows$,$ from the difference in 
cardinalities$,$ that 
there are infinitely many members of $\Hyper (S_A)$ for which there does not 
exist a readable rule that will select such members. However$,$ this does not 
preclude the possibility that there is a set of purely unreadable sentences 
that do determine a specific member of 
$\hyper {S_A} - {^\sigma M(A)}.$ This might come about in the following 
manner. Suppose that $H$ is an infinite set of formal sentences that is 
interpreted to be a set of rules for the selection of distinct members of $M(A).$ 
Suppose we have a bijection $h\colon M(A) \to \b H$ that represents this 
selection process. Let $\Hyper {\cal M}$ be at least a polysaturated 
enlargement of $\cal M,$ and consider $^\sigma f: {^\sigma(M(A))} \to {^\sigma 
{\b H}}.$ The map $^\sigma f$ is also a bijection and $^\sigma f: 
{^\sigma(M(A))} \to \Hyper {\b H}.$ Since $\vert {^\sigma(M(A))} \vert <
\vert {\cal M} \vert,$ it is  well-known that there exists an internal map 
$h\colon A^\prime \to \Hyper {\b H}$ such that $h\mid {^\sigma (M(A))} = {^\sigma 
f},$
 and $A^\prime,\ h[A^\prime]$ are internal. Further$,$ for internal $A^\prime \cap 
\Hyper {(S_A)} =B,\ {^\sigma (M(A))} \subset B.$ However$,$ $^\sigma (M(A))$ is 
external.  This yields that $h$ is defined on $B$ and 
$B\cap(\Hyper {S_A} - {^\sigma (M(A))}) \not= \emptyset.$ Also$,$ $^\sigma {\b 
H} \subset h[B] \subset \Hyper {\b H}$ implies$,$ since $h[B]$ is 
internal$,$ that $^\sigma H \not= h[B].$ Consequently$,$ in this case$,$ $h[B]$ may 
be interpreted as a set of *-rules
that determine the selection of members of 
$B.$ That is to say that there is some $[g] \in h[B] - {^\sigma H}$ and a $[k] 
\in \Hyper {S_A} - {^\sigma (M(A))}$ such that $([k],[g]) \in h.$  As it will 
be shown in the next section$,$ the set $H$ can be so constructed that if $[g] 
\in h[B] - {^\sigma H},$ then $[g]$ is unreadable. \par
\bigskip
\leftline{\bf 7.3 Ultrawords}
\medskip
Ordinary propositional logic is not compatible with deductive quantum logic$,$ 
intuitionistic logic$,$ among others. In this section$,$ a subsystem of 
propositional logic is investigated which rectifies this incompatibility. I 
remark that when a standard propositional language L or an informal language 
P isomorphic to L is considered$,$ it will always be the case that the L or P is  
minimal relative to its applications. This signifies that if L or P is 
employed  
in our investigation for a developmental paradigm$,$ then L or P is constructed 
only from those distinct propositional atoms that correspond to distinct 
members of d$,$ etc. The same minimizing process is always assumed for the 
following constructions. \pars
Let B be a formal or$,$ informal nonempty set of propositions. Construct the 
language $\r P_0$ in the usual manner from B (with superfluous parentheses 
removed) so that $\r P_0$ forms the smallest set of formulas that contains B
and such that $\r P_0$ is closed under the two binary operations $\land$ and 
$\to$ as they are formally or informally expressed. Of course$,$ this language 
may be constructed inductively or by letting $\r P_0$ be the intersection of 
all collections of such formula closed under $\land$ and $\to.$ \pars
We now define the deductive system $S.$ Assume substitutivity,
parenthesis 
reduction and the like.  Let $\r d 
= \{\r F_i \mid i \in \nat \}= \r B$ be a 
development paradigm$,$ where each $\r F_i$ is a 
readable frozen segment and describes the behavior of a Natural system over a 
time subinterval. Let the set of axioms be the schemata\parm
\line{\hskip 0.5in (1)\hfil $({\cal A} \land {\cal B}) \to {\cal A},\ {\cal A} \in \r B$ \hfil \hskip 0.5in }\pars
\line{\hskip 0.5in (2)\hfil $({\cal A} \land {\cal B}) \to {\cal B}$ \hfil \hskip 0.5in }\pars
\line{\hskip 0.5in (3)\hfil ${\cal A} \land ({\cal B} \land {\cal C}) \to ({\cal A}\land 
{\cal B}) \land {\cal C},$ \hfil \hskip 0.5in }\pars
\line{\hskip 0.5in (4)\hfil $ ({\cal A} \land {\cal B}) \land {\cal C} \to {\cal A} \land ({\cal 
B} \land {\cal C}).$ \hfil \hskip 0.5in }\parm
If $\r P_0$ is considered as informal$,$ which appears to be necessary for some 
applications$,$ where the parentheses are replaced by the concept of symbol 
strings being to the   ``left'' or   ``right'' of other symbol strings and the 
concept of strengths of connectives is used (i.e. $\r A \land \r B \to \r C$ 
means $((\r A \land \r B) \to \r C$)$,$ then axioms 3 --- 4 and the parentheses 
in (1) and (2) may be omitted. The one rule of inference is Modus Ponens (MP). 
Proofs or demonstrations from hypotheses $\Gamma$ contain finitely many steps$,$ 
hypotheses may be inserted as steps and the last step in the proof is either a 
theorem if $\Gamma = \emptyset$ or if $ \Gamma \not= \emptyset,$ then the last 
step is a consequence of ( a deduction from ) $\Gamma.$ Notice that repeated application of (4)
along with (MP) will allow all left parentheses to be shifted to the right with the exception of the
(suppressed) outermost left one. Thus this leads to the concept of left to right ordering of a
formula. This allows for the suppression of such parentheses. In all the following$,$ this
suppression will be done and replaced with formula left to right ordering.  \pars
For each $\Gamma \subset \r P_0,$  let $S(\Gamma)$ denote the set of all 
formal theorems and consequences obtained from the above defined system $S.$ 
Since hypotheses may be inserted$,$ for each $\Gamma \subset \r P_0,\ 
\Gamma \subset S(\Gamma) \subset \r P_0.$ This implies that $S(\Gamma) \subset 
S(S(\Gamma)).$ So$,$ let $\r A  \in S(S(\Gamma)).$ The general concept of 
combining together finitely many steps from various proofs to yield another 
formal proof leads to the result that $\r A \in S(\Gamma).$ Therefore$,$ 
$S(\Gamma) = S(S(\Gamma)).$ Finally$,$ the finite step requirement also yields 
the result that if $\r A \in S(\Gamma),$ then there exists a finite $\r F 
\subset \Gamma$ such that $\r A \in S(\r F).$ Consequently$,$ $S$ is a finitary 
consequence operator and observe that if $C$ is the propositional consequence 
operator$,$ then $S(\Gamma) {\subset\atop \not=} C(\Gamma).$ Of course$,$ we may 
now apply the nonstandard theory of consequence operators to $S.$ \pars
It is well-known that the axiom schemata chosen for $S$ are theorems in 
intuitionistic logic. Now consider quantum logic with the {Mittelstaedt 
conditional} $i_1(\r A,\r B) = \r A^\perp \lor (\r A \land \r B).$ [1]  Notice 
that $i_1(\r A \land  \r B,\r B)= (\r A \land \r B)^\perp \lor ((\r A \land \r B) 
\land \r B)= (\r A \land \r B)^\perp \lor (\r A \land  \r B)= I$ ( the upper 
unit.) Then $i_1(\r A \land \r B),\r A)= (\r A \land \r B)^\perp \lor
((\r A \land \r B) \land  \r A) = (\r A \land \r B)^\perp \lor (\r A \land \r 
B) = I;\ i_1((\r A \land \r B) \land \r C),\r A \land (\r B \land \r C)) = 
((\r A \land \r B) \land \r C)^\perp \lor (\r A \land (\r B \land \r C)) = I=
i_1(\r A \land (\r B \land \r C), (\r A \land \r B) \land \r C).$ Thus with 
respect to the interpretation of ${\cal A} \to {\cal B}$ as  conditional  
$i_1$ the axiom schemata  for the system $S$ are theorems and the system $S$ 
is compatible with deductive quantum logic under the Mittelstaedt 
conditional.\pars
Recall that $\r d 
= \{\r F_i \mid i \in \nat \}$ is a 
development paradigm$,$ where each $\r F_i$ is a 
readable frozen segment, and describes the behavior of a Natural system at each moment of 
a time interval.  For the next construction 
a formal language that is, of course,  isomorphic to the informal 
language is employed. Each $\land$ [resp. $\r F_i$] corresponds to a specific $\rm \sp and\sp$ [resp. a propositional atom that corresponds to a specific word] when  embedded.  This eliminates confusion, when $\rm \sp and\sp$ appears in the $\r F_i$.  Let $\r M_0 = \r d.$ Define $\r M_1 = \{\r F_0\sp{\rm and}\sp\r 
F_1\}.$ Assume that $\r M_n$ is defined. Define $\r M_{n+1} = \{\r x\sp{\rm 
and}\sp\r F_{n+1}\mid \r x \in \r M_n\}.$ From the fact that $\r d$ is a 
developmental paradigm$,$ where the last two symbols in each member of $\r d$ is 
the time indicator  ``i.''$,$ it follows that no member of $\r d$ is a member 
of $\r M_n$ for $n > 0.$ Now let ${\rm M_d} = \bigcup \{\r M_n \mid n 
\in \nat\}.$ Intuitively$,$ $\sp{\rm and}\sp$ behaves as a conjunction and each 
$\r F_i$ as an atom within our language. Notice the important formal demonstration fact that for
an hypothesis consisting of any member of $\r M_n,\ n>0,$ repeated applications of
(1)$,$ (MP)$,$ (2)$,$ (MP) will lead to the members of $\r d$ appearing in the proper time ordering at
increasing (formal) demonstration step numbers.\parm
{\bf Theorem 7.3.1} {\sl For $\b d = \{\b F_i \mid i\in \nat \},$ there exists 
an ultraword $w \in \Hyper {\bf M_d} - \Hyper {\b d}$ such that $\b F_i \in 
\Hyper {\b S}(\{w\})$ {\rm (}i.e. $w \Hyper {\vdash_S} \b F_i{\rm )}$ for each $i \in 
\nat.$}\pars 
Proof. Consider the binary relation $G = \{(x,y)\mid (x \in \b d) \land (y \in 
{\bf M_d - d}) \land ( x \in \b S(\{y\}) \}.$ Suppose that $\{(x_1,y_1),\ldots 
(x_n,y_n)\}\subset G.$ For each $i=1,\ldots,n$ there is a unique $k_i \in \nat$ such 
that $x_i = \b F_{k_i}.$ Let $m = \max \{k_i\mid (x_i = \b F_{k_i})\land
(i = 1,\ldots,n) \}.$ Let $b \in \b M_{m+1}.$ It follows immediately that 
$x_i \in \b S(\{b\})$ for each $i = 1,\ldots, n$ and$,$ from the construction of
$\r d,\ b \notin \b d.$ Thus $\{(x_1,b),\ldots,(x_n,b)\} \subset G.$  
Consequently$,$ $G$ is a concurrent relation. Hence$,$ there exists some $w \in 
\Hyper {\bf M_d}- \Hyper {\b d}$ such that $^\sigma \b F_i = \b F_i \in \Hyper {\b 
S}(\{w\})$ for each $i \in \nat.$ This completes the proof. [See note 4.] \qed
Observe that $w$ in Theorem 7.3.1 has all of the formally expressible 
properties of a readable word. For example$,$ $w$ has a hyperfinite 
length$,$ among other properties. However$,$ since $\r d$ is a denumerable set$,$ 
each ultraword has a very special property. \pars
Recall that for each $[g] \in \cal E$ there exists a unique $m \in \nat$ and 
$f^\prime \in  T^m$ such that $[f^\prime] = [g]$ and for each $k$ such that $m 
< k \in \nat,$ there does not exist $g^\prime \in T^k$ such that $[g^\prime] = 
[g].$ The function $f^\prime \in T^m$ determines all of the alphabet symbols$,$ 
the symbol used for the blank space$,$ and the like$,$ and determines there 
position within the intuitive word being represented by $[g].$ Also for each 
$j$ such that
$0\leq j\leq m,\ f^\prime(j) = i(\r a) \in i[{\cal W}] = T,$ where 
$i(\r a)$ is the ``encoding'' in $T$ of the symbol ``a''. For each $m \in 
\nat,$ let $P_m = \{f \mid ( f \in T^m) \land (\exists z ((z \in {\cal E}) 
\land (f \in z)\land \forall x((x \in \nat)\land (x > m)\to \neg \exists y((y \in T^x)\land 
(y \in z)))))\}.$ An element $n \in \Hyper T$ is a {{\it subtle 
alphabet symbol}} if there exists $ m \in \nat$  and $f \in \Hyper (P_m)  = (\hyper P)(\Hyper m) = (\hyper P)_m = \hyper P_m\ (\ m = \Hyper m)$ or
if $\delta \in \nat_\infty$ and $f \in \hyper {P_\delta},$ and some $j \in \hypernat$ such that 
$f(j) = n.$ A symbol is a {{\it pure subtle alphabet symbol}} if 
$f(j) = n \notin i[{\cal W}].$ Subtle {alphabet symbols can be characterized}
in $\hyper {\cal E}$ for they are singleton objects. A $[g] \in \hyper {\cal 
E}$ represents a subtle alphabet symbol iff there exists some $f\in (\Hyper 
T)^0$ such that $[f] = [g]= [(0,f(0))],\ f = \{(0,f(0))\}.$ \parm

{\bf Theorem 7.3.2} {\sl Let $\r d = \{\r F_i \mid i\in \nat\}$ be a denumerable 
developmental paradigm and use ${\cal M}_1$ of 9.1. For each ultraword $w$, that yields $\b d$ via $\Hyper {\b S}$, as in Theorem 7.3.1, the external cardinality of the collection of all pure subtle alphabet symbols represented in each $w$ is greater than or equal to $2^{\aleph_0.}$}\pars

Proof. Consider the conceptual Kleene ``tick'' notation for the natural numbers (i.e $_{|}, _{||}. _{|||}, \ldots$).
For this proof, let $_{|}$ correspond to 0. Every member of $\r d = \{\r F_i \mid i\in \nat\}$ contains a distinct symbol-string  $\r b_i$ that represents the natural number followed by the ``period'' symbol that appears as the last two symbols in a member of $\r d.$ Consider the single $\r W_n \in \r M_n,\ n > 0.$ Then $n +2$ of these distinct symbol-strings, the Kleene symbols and a ``period'' symbol, appear in $\r W_n$ along with other alphabet systems. Hence, in $\r W_n$, there are more than $n+2$ alphabet symbols. \pars
For the embedding $\cal E,$ there is a $\r W_n$ representation $[g] \in \cal E$ and two unique mappings $f_k \sim f_0 \sim g,$ where the inverse of the embedding $i$ yields the entire word for $f_0$ and, for $f_k \in P_k,$ yields the entire word as it is join constructed from individual symbols (eq. 1.2.4). In this case, $k > n +2.$ \pars

Consider the *-transform.  Let $w = [g]$ be an ultraword such that for each $i \in \nat$, $\b F_i \in \Hyper {\b S}(\{w\}).$ Theorems 7.3.1 shows that such ultrawords exist. From the definition of $\r S$, $w \in \Hyper {\bf M_d} - \sig {\bf M_d}.$ Hence, there is a $\nu,\ \delta \in \nat_\infty$ and $\Hyper {\b M}_\nu \in \{\Hyper {\b M}_x\mid x \in \hypernat\}$ such that $[f_\delta] = [g] \in \Hyper {\b M}_\nu,\  \delta > \nu + 2,\ f_\delta \in \hyper P_\delta$ \pars

Let $\rm K= \{[1,n +2]\mid n \in \nat\}.$ Then there exists a mapping $\rm C\colon K \to \nat$ such that $\rm C([1,n+2]) = n+2.$ The mapping $\rm C$ is considered as yielding the intuitive cardinality of $[1,n+2].$ Hence, $\Hyper {\b C}([1,\nu +1]) = \nu +2.$ 
To get an idea as to the external cardinality $\vert [1,\nu +2]\vert$ of $[1,\nu +2]$, consider Theorem 3.1 in [16, p. 201], where it is shown that $\vert [1,\nu +2] \vert \geq 2^{\aleph_0}.$ Since $\vert [1,\delta]\vert \geq \vert [1, \nu +2]\vert,$ and the set of all subtle alphabet symbols that yields members of $\cal W$ is denumerable, then it follows that for $w$ the set of all pure subtle alphabet symbols also has an external cardinality great than or equal to $2^{\aleph_0}.$ This completes the proof. \qed

With respect to the proof of Theorem 7.3.2$,$ the function $f_\delta$ determines the 
alphabet composition of the ultraword $w.$ The word $w$ is unreadable not only 
due to its infinite length but also due to the fact that it is composed of 
infinitely many purely subtle alphabet symbols.\pars
The developmental paradigm d utilized for the two previous theorems is 
composed entirely of readable sentences. We now investigate what happens if a 
developmental paradigm contains countably many unreadable sentences. Let the nonempty developmental paradigm $d^\prime$ be composed of at most 
countably many members of $\hyper {\cal E} - \cal E$ and, for countable $\r B$, let $d^\prime \subset 
\Hyper {\b B} \subset \Hyper {\b P_0}.$ Construct$,$ as previously$,$ the set 
${\rm M_B}$ from $\r B,$ rather than from d and suppose that $\r B \cap 
\r M_i = \emptyset,\ i \not= 0.$ [See Note [2] on page 82.] Let $\empty \not= \lambda \subset \nat.$\parm

{\bf Theorem 7.3.3} {\sl Let $d^\prime = \{[g_i]\mid i\in \lambda  \}.$ Then 
there exists an ultraword $w\in \Hyper {\bf M_B} - \Hyper {\b B}$ such that for 
each $ i\in \nat,\ [g_i] \in \Hyper {\b S}(\{w\}).$}\pars
Proof. Consider the internal binary  relation $G = \{(x,y)\mid (x \in \Hyper 
{\b B})\land (y \in \Hyper {\bf M_B} - \Hyper {\b B})\land (x \in \Hyper {\b 
S}(\{y\})\}.$ Note that members of $d^\prime$ are members of $^\sigma {\cal E}$ 
or$,$ at the most$,$ denumerably many members of $\hyper {\cal E} -{^\sigma{\cal 
E}}.$ From the analysis in the proof of Theorem 7.3.1$,$ for a finite $\r F \subset \r B,$ 
there exists some ${\rm y \in M_B- B}$ such that $\r F 
\subset S(\{\r y\}).$ It follows by *-transfer
that if $F$ is a finite or *-finite subset of $\Hyper {\b B},$ then there exists 
some $y \in \Hyper {\bf M_B} - \Hyper {\b B}$ such that $F \subset \Hyper {\b 
S}(\{y\}).$ As in the proof of Theorem 7.3.1$,$ this yields that $G$ is at least 
concurrent on $\Hyper {\b B}.$ However$,$ $d^\prime \subset \Hyper {\b B}$ and  
$\vert d^\prime \vert \leq \aleph_0.$ From $\aleph_1$-saturation$,$ there exists 
some $w \in\Hyper {\bf M_B} -\Hyper {\b B}$ such that for each $[g_i] \in 
d^\prime,\ [g_i] \in \Hyper {\b S}(\{w\}).$ This completes the proof. \qed

Let $\emptyset \not= \lambda, \gamma \subset \nat$, $j \in \gamma$, ${\cal D}_j = \{d_{ij} \mid i \in \lambda \}$, and for each $j \in \gamma$, $i \in \lambda,\ d_{ij} \subset \Hyper {\b B}$ 
is considered to be a developmental paradigm either of type $d$ or type 
$d^\prime$ and $\r B \cap \r M_i = \emptyset,\ i \not=0.$ Notice that ${\cal D}_j$ may be either a finite
or denumerable set and 
Theorem 7.3.1 holds for the case that $\r d \subset \r B,$ where $w \in 
\Hyper {\bf M_B} -\Hyper {\b B}.$ 
For each $d_{ij} \in {\cal D}_j,$ use the Axiom of Choice 
to select an ultraword $w_{ij} \in \Hyper {\bf M_B} - \Hyper {\b B}$ that exists 
by Theorems 7.3.1 (extended) or 7.3.3. Let $\{w_{ij} \mid i \in \lambda\}$ be
ultrawords. \parm
\vfil\eject

{\bf Theorem 7.3.4} {\sl For $j \in \gamma,$ there exists an ultimate ultraword $w_j' \in \Hyper {\bf 
M_B} - \Hyper {\b B}$ such that for each $i \in \lambda,\ w_{ij} \in \Hyper {\b 
S}(\{w_j' \})$ and$,$ hence$,$ for each $d_{ij} \in {\cal D}_j,\ d_{ij} \subset \Hyper {\b S}(\{w_{ij}\})\subset \Hyper {\b S}(\{w_j'\}).$ }\pars
Proof. For each finite ${\rm \{F_1,\ldots,F_n\} \subset {M_B} - {B}}$ 
there is a natural number$,$ say $m,$ such that for each $i = 1,\ldots, n,
\ \r F_i \in \r M_j$ for some $j \leq m.$ Hence$,$ taking $\r b \in \r M_{m+1},$ we 
obtain that each $\r F_i \in S(\{\r b\}).$ Observe that $\r b \notin {\r B.}$ 
By *-transfer$,$ it follows that the 
internal relation $G =\{(x,y)\mid (x \in \Hyper {\bf M_B} -\Hyper 
{\b B})\land (y \in \Hyper {\bf M_B} - \Hyper {\b B})\land (x \in \Hyper {\b 
S}(\{y\})\}$ is concurrent on internal $\Hyper {\bf M_B} -\Hyper 
{\b B}$ and $\{w_{ij}\mid i \in \lambda \} \subset \Hyper {\bf M_B} -\Hyper 
{\b B}.$ Again $\aleph_1$-saturation yields that there is some $w_j' \in 
\Hyper {\bf M_B} -\Hyper {\b B}$ such that for each $i \in \lambda,\ w_{ij} \in 
\Hyper {\b S}(\{w_j'\}).$ The last property is obtained from $d_{ij} \subset 
\Hyper {\b S}(\{w_{ij}\}) \subset \Hyper {\b S}(\Hyper {\b S}(\{w_j^\prime\}))= 
\Hyper {\b S}(\{w_j^\prime\})$ since $\{w_{ij}\}$ is an internal subset of $\Hyper
{\b P_0}.$ This completes the proof.\qed 

{\bf Corollary 7.3.4.1} {\sl There exists an ultimate ultraword $w' \in \Hyper {\bf M_B} - \Hyper {\b B}$ such that for each $j \in \gamma,\ w_j' \in \Hyper {\b 
S}(\{w' \})$ and$,$ hence$,$ for each $d_{ij} \in \bigcup{\cal D}_j,\ d_{ij} \subset \Hyper {\b S}(\{w_j'\})\subset \Hyper {\b S}(\{w'\}).$ }\parm

The same analysis used to obtain Theorem 7.3.2  
can be applied to the 
ultrawords of Theorems 7.3.3 and 7.3.4. (See note [5].)\par
\bigskip
\leftline{\bf 7.4 Ultracontinuous Deduction}
\medskip
In 1968$,$ a special {topology on the set of all nonempty  subsets} of a given set 
$X$ was constructed and investigated by your author. We apply a similar 
topology to subsets of $\cal E.$ \pars
Suppose that nonempty $X \subset \cal E.$ Let $\tau$ be the discrete topology 
on $X.$ In order to topologize $\power X,$ proceed as follows: for each $G \in 
\tau,$ let $N(G) = \{A\mid (A \subset X)\land (A \subset G)\}=\power G.$ 
Consider ${\cal B} = \{N(G)\mid G \in \tau \}$ to be a base for a topology $\tau_1$ on $\power X.$ Let $A \in N(G_1) 
\cap N(G_1).$ The discrete topology implies that $N(A)$ is a base element and 
that $N(A) \subset N(G_1) \cap N(G_2).$ 
There is only one 
member of $\cal B$ that contains $X$ and this is $\power X.$ Thus if $\power 
X$ is covered by members of $\cal B,$  then $N(X)=\power X$ is one of these covering 
objects. Thus $(\power X, \tau_1)$ is a compact space. Further$,$ since 
$N(\emptyset) \subset N(G)$ for each $G \in \tau,$ the space  $(\power X, 
\tau_1)$ is connected. The topology $\tau_1$ is a special case of a more 
general topology with the same properties. [2] Suppose that $D \subset X.$ Let 
$D \in N(G)= \power G, G \in \tau.$ Then $D \in N(D) \subset N(G).$ This 
yields that the nonstandard monad is $\monad D = \bigcap \{\hyper N(G) \mid N(G) \in {\cal B}\}= \Hyper 
{(\power D)} = \Hyper {\power {\hyper D}}.$ \parm
{\bf Theorem 7.4.1} {\sl Any consequence operator $C \colon (\power X, 
\tau_1) \to (\power X,\tau_1)$ is continuous.}\pars
Proof. Let $A \in \power X$ and $H \in \Hyper {\b C}[\monad A].$ Then there 
exists some $B \in \monad A$ such that $\Hyper {\b C}(B) = H.$ Hence$,$ $B \in 
\Hyper {\power {\hyper A}}.$ By *-transfer of a basic property of our 
consequence operators$,$ $\Hyper {\b C}(B) \subset \Hyper {\b C}(\hyper A) = 
\Hyper {(\b C(A))}.$ Thus $\Hyper {(\b C(B))} \in \Hyper {(\power {\b 
C(A))}}$ implies that $\Hyper {\b C}(B) \in \monad A.$ Therefore$,$ $\Hyper {\b 
C}[\monad A] \subset \monad {[\b C(A)]}.$ Consequently$,$ $\b C$ is continuous. \qed
{\bf Corollary 7.4.1.1} {\sl For any $X \subset \cal E$$,$ and any consequence 
operator $\b C \colon \power X \to \power X,$ the map $\Hyper {\b C}\colon 
\Hyper {(\power X)} \to  \Hyper {(\power X)}$ is ultracontinuous.}\parm
{\bf Corollary 7.4.1.2} {\sl Let $\r d$ {\r [}resp. $d^\prime$$,$ $\r d$ or 
$d^\prime${\r ]} be a developmental 
paradigm as defined for Theorem 7.3.1 {\r [}resp. Theorem 7.3.3$,$ 7.3.4{\r ]}. Let 
$w$ be a 
ultraword that exists by Theorem 7.3.1 {\r [}resp Theorem 7.3.3$,$ 7.3.4{\r ]}. Then 
$\r d$ {\r [}resp. $d^\prime$$,$ $\r d$ or $d^\prime${\r ]} is obtained by means of a 
ultracontinuous subtle deductive process applied to $\{w\}.$ }\parm
Recall that in the real valued case$,$ a function $f \colon [a,b] \to \real$ is 
uniformly continuous on $[a,b]$ iff for each $p,\ q \in \Hyper {[a,b]}$ such 
that $ p -q \in \monad 0,$ then $f(p) - f(q) \in \monad 0.$ If $D \subset [a,b]$ 
is compact$,$ then $p,\ q \in \Hyper D$ and $p-q \in \monad 0$ imply that there is 
a standard $ r \in D$ such that $p, \ q \in \monad r.$ Also$,$ for each $ r\in D$ 
and any $p,\ q \in \monad r,$ it follows that $p-q \in  \monad r.$ Thus$,$ if 
compact $D \subset [a,b],$ then $f:D \to \real$ is uniformly continuous iff 
for every $r \in D$ and each $p,\ q \in \monad r,\ \hyper f(p),\ \hyper f(q) 
\in \monad {f(r)}.$ With this characterization in mind$,$ it is clear that any 
consequence operator $\b C \colon \power X \to \power X$ satisfies the 
following statement. For each $A \in \power X$ and each $ p,\ q \in \monad A,\ 
\Hyper {\b C(p)},\ \Hyper {\b C(q)} \in \monad {\b C(A)}.$\pars
From the above discussion$,$ one can think of ultracontinuity as being a type 
of {ultrauniform continuity}.\par
\bigskip
\leftline{\bf 7.5 Hypercontinuous Gluing}
\medskip 
There are various methods that can be used to investigate the behavior of 
adjacent frozen segments. All of these methods depend upon a significant 
result relative to discrete real or vector valued functions. The major goal in 
this section is to present a complete proof of this major result and to 
indicate how it is applied. \pars
First$,$ as our standard structure$,$ consider either the intuitive real numbers 
as atoms or axiomatically a standard structure with atoms ${\bf ZFR} = {\bf ZF} 
+ {\bf AC}  +A_1({\rm atoms}) +A{\rm (atoms)} + \vert  A \vert = \r c,$ where $A$ 
is isomorphic to the real 
numbers and $A_1 \cap A = \emptyset.$ 
Then$,$ as done previously$,$ there is a model 
$\langle C,\in,= \rangle $ within our ${\bf ZF} + {\bf AC}$ 
model for  ${\bf ZFR},$ where $A$ has all of the 
ordered field properties as the real numbers. A superstructure $\langle {\cal 
R}, \in, = \rangle $ is constructed in the usual manner$,$ where the  superstructure 
$\langle {\cal N}, \in, = \rangle $ is   a substructure. Proceeding as in 
Chapter 2$,$ construct $\hyper {\cal M}_1 = \langle \Hyper {\cal R}, \in = 
\rangle$ and ${\cal Y}_1.$ The structure ${\cal Y}_1$ is called\break
 the {{\it
Extended Grundlegend Structure}} --- the EGS. The Grundlegend 
Structure is a substructure of ${\cal Y}_1.$\pars
It is important to realized in  what follows that the objects utilized for the 
G-structure {\it interpretations} are nonempty finite equivalence classes of 
partial sequences. Due to this fact$,$ the following results should not lead to 
ambiguous interpretations. \pars
As a preliminary to the technical aspects of this final section$,$ we 
introduce the following 
definition. A function $f\colon [a,b] \to \realp m$ is {{\it 
differentiable-C}} on $[a,b] $ if it is continuously differentiable on $ 
(a,b)$ except at finitely many removable discontinuities. This definition is 
extended to the end points $\{a,b\}$ by application of one-sided 
derivatives.  For any $[a,b],$ consider a partition $ P = 
\{a_0,a_1,\cdots,a_n,a_{n+1}\},\ n\geq 1,\ a = a_0,\ b = a_{n+1} $ and
$ a_{j-1} < a_j,\ 1\leq j\leq n+1.$ For any such
partition $ P,$ let the real valued function $ g $ be
defined on the set $ D = [a_0,a_1)\cup  (a_1,a_2) \cup \cdots \cup 
(a_n,a_{n+1}] $ as follows: for each  $x\in [a_0,a_1),$ let $ g(x) = 
r_1 \in \real$; for each $ x \in (a_{j-1},a_j),$ let $ g(x) = r_j \in 
\real ,\ 1 < i \leq n$; for each $ x \in (a_n,b],$ let $ g(x) = r_{n+1} 
\in \real .$ 
It is obvious that $ g $ is a type of simple step 
function. Notationally$,$ let
$ {\cal F}(A,B) $ denote the set of all functions with domain $ A $
and codomain $ B.$  \parm
\indent {\bf Theorem 7.5.1} {\sl There exists a function 
$ G\in \Hyper ({\cal F}([a,b],\real))  $ with the following properties.\pars
\indent\indent (i) The function $ G $ is *-continuously *-differentiable and *-uniformly
*-continuous on $ \hyper {\eskip [a,b]},$\hfil\break
\indent\indent (ii) for each odd $ n\in \hypernat ,\ (n \geq 3),\ G $ is *-
differentiable-C of order $ n $ on $ \hyper {\eskip [a,b]},$\hfil\break
\indent\indent (iii) for each even $ n\in \hypernat ,\  G $ is *-continuously 
*-differentiable in $\hyper {[a,b]} $ except at finitely many 
points,\hfil\break
\indent\indent (iv) if $ c = \min \{r_1,\cdots,r_{n+1}\},\ d=\max 
\{r_1,\cdots,r_{n+1}\} $$,$ then the range of $ G = \hyper {\eskip [c,d]},\ \st 
G $
at least maps $ D $ into $ [c,d] $ and $ (\st G)\vert D = g.$}\pars
   Proof. First$,$ for any real $ c,d,$ where $ d \not= 0,$ consider the 
finite set of functions $$h_j(x,c,d) = (1/2)(r_{j+1}-r_j)\Bigl(\sin \bigl((x 
-c)\pi/(2d)\bigr) +1\Bigr) + r_j,\eqno (7.5.1)$$
$1\leq j \leq n.$ Each $ h_j $ is continuously differentiable for any order 
at each $ x \in \real .$ Observe that for each odd $ m\in \nat,$ each
m'th derivative $ h_j^{(m)} $ is continuous at $ (c+d) $ and $ (c-d) $ 
and $ h_j^{(m)}(c+d) = h_j^{(m)}(c-d) = 0 $ for each $ j.$ \pars
Let positive $ \delta \in \monad 0.$ Consider the finite set of internal 
intervals $ \{[a_0,a_1 - \delta),(a_1 + \delta , a_2 - \delta),\cdots,
(a_n + \delta,b]\} $ obtained from the partition $ P.$ Denote these intervals 
in the expressed order by $ I_j,\ 1\leq j\leq n+1.$ Define the internal 
function $$G_1 =\{(x,r_1)\vert x\in I_1 \} \cup \cdots \cup \{(x,r_{n+1})\vert 
x\in I_{n+1} \}.\eqno (7.5.2)$$
Let internal $ I_j^\dagger = [a_j - \delta, a_j + \delta],\ 1 \leq j \leq n,$ 
and for each $ x\in I_j^\dagger,$ let internal      
$$G_j(x) = (1/2)(r_{j+1}-r_j)\Bigl(\hyper {\sin}
 \bigl((x-a_j)\pi/(2\delta)\bigr)+1\Bigr) + r_j.\eqno (7.5.3)$$
Define the internal function
$$G_2 =\{(x,G_1(x))\vert x\in I_1^\dagger \}\cup \cdots \cup \{(x,G_n(x))\vert 
x\in I_n^\dagger \}.\eqno (7.5.4)$$
The final step is to define $ G = G_1 \cup G_2$. Then $ G\in \Hyper ({\cal 
F}([a,b],\real)).$\pars
By *-transfer$,$ the function $ G_1 $ has an internal *-continuous *-
derivative $ G_1^{(1)} $ such that $ G_1^{(1)}(x) = 0 $ for each 
$ x\in I_1 \cup\cdots\cup I_{n+1}.$ Applying *-transfer to the properties of 
the functions $ h_j(x,c,d),$ it follows that $ G_2 $ has a unique internal 
*-continuous *-derivative $$G_2^{(1)} =(1/(4\delta))(r_{j+1} - r_j)\pi
\Bigl(\hyper \cos \bigl( (x - a_j)\pi/(2\delta)\bigr)\Bigr)\eqno (7.5.5)$$
for each $ x\in I_1^\dagger \cup\cdots\cup I_n^\dagger.$  The results that 
the *-left limit for the internal $ G_1^{(1)} $ and the *-right limit for 
internal $ G_2^{(1)} $ at each $ a_j - \delta $ as well as the *-left 
limit of $ G_2^{(1)} $ and *-right limit of $ G_1^{(1)} $ at each $ a_j 
+\delta $ are equal to $ 0 $ and $ 0 = G_2^{(1)}(a_j - \delta) = 
G_2^{(1)}(a_j+\delta) $ imply that internal $ G $ has a *-continuous *-
derivative $ G^{(1)} = G_1^{(1)} \cup G_2^{(1)} $ defined on $ \hyper 
{\eskip [a,b]}.$\pars
A similar analysis and *-transfer yield that for each $ m \in \hypernat,\ m 
\geq 2,\ G $ has an internal *-continuous *-derivative $ G^{(m)} $ defined 
at each $ x\in \hyper {\eskip [a,b]} $ except at the points $ a_j\pm \delta $ 
whenever $ r_{j+1} \not= r_j.$ However$,$ it is obvious from the definition of 
the functions $ h_j $ that for each odd $ m\in \hypernat, \ m\geq 3,$ each 
internal $ G^{(m)} $ can be made *-continuous at each $ a_j\pm \delta $ by 
simply defining $ G^{(m)}(a_j\pm \delta) = 0 $ and with this parts (i)$,$ 
(ii)$,$ and (iii) are established.\pars
For part (iv)$,$ assume that $ r_j \leq r_{j+1}.$ From the definition of the 
functions $ h_j,$ it follows that for each $ x \in I_j\cup I_j^\dagger\cup
I_{j+1},\ r_j \leq G(x) \leq r_{j+1}.$ The nonstandard intermediate value 
theorem implies that $ G\bigl[ \hyper {\eskip [a_j,a_{j+1}]}\bigr] = \hyper 
{\eskip [r_j,r_{j+1}]} $ and in like manner for the case that $ r_{j+1}< r_j.$
Hence$,$ $ G\bigl[\hyper {\eskip [a,b]}\bigr] = \hyper {\eskip [c,d]}.$ Clearly$,$ $\st {\Hyper {D}} =
[a,b].$ If $ p\in D $ and $ x\in \monad p \cap \Hyper {D},$ then $ G(x) = r_j=g(p) $
for some $ j $ such that $ 1 \leq j \leq n+1.$ This completes the 
proof. \qed
The nonstandard approximation theorem 7.5.1 can be extended easily to functions 
that map $ D $ into $\realp m.$ For example$,$ assume that $ F\colon D \to 
\realp 3,$  the component functions $ F_1,\  F_2 $ are continuously 
differentiable on $ [a,b];$ but that $ F_3 $ is a $ g $ type step 
function on $ D.$ Then letting $ H = (\hyper F_1,\hyper F_2,G),$
on $ \hyper {\eskip [a,b]},$ where $ G $ is defined in Theorem 4.1$,$ we 
have an internal *-continuously *-differentiable function $ H\colon \hyper {\eskip [a,b]} 
\to \hyperrealp 3,$ with the property that $ \st H\vert D = F.$\pars
With respect to Theorem 7.5.1$,$ it is interesting to note that if $h_j$ is 
defined on $\real,$ then for even orders $ n \in \nat,$ 
$$\vert\  h_j^{(n)}(c \pm d)\ \vert = \bigg|{{(r_{j+1} -r_j)\pi^n}\over {2^{n+1}d^n}}
\biggl| = 0 \leqno (7.5.6)$$
for $ r_{j+1} = r_j$ but not 0 otherwise. If $r_{j+1} - r_j \not= 0,$ then 
$G_2^{(n)}(a_j \pm \delta)$ is an infinite nonstandard real number. Indeed$,$ if 
$m_i$ is an increasing sequence of even numbers in $\hypernat$  and $r_{j+1} 
\not= r_j,$ then $\vert G_2^{(m_i)}(a_j \pm \delta) \vert$ forms a decreasing 
sequence of nonstandard infinite numbers. The next result is obvious from the 
previous result. \parm
{\bf Corollary 7.5.1.1} {\sl For each $ n \in \hypernat,$ then internal 
$G^{(n)} = G_1^{(n)} \cup G_2^{(n)}$ is *-bounded on $\Hyper [a,b].$}\parm
Let $D(a,b)$ be the set of all bounded and piecewise continuously 
differentiable functions defined on $[a,b].$ By considering 
all of the possible (finitely many) subintervals$,$ where $f \in D(a,b),$ it 
follows from the Riemann sum approach that for each real $\nu > 0,$ there 
exists a real $\nu_1 > 0$ such that for each real $\nu_i,\ 0<\nu_i < \nu_1,$ 
a sequence of partitions $P_i = \{a = b_0^i < \cdots < b_{k_i}^i = b \}$ can 
be selected such that the ${\rm mesh}(P_i) \leq \nu_i$ and 
$$ \vert (f(b) - f(a)) - \sum_{n=1}^{k_i} f^\prime(t_n)(b_n^i - b_{n-1}^i) 
\vert < \nu \leqno (7.5.7)$$ for any $t_n \in (b_{n-1}^i,b_n^i).$ \pars
Moreover$,$ for any given number $M,$ the sequence of partitions can be so 
constructed such that there exists a $j$ such that for each $i > j,\ k_i >M,$ 
where $P_i$ and $P_j$ are partitions within the sequence of partitions. By 
*-transfer of these facts and by application of Theorem 7.5.1 and its 
corollary we have the next result.\parm
{\bf Corollary 7.5.1.2} {\sl For each $n \in \nat$ and each internal 
$G^{(n)},$ 
the difference $G^{(n)}(b) - G^{(n)}(a)$ is infinitesimally close to an 
(externally) infinity *-finite sum of infinitesimals.}\parm
A developmental paradigm is a very general object and$,$ therefore$,$  can be 
used for numerous applications. At present$,$ developmental paradigms are still 
being viewed from the {{\it substratum}} or external world. For what follows$,$ 
it is assumed that a developmental paradigm d traces the evolutionary history 
of a specifically named natural system or systems. In this first application$,$ 
let each $\r F_i \in \r d$ have the following property ({\b P}).\pars
 {\leftskip=0.5in \rightskip= 0.5in \noindent  $\r F_i$ describes 
``the general behavior and characteristics of the named natural system $S_1$ 
as well as the behavior and characteristics of named constituents contained 
within $S_1$ at time $t_i.$''\par}\pars
Recall that for $\r F_i,\ \r F_{i+1} \in \r d,$ there exist unique functions 
$f_0 \in \b F_i = [f],\ g_0 \in \b F_{i+1} = [g]$ such that $f_0,\ g_0 \in 
T^0$ and $\{(0,f_0(0))\} \in [f],\ \{(0,g_0(0))\} \in [g].$ Thus$,$ to each 
$\r F_j \in \r d,$ correspond the unique natural number $f_0(0).$ Let $D = 
[t_{i-1}, t_i) \cup (t_i,t_{i+1}]$ and define $f_1\colon D \to \nat$ as 
follows: for each $x \in  [t_{i-1}, t_i),$ let $f_1(x) = f_0(0);$ for each $ x 
\in (t_i,t_{i+1}],$ let $f_1(x) = g_0(0).$ Application of theorem 7.5.1 yields 
the internal function $G$ such that $G\vert D = f_1.$ For these physical 
applications$,$ utilize the term ``substratum'' in the place of the technical 
terms  ``pure nonstandard.'' [Note: Of course$,$ elsewhere$,$ the term  
``pure NSP-world'' or simply the ``{NSP-world}'' is used as a specific name for 
what has here been declared as the substratum.] This yields the following 
statements$,$ where the symbols $\r F_i$ and $\r F_{i+1}$ are defined and 
characterized by the expression inside the quotation marks in property (\b 
P).\parm
{\leftskip=0.5in \rightskip=0.5in \noindent (\b A): There exists a substratum 
hypercontinuous$,$ hypersmooth$,$ hyperuniform process $G$ that 
binds together $\r F_i$ and $\r F_{i+1}.$\par}\pars
{\leftskip=0.5in \rightskip=0.5in \noindent (\b B): There exists a substratum 
hypercontinuous$,$ hypersmooth$,$ hyperuniform alteration process $G$ that 
transforms $\r F_i$ into $\r F_{i+1}.$\par}\pars
{\leftskip=0.5in \rightskip=0.5in \noindent (\b C): There exists an 
ultracontinuous  subtle force-like (i.e. deductive) process that yields 
$\r F_i$ for each time $t_i$ within the development of the natural 
system.\par}\parm
In order to justify (\b A) and (\b B)$,$ specific measures of physical properties 
associated with constituents may be coupled together. Assume that for a subword 
$\r r_i \in F_i \in \r d,$ the symbols $\r r_i$ denote a numerical quantity 
that aids in characterizing the behavior of an object in a system $S_1$ or the 
system itself.  Let $(\b M_1)$ be the statement:\parm
{\leftskip=0.5in \rightskip=0.5in \noindent  ``There exists a substratum  
hypercontinuous$,$ hypersmooth$,$ hyperuniform functional process $G_i$ such that 
$G_i$ when restricted to the standard mathematical domain it is $f_i$ and such 
that $G_i$ hypercontinuously changes $r_i$ for system $S_1$ at time $t_i$ into 
$r_{i+1}$ for system $S_1$ at time $t_{i+1}.$''\par}\parm
\noindent This modeling procedure yields the following interpretation:\parm
{\leftskip=0.5in \rightskip=0.5in \noindent (\b D) If there exists a continuous or 
uniform [resp. discrete] functional process $f_i$ that changes $r_i$ for $S_1$
at time $t_i$ into $r_{i+1}$ for $S_1$ at time $t_{i+1},$ then $(\b 
M_1).$\par}\parm
At a particular moment $t_i,$ two natural systems $S_1$ and $S_2$ may 
interface. More generally$,$ two very distinct developmental paradigms may  
exist one $\r d_1$ at times prior to $t_i$ (in the $t_i$ past) and one $\r d_2$ 
at time after $t_i$ (in the $t_i$ future). We might refer to the time $t_i$ as 
a {{\it standard time fracture.}} Consider the developmental paradigm 
$\r d_3 = \r d_1 \cup \r d_2.$ In this case$,$ the paradigms may be either of type
$\r d$ or $d^\prime.$ For the type $d^\prime,$ the corresponding system need 
not be considered a natural system but could be a pure substratum system.\pars
At $t_i$ an $\r F_i \in  \r d_3$ can be characterized by statement (\b P) 
(with the term natural removed if $\r F_i$ is a member of a $d^\prime$).  
In like manner$,$ $\r F_{i+1}$ at time $t_{i+1}$ can be characterized by (\b 
P). Statements (\b A)$,$ (\b B)$,$ (\b C) can now be applied to $\r d_3$ and a 
modified statement (\b D)$,$ where the second symbol string $S_1$ is changed to 
$S_2.$ Notice that this modeling applies to the actual human ability that only 
allows for two discrete descriptions to be given$,$ one for the interval 
$[t_{i-1},t_i)$ and one for the interval $(t_i,t_{i+1}].$ From the modeling 
viewpoint$,$ this is often sufficient since the length of the time intervals can 
be made smaller than Planck time. \pars
Recall that an analysis of the scientific method used in the investigation of 
natural system should take place exterior to the language used to describe the 
specific system development. Suppose that $\cal D$ is the language accepted 
for a scientific discipline and that within $\cal D$ various expressions from 
mathematical theories are used. Further$,$ suppose that enough of the modern 
theory of sets is employed so that the EGS can be constructed. The following 
statement would hold true for $\cal D.$ \parm
{\leftskip=0.5in \rightskip=0.5in \noindent {\it
If by application of first-order logic to a set of non-mathematical premises 
taken from $\cal D$ it is claimed that it is not logically possible for 
statements such as \r (\b A\r )$,$ \r (\b B\r )$,$ \r (\b C\r ) and \r (\b D\r ) to hold$,$ then the set of 
premises is inconsistent.}\par}\parm
\vfil\eject
\centerline{\bf CHAPTER 7 REFERENCES}
\medskip
\noindent {\bf  1} Beltrametti$,$ E. G.and G. Cassinelli$,$ The logic of quantum 
mechanics$,$ in {\it Encyclopedia of Mathematics and its Application}$,$ 
Vol. 15$,$ Addison-Wesley$,$ Reading$,$ 1981.\pars
\noindent {\bf 2} Herrmann$,$ R. A.$,$ {\it Some Characteristics of Topologies on 
Subsets of a Power Set}$,$ University Microfilm$,$ M-1469$,$ 1968.\pars 
\noindent {\bf 3} Herrmann$,$ R. A.$,$ {\it The Theory of Infinitesimal Light-clocks,} (1993)
http://arxiv.org/abs/math/0312189 \parm
\centerline{NOTES}\pars
 [1] Note that this last requirement for $\r B$ can be achieved as follows: construct a special symbol not originally in $\cal A$. Then this symbol along with $\cal A$ is considered the alphabet. Next only consider a $\r B$ that does not contain this special symbol within any of its members. Now using this special symbol in place of the $\land,$ consistently construct $\r M_i,\ i \not= 0.$
 Of course$,$ $\land$ is interpreted as this special symbol in the axiom system $S.$\pars
[2] The actual members, $\r F_i,$  of a developmental paradigm $\rm d$ need not be unique. However, the specific information contained in each readable word used for a specific $\r F_i \in \nat$  is unique. Other readable sentences can be used in place of a specific $\r F_i$ as long as they are ``equivalent'' in the sense that the specific information being displayed by each is the same information.\pars
[3] Depending upon the application, a single standard word may also be termed as an ultraword. \pars
[4] For a method to obtain an ultraword for a refined developmental paradigm, see pages 4 - 7 in http://arxiv.org/abs/math.0605120 \pars
[5] (Added 9/20/2009.) A concurrent relation is not needed to obtain important ``ultrawords.''\pars
{\bf Theorem 7.3.5} {\sl For $\b d = \{\b F_n\mid n \in \nat\}$ and each infinite $\lambda \in \nat_\infty$, there exists one and only one $w_\lambda \in \Hyper {\b M}_\lambda$ and hyperfinite $d_\lambda$ such that $\b d \subset d_\lambda \subset \Hyper {\b S}(\{w_\lambda\}),$ and $d_\lambda \subset \hyper {\b d}.$}\pars
Proof. For each, $n \in \nat$, let $G(n) = \{\b F_i\mid 0\leq i \leq n\}\subset \b d.$ Thus,
$G\colon \nat \to {\cal F}(\b d)$ the set of all finite subsets of $\b d.$ Let $n > 0.$ Then $\b M_n$ has one and only one member and by definition $\b w_n \in \b M_n$ has the property that $G(n) \subset \b S(\{\b w_n\}).$ Hence, by *-transfer, for the function $\Hyper {\b M}$, and each $\lambda \in \nat_\infty,$ there is one and only one $w_\lambda \in \Hyper {\b M}_\lambda$ such that hyperfinite $\Hyper {G}(\lambda) \subset \Hyper {\b S}(\{w_\lambda\}).$ Finally, by defintion of $G$, $\b d \subset \Hyper G(\lambda) \subset \hyper {\b d}.$ \qed
Note that theorems that generate or use ultrawords may need to be trivially modified or not used depending upon the definition for d. For example, for the ordering used in [1], then the ultraword used and its location would be an ultimate ultraword as generated in theorem 7.3.4, where each $w_{ij}$ is an ultraword or ultimate ultraword that, upon application of $\Hyper {\b S}$, yields the developmental paradigm for each interval $[c_i, c_{i+1}).$\parm
\noindent [1] The GGU-model and Generation of the Developmental Paradigms\hfil\break http://arxiv.org/abs/math/0605120  
\vfil\eject

\centerline{\bf 8. A SPECIAL APPLICATION}\par
\bigskip
\leftline{\bf 8.1 A Neutron Altering Process.}\par 
\medskip
The purpose of this chapter is to justify the interpretations utilized in 
reference [1]. Let $B$ be the set of all nondecreasing bounded real valued
functions defined on $D = [a,t^\prime) \cup (t^\prime, T].$ Let ${\cal Q} \in B
$ and ${\cal Q}(t) =  2$ for each $t \in 
[a,t^\prime);\ {\cal Q}(t) = 3$ for each $ t \in (t^\prime,T]$ be the {{\it 
discrete  neutron altering process}}. Application of Theorem 7.5.1 implies 
that there exists internal $G\colon \Hyper [a,T] \to \hyperreal$ such that 
$\st G\vert D = G\vert D={\cal Q},$ and $G$ is hypercontinuous$,$ hypersmooth$,$ hyperaltering 
process defined on the hyperinterval $\Hyper [a,T].$ Hence$,$ $G$  satisfies 
statement (A) in section (2) of [1]. Theorem 7.5.1 also implies that $G$ is 
hyperuniformly continuous on $\Hyper [a,T].$ \pars
Recall how a {*-special partition} for $\Hyper [a,T]$ is generated. Let $0< 
\Delta t \in \real^+.$ Then $P(\Delta t) = \{a = t_0 \leq \cdots t_n\leq t_{n+1} = 
T\},$ where $n$ is the largest natural number such that $a + n(\Delta t) \leq 
T$ and for $i = 0,\ldots, (n -1);\ t_{i+1} - t_i = \Delta t,$ and $t_{n+1} - 
t_n = b - (a + n(\Delta t)) < \Delta t.$ It is possible that $t_{n+1} = t_n.$ 
If $\cal P$ is the set of all special partitions$,$ then letting $dt \in \monad 
0^+$ (the set of all positive infinitesimals) it follows that 
$P(dt) \in \hyper {\cal P}$ and $P(dt)$ has the same first-order properties as does 
$P(\Delta t).$ \parm
{\bf Theorem 8.1.1} {\sl Let internal $G$ be hypercontinuous on $\Hyper C = 
\Hyper [a,T]$ and $z_j \in (a,T),\ j = 1,\ldots,m.$ For any $dx \in \monad 
0^+,$ there exists a $dy \in \monad 0^+$ such that for each $x,\ y \in \hyper 
D$ such that $\vert x-y \vert < dy $ it follows that $\vert G(x) -G(y) \vert < 
dx,$ and there is a hyperfinite partition $\{a= t_0^\prime < \cdots < 
t_{\nu+1}^\prime = T\}$ such that for $ i=0,\ldots, \nu+1;\ j = 1,\ldots,m$ we 
have $t_j^\prime \not= z_j,\ G(t_{i+1}^\prime) - G(t_i^\prime) \in \monad 0,\ 
t_{i+1}^\prime - t_i^\prime \in \monad 0$ and $G(T) -G(a) = \sum_0^\nu 
(G(t_{i+1}^\prime) - G(t_i^\prime)).$}\pars
Proof. Since internal $G$ is *-uniformly continuous$,$ it follows that for any 
$dy \in \monad 0^+$ there exists some $\delta$ such that $0 < \delta \in 
\hyperreal$ and for each $x,\ y \in  \Hyper [a,T]$ such that $\vert x - y 
\vert < \delta,$ it follows that $\vert G(x) - G(y) \vert < dx.$ Now let $dy < 
\delta$ and $dy \in \monad 0^+$ and consider the *-special partition 
$P(dy/3).$ Let $y \in [t_i,t_{i+1}],\ x \in [t_{i+1},t_{i+2}],\ i= 
0,\ldots,\nu-1$ and $x,\ y \not= t^\prime.$ Then $\vert y-x \vert < dy$ and 
each *-closed interval is nonempty. By means of internal first-order 
statements that imply the existence of certain objects and the choice axiom$,$ 
select $t_0^\prime = a, \ t_{\nu+1}^\prime = T$ and if $t_{\nu+1} = 
t_\nu,$ then for $i = 1,\ldots, \nu -2$ select some $t_i^\prime \in 
[t_i,t_{i+1}]$ such that $t_i^\prime \not= z_j$ for $j = 1,\ldots,m;$ or if 
$t_{\nu+1} \not= t_\nu,$ then for $i=1,\ldots,\nu-1$ select some $t_i^\prime
\in [t_i,t_{i+1}]$ such that $t_i^\prime \not= z_j$ for $j = 1,\ldots,m.$ This 
yields a hyperfinite internal partition with the properties listed in the 
hypothesis and by\noindent  *-transfer of the properties of a finite telescoping series$,$ 
we have that $G(T) - G(a) = \sum_0^\nu((G(t_{i+1}^\prime) - G(t_i^\prime)),$ 
and $\vert G(t_{i+1}^\prime) - G(t_i^\prime)\vert < dx$ for $i = 0,\ldots,\nu$ implies 
that $G(t_{i+1}^\prime) - G(t_i^\prime) \in \monad 0$ and each $t_i^\prime \in 
\Hyper D$ has the property that $\vert t_{i+1}^\prime - t_i^\prime \vert < 
dy.$ This complete the proof. \qed
 We now apply Theorem 8.1.1 to the discrete altering function $\cal Q.$ Let 
$Q$ be the set of all finite partitions of $D.$ Then$,$ for $n > 0$ and the 
partition $\{ a = t_0 <t_1 < \cdots < t_n \leq t_{n+1}= T \},$ consider the partial 
sequence $S\colon [0,n+1] \to \real $ defined by $S(i) = t_i,\ i = 0, \ldots, 
n+1.$ Define $T_i = [t_i, t_{i+1}],\ i = 0, \ldots, n.$ Consider the set\break
\centerline{}
\vfil
\eject
 \noindent$H = 
\{T_i \mid i = 0, \ldots, n\}.$ Then $H\in \power C,$ where $C = [a, T].$ There 
is an $N \in \real^{\power C - \emptyset}$ such that $N(T_i) = {\cal 
Q}(t_{i+1)} - {\cal Q}(t_i), \ {\cal Q}(t_i) = r_i,\ i=0,\ldots, n.$ The 
function $N$  is a {{\it resolving process}} for the function $\cal Q$ and 
each $r_i$ is a degree for the {{\it constituent}} $N(T_i).$ Let $M(\cal Q)$
be the set of all such resolving processes generated by the infinite set of 
finite partitions of $D$ for a fixed $\cal Q.$ Consider the *-finite partition 
$P(dy/3)$ of $D$ generated in the proof of Theorem 8.1.1. Now modify this 
*-finite partition in the following manner. Consider the standard finite 
partition generated by $S\colon [0,n+1] \to \real.$ Let $T_i = [t_i, t_{i+1}]= 
[S(i),S(i+1)], \ t_i \not= t^\prime,\ i = 0, \ldots, n;\ H = \{T_i \mid 
i=0,\ldots, n\}$ and$,$ for the fixed ${\cal Q},\ N(t_i)= {\cal Q}(t_{i+1}) - 
{\cal Q}(t_i) = G(t_{i+1}) - G(t_i),\ i =0, \ldots, n,$ where $G$ is in 
the statement of Theorem 8.1.1. This sequence $S$ extends in the usual manner 
to $\Hyper S\colon [0, \nu +1] \to \hyperreal.$ \pars
Since $0 < \vert t_i - t_j \vert,$ where $i \not= j,$ for each $t_i$ there 
exists a *-closed interval $[v_j, v_{j+1}]$ generated by $P(dy/3)$ such that 
$t_i \in (v_j, v_{j+1}),$ or $t_i = v_j$ or $v_{j+1}$ not both. In the case 
that $t_i \in (v_j,v_{j+1}),$ the interval is unique. Moreover$,$ there are only 
finitely many such $t_i$ in the standard partition. Hence for these finitely 
many real number cases$,$ where $t_i \in (v_j, v_{j+1}),$ modify the partition 
by subdividing $[v_j,v_{j+1}]$ into two intervals $[v_j,t_i] \cup 
[t_i,v_{j+1}].$ This process can$,$ obviously$,$ be defined by a finite set of 
first-order statements. This adds an additional finite number of intervals to 
our hyperfinite partition and yields a partition number $\lambda \in 
\hypernat$ to replace $\nu.$ Since the infinitesimal length of these adjoined 
intervals is $< dy/3,$ Theorem 8.1.1 still holds with $\lambda$ replacing 
$\nu.$ This yields an internal sequence $S^\prime \colon [0,\lambda +1] \to
\hyperreal$ such that $S^\prime(i) = t^\prime_i$ as defined in the proof of 
Theorem 8.1.1 with the condition that finitely many of the $t^\prime_j$ 
correspond to the standard partition elements $t_i.$ (Notice that$,$ for all of 
this construction$,$ the assumed $n$ is fixed.) \pars
From  the above$,$ we have the *-closed  intervals $T^\prime_i = 
[t^\prime_i,t^\prime_{i+1}],\ i = 0,\ldots,\lambda$ as well as the internal
$H^\prime = \{z\mid (z \in \Hyper{\power { \Hyper {[a,T]}}}\land (\exists i(i \in 
[0,\lambda]) \land \forall x(x \in z \iff S^\prime(i) \leq x \leq 
S^\prime(i+1)))\}.$ Obviously$,$ $H^\prime \in \Hyper (\power C).$ 
Thus there is in $\hyper M({\cal Q})$ an internal hyperresolving process 
$N^\prime$ such that $N^\prime(T^\prime_i) = \hyper {\cal Q}(t^\prime_{i+1}) - 
\hyper {\cal Q}(t^\prime_{i+1}) = G(t^\prime_{i+1}) - G(t_i),$ where $T^\prime 
= [S^\prime(i),S^\prime(i+1)] \in H^\prime$ and $i = 0,\ldots,\lambda.$ \pars
Technically$,$ it is not true that $H \subset H^\prime.$ Thus define the 
{{\it standard restriction}} of $N^\prime$ to $N,$ where $N$ is generated by 
the standard sequence $S\colon [0,n+1] \to \real$ that is obtained as follows: 
consider the set $\{S^\prime(i) \mid i=0,\ldots,\lambda +1 \} \cap [a,T] = 
P_0.$ Since $P_0$ is a finite standard set$,$ it can be ordered by the $<$ of 
the reals and let $P_0 = \{a = t_0 < t_i < \cdots t_n \leq t_{n+1} = T.\}$ 
This yields a sequence $S"\colon [0,n+1] \to \real,\ S"(i) = t_i,\ 
i=0,\ldots, n+1.$ Let $S" = S.$ Utilizing $S",$ generate the original 
resolving process $N$ from $N^\prime.$\pars
Application of Theorem 8.1.1 yields the following description. There exists a 
hyperpartition (generated by) $S^\prime$ for the hyperinterval $\hyper D$ 
(since $t^\prime \in {\rm range}\ S^\prime$) and $S^\prime$ (generates) the 
hyperresolution $N^\prime$ for the hyperaltering process $G.$ The 
hyperresolution $N^\prime$ is defined on the hyperfinitely many internal 
subintervals of $\Hyper [a,T]$ and the range of $N^\prime$ is composed of 
hyperfinitely many hyperconstitutents $G(t^\prime_{i+1}) - G(t^\prime_i)$ 
that$,$ by *-transfer of the standard supremum function defined on nonempty 
finite sets of real numbers$,$ yields a maximum degree among all of the degrees 
of the hyperconstitutents. This maximum degree is infinitesimal$,$ by Theorem 
8.1.1$,$ and since $G(T) - G(a) \in\real^+$ and taking $G$ as nondecreasing$,$ this maximum 
degree is a positive infinitesimal. By the above restriction process$,$ $N$ is 
the restriction of $N^\prime$ to the standard world.  Consequently$,$ $N^\prime$ 
satisfies statement (B) in section 4 of [1]. \pars
Finally$,$ the length function $L$ defined on the set of all closed intervals
extends to the set of all *-closed intervals that are subsets of $\hyperreal.$
Then $\hyper L(\Hyper [a,T]) = T - a = L([a,T]).$ Thus (C) of section 4 in [1]
holds. (D) in section 4 of [1] follows from the unused conclusions that appear 
in Theorem 8.1.1$,$ among others.\pars
For the nondecreasing bounded classical neutron altering process $\cal CQ,$ 
there is assumed to exist a standard smooth function $f$ defined on $[a,t]$ 
such that $f| D = {\cal CQ}.$ Now define standard $G\colon [a,T] \to \real$ 
as follows: let $G_0(t) = f(t),\ t \in [a,t^\prime);\ G_1(t) = f(t),\ t \in 
(t^\prime,T].$ Then since $G_0(t_0) \leq G_1(t_1),$ for $t_0 \in [a,t^\prime)$ 
and $t_1 \in (t^\prime,T],$ it follows that $h = \sup \{G_0(t)\mid t\in 
[a,t^\prime)\}$ exists$,$ and we can let $G = G_0 \cup G_1 \cup
\{(t^\prime,h)\}.$ Obviously$,$ $G(t_0) \leq G(t^\prime) \leq G(t_1),\ t_0 \in 
[a,t^\prime),\ t_1 \in (t^\prime,T],$ and $G| D = {\cal CQ}.$ It follows from 
left and right limit considerations that $G=f.$ (Note: G is defined in this 
manner only to conform to the discrete case.) Theorem 8.1.1 holds for $\hyper 
G$ and$,$ in this case$,$ we simply repeat the entire discussion that appears after 
that statement of Theorem  8.1.1 and replace the $G$ that appears in that discussion 
with $\hyper G = \hyper f.$ This yields a model for statements (E)$,$ (F)$,$ (G) and  
(H) in section 4 of reference [1].\parm
\centerline{\bf CHAPTER 8 REFERENCE}
\medskip
\noindent {\bf 1} Herrmann$,$ R. A. Mathematical philosophy and developmental 
processes$,$ {\it Nature and System,} 5(1/2) (1983)$,$ 17---36.\par                                                         
 \vfil
\eject
\centerline{NOTES}                                                                                          
\vfil\eject
 \centerline{\bf 9. NSP-WORLD ALPHABETS}\par
\bigskip
\leftline{\bf 9.1 An Extension.}\par 
\medskip
Although it is often not necessary$,$ we assume when its useful that we are 
working 
within the {EGS}. Further$,$ this structure is assumed to be ${\vert {\cal 
M}_1\vert}^+$-saturated$,$ and a polyenlargement [5, p. 35], where ${\cal M}_1 =\langle {\cal R}, \in,= \rangle,$ (or ${\cal M}_1 =\langle {\cal Q}, \in,= \rangle,$ where $\cal Q$ is  the set of rational numbers). 
Referring to the paragraph prior to Theorem 7.3.3$,$ it can be assumed 
that the developmental paradigm $d^\prime \subset \Hyper {\b B} \subset \Hyper
{\b P_0.}$ It is not assumed that such a developmental paradigm is obtained 
from the process discussed in Theorem 7.2.1$,$ although a modification of the 
proof of Theorem 7.2.1 appears possible in order to allow this method of 
selection. \parm
{\bf Theorem 9.1.1} {\sl Let $d^\prime = \{[g_i] \mid i \in \lambda \},\ \vert 
\lambda \vert < {\vert {\cal M}_1 \vert}^+.$ There exists an ultraword $w \in 
\Hyper {\b M_B} - \Hyper {\b B}$ such that for each $ i \in \lambda,\ 
[g_i] \in \Hyper {\b S}(\{w\}).$}\pars
Proof. The same as Theorem 7.3.3 with the change in saturation.\qed
Let ${\cal D} = \{d_i \mid i \in \lambda\}, \ \vert \lambda \vert< {\vert 
{\cal M}_1\vert}^+,\ \vert d_i \vert < {\vert {\cal M}_1\vert}^+$ and each $d_i 
\subset \Hyper {\b B}$ 
is considered to be a developmental paradigm either of type $d$ or type 
$d^\prime.$ 
For each $d_i \in \cal D,$ use the Axiom of Choice 
to select an ultraword $w_i \in \Hyper {\b M_B} - \Hyper {\b B}$ that exists 
by Theorems 9.1.1. Let $\{w_i \mid i \in \lambda\}$ be such 
a set of ultrawords. \parm
{\bf Theorem 9.1.2} {\sl There exists an ultraword $w^\prime \in \Hyper {\b 
M_B} - \Hyper {\b B}$ such that for each $i \in \lambda,\ w_i \in \Hyper {\b 
S}(\{w^\prime \})$ and$,$ hence$,$ for each $d_i \in {\cal D},\ d_i \subset \Hyper {\b 
S}(\{w^\prime\}).$ }\pars
Proof. The same as Theorem 7.3.4 with the change in saturation.\qed
\smallskip
\leftline{\bf 9.2 NSP-World Alphabets.}
\medskip
First$,$ recall the following definition.
$P_m = \{f \mid ( f \in T^m) \land (\exists z ((z \in {\cal E}) 
\land (f \in z)\land \forall x((x \in \nat)\land (x > m)\to \neg \exists y((y \in T^x)\land 
(y \in z)))))\}.$ The set $T=i[{\cal W}].$ The set $P_m$ determines the 
unique partial sequence $f \in [g] \in {\cal E}$ that yields$,$ for 
each $ j \in \nat$ such that $0\leq j \leq m,\ f(j) = i(\r a),$ where 
$i(\r a)$ is an ``encoding'' in $A_1$ of the alphabet symbol  ``a'' used to 
construct our intuitive language $\cal W.$ The set $[g]$ represents an 
intuitive word constructed from such an {alphabet of symbols.} \pars
Within the discipline of Mathematical Logic$,$ it is assumed
that there exists 
symbols --- a sequence of variables --- each one of which corresponds$,$ in a one-to-one 
manner$,$ to a natural number. Further$,$ under the subject matter of {generalized 
first-order theories} [2]$,$ it is\noindent 
assumed that the cardinality of the set of 
constants is greater than $\aleph_0.$ In the forthcoming investigation$,$ it may 
be useful to consider an alphabet that injectively corresponds to the real 
numbers $\real.$ This yields a new alphabet ${\cal A}^\prime$ containing 
our original 
alphabet. A new collection of words ${\cal W}^\prime$ composed of 
nonempty finite strings of such alphabet symbols may be constructed. It may also be 
useful to well-order $\real.$ The set $\cal E$ also exists with respect to the 
set of words ${\cal W}^\prime.$ Using the {ESG}$,$ many  
 previous results in this book now hold with respect to 
${\cal W}^\prime$ and for the case that we are working in a ${\vert{\cal 
M}_1\vert}^+$-saturated polyenlargement.\pars 
 With respect to this {extended language}$,$ if you 
wish to except the possibility,
a definition as to what constitutes a {purely 
subtle alphabet symbol} would need to be altered in the obvious fashion. 
Indeed$,$ for $T$ in the definition of $P_m,$ we need to substitute $T^\prime = i[{\cal 
W}^\prime].$ Then the altered definition would read that $r \in \hyper A_1 
\simeq \hyperreal$ is 
a {{\it pure subtle alphabet symbol}} if there exists an $m \in \nat$ and $f 
\in \Hyper {(P_m)},$ or if $m \in \hypernat - \nat$ an $f \in P_m,$ and some 
$j \in \hyperreal$ such that $f(j) = r \notin i[{\cal W}^\prime].$  Further$,$ 
some of the previous theorems also hold when the proofs are modified.\pars
Although these extended languages are of interest to the mathematician$,$ most 
of science is content with approximating a real number by means of a rational 
number. In all that follows$,$ the cardinality of our language$,$ if not 
denumerable$,$ will  be 
specified. All theorems from this book that are used to establish a result 
relative to a denumerable language will be stated without qualification. If 
a theorem has not been reestablished for a higher language but can be so 
reestablished$,$ then the theorem will be termed an {{\it extended}} theorem.
\pars
\bigskip
\leftline{\bf 9.3 General Paradigms.}
\medskip
There is the developmental paradigm$,$ and for nondetailed descriptions the 
{{\it general developmental paradigm.}} But now we have something totally new 
--- the {{\it general paradigm}}. It is important to note that the general 
paradigm is considered to be distinct from developmental paradigms$,$ although 
certain results that hold for general paradigms will hold for developmental 
paradigms and conversely. For example$,$ associated with each general paradigm 
${\rm G_A}$
is an {ultraword} $w_g$ such that the set ${\bf G_A} \subset  
\Hyper {\b S}(\{w_g\})$ and all other theorems relative to such 
ultrawords hold for general paradigms. The {general paradigm} is a collection 
of words that discuss$,$ in general$,$ the behavior of entities and other 
constituents of a natural system. They$,$ usually$,$ do not contain a time 
statement $\r W_i$ as it appears in section 7.1 for developmental paradigm 
descriptions. Our interest in this section is relative to only two such general 
paradigms. The reader can easily generate many other general paradigms.  
\pars
Let $\r c^\prime$ be a symbol that denotes some fixed real number and $ 
n^\prime $ a symbol that denotes a natural number. [Note: what follows is 
easily extended to an {extended language}.] Suppose that you have a theory 
which includes each member of the following set ($i$ suppressed).\parm
\line{\hfil ${\rm G_A} = \{{\rm An\sp elementary\sp particle\sp 
\alpha(\r n^\prime)\sp with\sp}$\hfil}\pars
\line{(9.3.1)\hfil  ${\rm kinetic\sp energy\sp \r c^\prime}{+}1/(\r n^\prime).\mid 
n \in G \land n \not= 0\},$\hfil }\parm
\noindent where $G$ is, at the least, a denumerable subset of the real numbers. 
\pars
Of particular interest is the composition of members of $\Hyper {\bf G_A - 
G_A}.$  Notice that $\vert {\bf G_A} \vert = \vert G\vert$ 
since  $z_1,\ z_2 \in {\bf G_A}$  and $z_1 \not= z_2$ iff $[x_1] = z_1,\ 
[x_2] = z_2,\ x_1(30)= x_1(2) \not= x_2(30)= x_2(2),\ x_1(2),\ x_2(2) \in 
G.$ Now consider the bijection 
$K\colon {\bf G_A} \to G.$ \parm
\vfil
\eject
{\bf Theorem 9.3.1} {\sl The set $[g] \in \Hyper {\bf G_A} - {\bf G_A}$ iff
there exists a unique  $f\in \Hyper (P_{55})$ and $\nu \in \Hyper G-G$ such that $[g] = [f],$ and 
$f(55) = i(\r A),\ f(54) = i(\r n),\ f(53) = i(\sp), \cdots, f(30) = f(2), 
\cdots, f(3) = i((), f(2)=\nu \in \Hyper G - G \subset \hyperreal - \real,\ 
f(1) 
= i()),\ f(0) = (.).$ }\pars
Proof. From the definition of $\rm G_A$ the sentence \pars
\line{\hfil $ \forall z(z \in {\cal E}\to (( z \in {\bf G_A}) \iff \exists ! x\exists ! w(
(w \in G) \land (x \in P_{55}) 
\land (x \in z) \land $\hfil}\pars
\line{\hfil $((55,i(\r A)) \in x) \land ((54,i(\r n)) \in x) \land \cdots 
\land (x(30) = x(2)) \land 
\cdots \land $ \hfil}\pars
\line{\hfil $((3,i(()) \in x)\land (x(2)=w)\land(K(z) = w) \land$\hfil}\pars 
\line{\hfil\qquad $((1,i())) \in x) \land ((0,i(.)) \in x)))).$\hfil (9.3.2)}\parm
\noindent holds in $\cal M,$ hence in $\Hyper {\cal M}.$ From the fact that 
$K$ is a bijection$,$ it follows that $\hyper K[\Hyper {\bf G_A} - {\bf G_A}] = 
\Hyper G -G \subset \hyperreal - \real.$ The result now follows from 
*-transfer. \qed
Using Theorem 9.3.1$,$
each member of $\Hyper 
{\bf G_A} - {\bf G_A},$ when interpreted by considering $i^{-1},$ has only two 
positions with a single missing object since positions 30 and 2 do not 
correspond to any symbol string in our language ${\cal W}.$ 
This 
interpretation still retains a {vast amount of content}$,$ however. For a specific 
member$,$ you could substitute a {new constructed symbol}$,$ not in ${\cal W}',$ 
into these 
two missing positions. Depending upon what type of pure nonstandard number 
this inserted symbol represents$,$ the content of such a sentence could be 
startling. Let $\Gamma^\prime$ be a nonempty set of new symbols disjoint from ${\cal 
W}'$ and assume that $\Gamma^\prime$ is injectively mapped by $H$ into 
$\Hyper G - G.$ \pars
Although human ability may preclude the actual construction of more 
than denumerably many new symbols$,$ you might consider this mapping to be 
onto if you accept the ideas of {extended languages} with a greater cardinality.  
As previously$,$ denote these new symbols by $\zeta^\prime.$  Now let \par
\line{\hfil ${\rm G_A^\prime} = \{{\rm An\sp elementary\sp particle\sp 
\alpha(\zeta\prime)\sp with\sp}$\hfil}\pars
\line{\hfil \qquad\qquad${\rm kinetic\sp energy\sp \r c^\prime}{+}1/(\zeta^\prime).\mid 
H(\zeta^\prime) \in \Hyper G - G\},$ \hfil (9.3.3)}\parm
\noindent This leads to the following interpretation stated in terms of describing 
sets for the extended language.\parm
{\leftskip=0.5in \rightskip=0.5in \noindent (1) {\it The describing set $\rm G_A$ (mathematically) exists 
iff the describing set $\rm G_A^\prime$ (mathematically) exists.}\par}\par
\bigskip
\leftline{\bf 9.4 Interpretations}
\medskip
Recall that the Natural world portion of the NSP-world model may contain 
{{\it undetectable}} objects$,$ where  ``undetectable'' means that there does 
not appear to exist  human$,$ or humanly constructible machine sensors that 
directly detect the objects or directly measure any of the objects physical 
properties. The rules of the scientific method utilized within the 
micro-world of subatomic physics allow all such 
undetectable Natural objects to be accepted as existing reality.[1] The 
properties of such objects are indirectly deduced from the observed properties 
of gross matter. In order to have indirect evidence of the objectively real 
existence of such objects$,$ such indirectly obtained behavior will usually satisfy a 
specifically accepted model.\par
Although the numerical quantities associated with these undetectable Natural 
(i.e. standard) world objects$,$ if they really do exist$,$ cannot be directly and exactly
measured via any known instrumentation, these quantities are still represented by standard mathematical 
entities. By the rules of correspondence for interpreting pure NSP-world 
entities$,$ the members of ${\rm G_A^\prime}$ must be considered as {undetectable
pure NSP-world objects$,$} assuming any of them exist in this background 
world. On the other hand$,$ when viewed within the EGS$,$ any finite as well as 
many infinite subsets of $\rm 
G_A^\prime$ are internal sets. 
Consequently$,$ some finite collects of such objects {\it may be}  assumed to {indirectly 
effect behavior in the Natural world.} \par
The concept of {{\it realism}} often dictates that all interpreted members of 
a mathematical model be considered as existing in reality. The philosophy of 
science that accepts only {{\it partial realism}} allows for the following technique. 
One can stop at any point within a mathematically generated physical 
interpretation. 
Then proceed from that point to deduce an intuitive physical theory$,$ but only 
using other not interpreted mathematical formalism as {auxiliary constructs} or 
as {catalysts}.  With respect to the 
NSP-world$,$ another aspect of interpretation enters the picture. Assuming 
realism$,$ then the question remains which$,$ if any$,$ of these NSP-entities 
actually indirectly influence Natural world processes?  This interpretation 
process allows for the possibility that none of these pure NSP-world 
entities has 
any effect upon the standard world. These ideas should always be kept in 
mind.\parm
If you accept that such particles as described by $\rm G_A$ can exist in 
reality$,$ then the {philosophy of realism} leads to the next interpretation. 
\parm
{\leftskip 0.5in \rightskip 0.5in \noindent (2) {\it If there exist elementary 
particles with Natural system behavior described by $\rm G_A$$,$  then there 
exist pure NSP-world objects that display within the NSP-world behavior 
described by members of $\rm G_A^\prime.$}\par}\parm
\noindent The concept of {absolute realism} would require that the 
acceptance of the elementary 
particles described by $\rm G_A$ is indirect evidence for the existence 
of the $\rm G_A^\prime$ described objects. I caution the reader that the 
interpretation we apply to such sets of sentences as $\rm G_A$ are only to be 
applied to such sets of sentences.\pars
The EGS may$,$ of course$,$ be interpreted in infinitely many different ways. 
Indeed$,$ the {NSP-world model with its physical-type language} can also be 
applied in infinitely many ways to infinitely many scenarios. I have applied it 
to such models as the {MA-model| and the GGU-model among others. In this section$,$ I consider 
another possible interpretation relative to those {Big Bang cosmologies} that 
postulate real 
objects at or near infinite temperature$,$ energy or pressure. These theories 
incorporate the concept of the {{\it initial singularity(ies).}} \pars
One of the great difficulties with many Big Bang cosmologies is that no 
meaningful physical interpretation for formation of the initial singularity is 
forthcoming from the theory itself. The fact that a proper and acceptable 
theory for {creation of the universe} requires that consideration not only be 
given to the moment of {zero cosmic time} but to what might have occurred 
``prior'' to that moment in the nontime period is what partially influenced 
 {Wheeler} to consider 
the concept of a  {{\it pregeometry.}}[3]$,$ [4] It is totally unsatisfactory 
to dismiss such questions as ``unmeaningful'' simply because they cannot be 
discussed in your favorite theory. Scientists must search for a {broader 
theory} 
to include not only the question but a possible answer.\pars
Although the {initial singularity} for a Big Bang type of state of affairs 
apparently cannot be discussed in a meaningful manner by many standard 
physical theories$,$ unless one adjoins to the theory an ad hoc quantum field, it can be discussed by application of our NSP-world 
language. Let $\r c^\prime$  be a symbol that represents any fixed real 
number. Define \parm
\line{\hfil ${\rm G_B} = \{{\rm An\sp elementary\sp particle\sp 
\alpha(\r n^\prime)\sp with\sp}$\hfil}\pars
\line{ \hfil ${\rm total\sp energy\sp \r c^\prime}{+}\r n^\prime.\mid 
n \in \nat\},$\qquad \hfil (9.4.1)}\parm
Application of Theorem 9.3.1 to $\rm G_B$ yields the set \parm
\line{\hfil ${\rm G_B^\prime} = \{{\rm An\sp elementary\sp particle\sp 
\alpha(\zeta^\prime)\sp with\sp}$\hfil}\pars
\line{\hfil ${\rm total\sp energy\sp \r c^\prime}{+}\zeta^\prime.\mid 
\zeta \in \hypernat-\nat\},$ \hfil (9.4.2)}\parm
{\leftskip 0.5in \rightskip 0.5in \noindent (3) {\it If there exist elementary 
particles with Natural system behavior described by $\rm G_B$$,$  then there 
exist pure NSP-world objects that display within the NSP-world behavior 
described by members of $\rm G_B^\prime.$}\par}\parm
The particles being described by $\rm G_B^\prime$ have various {infinite 
energies}. These infinite energies {\bf do not} behave in the same manner as 
would the real number energy measures discussed in ${\rm G_B}.$ As is usual 
when a metalanguage physical theory is generated from a formalism$,$ we can 
further extend and investigate the properties of the $\rm G_B^\prime$ objects 
by imposing upon them the corresponding behavior of the {positive infinite 
hyperreal numbers}. This produces some interesting propositions. Hence$,$ we are 
able to use a nonstandard physical world language in order to give further 
insight into the state of affairs at or near a cosmic initial singularity. 
This gives {\it one} solution to a portion of the pregeometry problem. I point 
out that there are other NSP-world models for the beginnings of our universe$,$ 
if there was such a beginning. Of course$,$ the statements in $\rm G_B^\prime$  
need not be related at all to any Natural world physical scenario$,$ but could 
refer only to the behavior of pure NSP-world objects. \pars
Notice that Theorems such as 7.3.1 and 7.3.4 relative to the generation of 
developmental paradigms by ultrawords$,$ also apply to general paradigms$,$ where 
$\rm M, M_B,  P_0$ are defined appropriately.   
The 
following is a slight extension of Theorem 7.3.2 for general paradigms. 
Theorem 9.4.1 will also hold for developmental paradigms.\parm
{\bf Theorem 9.4.1} {\sl Let $\rm G_C$ be any denumerable general paradigm. 
Then there exists an ultraword  $w \in \Hyper {\b P_0}$ such that for each $\b 
F \in {\bf G_C},\ \b F \in \Hyper {\b S}(\{w\})$ and there exist infinitely 
many $[g] \in \Hyper {\bf G_C} - {\bf G_C} $ such that $[g] \in \Hyper {\b 
S}(\{w\}).$} \pars
Proof. In the proof of Theorem 7.3.2$,$ it is shown that there exists some $\nu 
\in \hypernat - \nat $ such that $\hyper h[[0,\nu]] \subset \Hyper {\b 
S}(\{w\})$  and $\hyper h[[0,\nu]] \subset \Hyper {\bf G_C}.$ Since 
$\vert \hyper h[[0,\nu]] \vert {\vert {\cal M}_1\vert}^+$, then $\vert \hyper h[[0,\nu]] 
-h[\nat] \vert \geq {\vert {\cal M}_1\vert}^+,$ for $h$ is a bijection. This completes the 
proof. \qed
{\bf Corollary 9.4.1.1} {\sl Theorem 9.4.1 holds$,$ where $\rm G_C$ is replaced 
by a developmental paradigm.}\parm
{\leftskip 0.5in \rightskip 0.5in \noindent (4) {\it Let $\rm G_C$ be a 
denumerable general paradigm. There exists an intrinsic ultranatural 
process$,$ $\Hyper {\b S},$ such that objects described by members of 
$\rm G_C$ are produced by $\Hyper {\b S}.$ 
During this production$,$ numerously many pure NSP-objects as described by 
statements in $\Hyper {\bf G_C} - {\bf G_C}$ are produced.}\par}\parm
\smallskip
\leftline{\bf 9.5 A Barrier To Knowledge.}       
\medskip
Our final discussion in this chapter deals with the use of ${\vert {\cal 
M}\vert}^+$-saturated models and our ability to analyze sets of sentences such 
as $\rm G_A^\prime.$ It is a very {strange property of the human mind} that it 
often produces an infinite cause and effect sequence.\pars
Consider the following comprehensible  and potentially infinite set of 
sentences $\{$I think.$,$ I think about my thinking.$,$ I think about my 
thinking about my thinking.$,$ I think about my thinking about my 
thinking about my thinking.$,\ldots\}.$ As far as comprehension is concerned$,$ 
one can ask what would be the  ``first cause'' for our thinking? For the 
Natural sciences$,$ we have {Engel's biological sequence of evolutionary causes
and effects} and$,$ of course$,$ our previously mentioned initial singularity problem 
or as Misner$,$ Thorne and Wheeler write$,$ ``No problem of cosmology digs more 
deeply into the foundations of physics than the question of what `proceeded' 
the `initial' state....''[3] A Natural science question and one which 
contains some logical difficulties might be  ``What precedes that which 
precedes?'' Regardless of whether or not this strange mental behavior persists 
when we analyze Natural system behavior$,$ the next result shows the existence 
of a possible Natural barrier to human knowledge.\pars
Each of our previous investigations is done with respect to a specific 
NSP-world structure $\Hyper {\cal M}$ based upon a infinite standard set $H$ 
(with a cardinality usually equal to $\aleph_0$) into which is 
mapped the symbols and words for all languages. The requirement that $H$ be a 
standard set is relative to the standard universe in which we function. 
Although there are infinitely 
many distinct nonisomorphic NSP-world structures$,$ each of our results is with 
respect to members of a subclass of the class of all such structures. In 
particular$,$ ${\vert {\cal M}\vert}^+$-saturated polyenlargement $,$ where $\cal M$ is 
based upon a standard set $H,$ where $\nat \subset H.$\pars

In order to analyze general paradigms $\rm G_A^\prime, \ G_B^\prime$ and the 
like$,$ we need to start$,$ I believe$,$ with a comprehensible set of sentences$,$ 
such as $\rm G_A,\ G_B,$ with nonempty content and insert new symbols but 
retain some of the content of the original sentences. What is shown next is 
that if we use any of our models based on $H$ and require them to be ${\vert {\cal M}\vert}^+$-saturation 
polyenlargement, then we cannot  
embed our new alphabet into the standard set $H$ and$,$ thus$,$ we 
cannot fully analyze sets of sentences such as $\rm G_A^\prime,\ G_B^\prime$ using 
our embedding procedures.  
\parm\vfil\eject
{\bf Theorem 9.5.1} {\sl Let $\Gamma^\prime$ be a set of symbols adjoined  to a 
countable alphabet ${\cal A},$ which is disjoint from $\cal A,$ and such that it is 
used to obtain 
the set of sentences in $\rm 
G_B^\prime.$ Let $\Hyper {\cal M} = \langle \Hyper {\cal H},\in,= \rangle$ be 
any $|{\cal M}|^+$-saturated polyenlargement of a superstructure 
based on the ground set $H,$ where here $\nat \subset H \subset \real.$ There does 
does not exist an injection from $\Gamma^\prime \cup {\cal A}$ into $\b L,$ where $\b L \in {\cal H}.$}\pars
Proof.   Suppose that there 
exists an injection $i\colon (\Gamma^\prime \cup {\cal A})\to \b L.$ 
Since the model is a polyenlargement, then $|\Gamma^\prime \cup {\cal A}| \geq |{\cal M}|^+.$ However, $|\b L| < |{\cal M}|^+.$ But under the assumption $|\Gamma^\prime \cup {\cal A}| \leq |\b L|.$ This contradiction implies that the
injection does not exist and this completes the proof. \qed
\bigskip
\centerline{\bf CHAPTER 9 REFERENCES}
\medskip
\noindent {\bf 1} Evans$,$ R.$,$ {\it The Atomic Nucleus}$,$ McGraw-Hill$,$ New York$,$ 
1955.\pars
\noindent {\bf 2} Mendelson$,$ E.$,$ {\it Introduction to Mathematical Logic}$,$ Ed. 
2$,$ D. Van Nostrand Co.$,$ New York$,$ 1979.\pars
\noindent {\bf 3} Misner$,$ Thorne and Wheeler$,$ {\it Gravitation}$,$ W. H. 
Freeman$,$ San Francisco$,$ 1973.\pars                                                         
\noindent {\bf 4} Patton C. and A. Wheeler$,$ Is physics legislated by 
cosmogony? in {\it Quantum Gravity}$,$ ed. Isham$,$ Penrose and Sciama$,$ Oxford 
University Press$,$ Oxford$,$ 1977$,$ 538---605.\pars
\noindent {\bf 5} Stroyan, K. D. and J. M. Bayod, {\it Foundations of Infinitesimal Stochastic Analysis}, North Holland, New York, 1986.\pars
\vfil
\eject
\centerline{NOTES}                 
\vfil\eject
\centerline{\bf 10. LAWS$,$ RULES AND OTHER THINGS}\par
\bigskip
\leftline{\bf 10.1 More About Ultrawords.}\par 
\medskip
Previously$,$ we slightly investigated the composition of an ultraword
$w \in  \Hyper {\bf M_d} - \b d.$ Using the idea of the minimum informal language 
$\r P_0\subset \r P,$  where $\r d$ is denumerable and $\r P$ is a propositional 
language$,$ our interest now lies in completely determining 
the composition of $\Hyper {\b S}(\{w\}).$ [Note: since our language is 
informal axiom (3) and (4) are redundant in that superfluous parentheses have 
been removed.] First$,$ two defined sets.\parm
\line{\hfil A $=\{\r x \mid \r x \in \r P_0$ is an instance of an axiom for 
S$\}$\hfil (10.1.1)}\pars  
\line{\hfil C $= \{\r x \mid \r x \in \r P_0$ is a finite $(\geq 1)$ conjunction of 
members of d$\}$\hfil (10.1.2)}\parm 
Notice that it is also possible to refine the set C by considering C to be an 
ordered conjunction with respect to the ordering of the indexing set used to 
index members of d. Further$,$ as usual$,$ we have that ${\rm A,\ C,\ d}$ are mutually disjoint.\parm
{\bf Theorem 10.1.1} { \sl Let $w \in \Hyper {\bf M_d} - \Hyper {\b d}$ 
be an ultraword for infinite $\b d \subset \Hyper {\b S}(\{w\}).$ Then $\Hyper 
{\b S}(\{w\}) = \Hyper {\b A} \cup Q_1 \cup d_1^\prime,$ where for internal *-
finite $d_1^\prime,\ \b d \subset d_1^\prime \subset \Hyper {\b d}$ and 
internal $Q_1 \subset \Hyper {\b C}$ is composed of *-finite $(\geq 1)$ 
conjunctions (i.e. $i({\rm\sp and\sp})$) of distinct members of $d_1^\prime$ and $w \in Q_1.$ Further$,$ each member of 
$d_1^\prime$ and no other *-proposition  is used to form the *-finite 
conjunctions in $Q_1,$ the only *-propositions in $\Hyper {\b S(\{w\})}$ are those in 
$w,$ and $\Hyper {\b A},\ Q_1$ and $d_1^\prime$ are mutually 
disjoint.}\pars 
Proof. The intent is to show that if ${\rm w \in M_d - d},$ then $\rm 
S(\{ w \})=A \cup Q \cup d^\prime,$ where $\rm Q \subset C,$ finite $ \rm 
d^\prime \subset d$ and $\rm Q$ is composed of finite $(\geq 1)$ conjunctions 
of members of $\rm d^\prime,$  each member of $\rm d^\prime$ is used to form 
these conjunctions and no other propositions.\pars
Let $\r J$ be the set of propositional atoms in the composite w. (\b 0) Then $\rm J 
\subset 
S(\{w\}).$ If K is the set of all propositional atoms in $\rm S(\{w\}),$ then 
$\rm J \subset K.$ Let $\rm b \in K - J.$ It is obvious that $ \rm b \notin 
\rm S(\{ w \})$ since otherwise $\rm \{w,b\} \subset S_0 (\{ w \})$ but
$\rm \not\models_{S_0} w \to b.$  
Thus$,$ $\rm J = K.$  Consequently$,$ $\rm J \subset 
S(\{w\}),\ J \subset d$ and there does not exists an $\rm F \in d - J,$ such 
that $\rm F \in S(\{w\}).$ (\b 1) {\it Let $\rm J = d^\prime.$ The only 
propositional atoms in                                               
 $\rm S(\{w\})$ are those in w.} Obviously $\rm A \subset S(\{\emptyset\}).$ 
\pars
Assume the language $\r P_0$ is inductively defined from the set of atoms 
d. Recall that for our axioms ${\cal X} = {\cal D} \to {\cal F}$$,$ the 
strongest connective in $\cal X$ is $\to.$ While in $\cal D,$ or $\cal F$ when 
applicable$,$ the strongest connective is $\land.$  
Since $\emptyset \subset 
\rm \{ w \},$ it follows that $\rm S(\{w\}) = S(\emptyset) \cup S(\{w\}).$ Let 
$\rm b \in S(\emptyset).$ The only steps in the formal proof for b contain 
axioms\noindent 
 or follows from modus ponens. Suppose that step $\rm B_k =b$ is the 
first modus ponens step obtained from steps $\rm B_i,\ B_j,\ i,j < k,$ 
where $\rm B_i =A \to b,\ B_j = A.$ The strongest connective for each axiom is 
$\to.$ However$,$ since $A \to b$ is an axiom$,$ the strongest connective in 
$\rm A$ is $\land.$ This contradicts the requirement that A must also be an axiom 
with strongest connect $\to.$ Thus no modus ponens step can occur in a formal 
proof for b. Hence$,$ (\b 2) $\rm A = S(\emptyset).$ (No modus ponens step can occur 
using two axioms.)\pars
\vfil
\centerline{}
\eject
Let $\rm B_k= b_1 \in P_0$ and suppose (a) that $\rm b_1 = w,$ or (b) $\rm b_1 
\not= w$ and is 
the first nonaxiom step that appears 
in a formal 
demonstration from the hypothesis w. Assume (b). Then all steps $\rm B_i 
\in \{w\} \cup A,
\ 0\leq i < k.$  
Then the  only way that $\rm b_1$ can be 
obtained is by means of modus ponens.  However$,$ all other steps$,$ not including 
that which is w$,$ are axioms. No modus ponens step can occur using two axioms. 
Thus one of the steps used for modus ponens must not be an axiom. The only
nonaxiom that occurs prior to the step $\rm B_k$ is the step $\rm B_m = w.$ 
Hence$,$ one of the steps required for $\rm B_k$  must be $\rm B_m = w.$ The 
other step must be an axiom of the form $\rm w \to b_1$ and $\rm b_1 \not= w.$
Thus$,$ from the definition of 
the axioms (\b 3) $\rm b_1$ is either a finite $(\geq 1)$ conjunction of 
atoms in $\r d^\prime,$ or a single member of $\r d^\prime.$           
 Assume strong 
induction. Hence$,$ for $\rm n > 1,$ statement ({\bf 3}) holds for all $\rm r,\ 
1\leq r\leq n.$  
A
similar argument shows that ({\bf 3}) holds for the $\rm b_{n+1}$ nonaxiom 
step. Thus by induction$,$ ({\b 3}) holds for all nonaxiom steps. \pars
Hence$,$ there exists a $\rm Q \subset C$ such that each member of $\r Q$ is 
composed of  finitely many $(\geq 1)$ distinct members of $\rm d^\prime$ and 
the set $\rm G(Q)$ of all the proposition atoms that appear in  any member of $\rm 
Q = d^\prime = J$ since $\rm w \in Q.$ Moreover$,$ (\b 4) $\rm S(\{w\}) = A \cup 
d^\prime \cup Q$ and (\b 5) $\rm A,\ d^\prime,\ Q$ are mutually 
disjoint.\pars
\line{\hfil$ \forall x(x \in {\bf M_d - d} \to \exists y\exists z((y \in 
F(\b d)) \land (z \subset \b C)\land (\b S(\{x\}) =$\hfil}\pars 
\line{\hfil $\b A \cup y\cup z)\land (\b A \cap y = \emptyset)\land 
(\b A \cap z = \emptyset)\land (x \in z)$\hfil}\pars
\line{(10.1.3)\hfil$\land(y\cap z= \emptyset)\land \b G(z) = y)).$\qquad \qquad\hfil}\parm
\noindent holds in $\cal M,$ hence also in $\Hyper {\cal M}.$ So$,$ let $w$ be 
an ultraword. Then there exists internal $Q_1 \subset \Hyper {\b C}, \ w \in Q_1$ and *-
finite $d_1^\prime \subset \hyper {\b d}$ such that
$\b d \subset \Hyper {\b S}(\{w\}) = 
\hyper {\b A} \cup d_1^\prime \cup Q_1;\ \hyper {\b A},\ d_1^\prime,\ Q_1$ are 
mutually disjoint and $\Hyper G(Q_1) =d_1^\prime = \hyper {\b J}.$ Hence$,$ $\b d \subset d_1^\prime.$\pars
Now to analyze the objects in $Q_1.$ Let $\rm d = \{F_i \mid i \in \nat\}.$ 
Consider a bijection $h\colon \nat \to \b d$ defined by $h(n) = {\bf F_n} = 
[f],$ where $f \in T^0$ is the special member of $\bf F_n$ such that 
$f = \{(0,f(0))\},\  f(0) =i({\rm F_n}) = q_n \in i[\r d].$ From the above 
analysis$,$ (A) $[g] \in \b S(\{w\}) - \b A - \b d,\ (w \in \bf M_d - d),$ iff there 
exist $k,\ j \in \nat$ such that $ k< j$ and $f_1^\prime \in i[\r P_0]^{2(j-
k)}$ such that $[f_1^\prime] = [g],$ and this leads to (B) that for each even 
$2p,\ 0\leq 2p \leq 2(j-k); \ f_1^\prime(2p) = q_{k+p} \in i[\r P_0] \subset 
A_1,\ [(0,q_{k+p})] \in \b d^\prime,$ all such $q_{k+p}$ being distinct. For 
each odd $2p+1$ such that $0\leq 2p+1 \leq 2(j-k),\ f_1^\prime(2p+1) = i({\rm
\sp and\sp}).$ Also (C) $h(p) \in h[[k,j]]$ iff there exists an even $2p$ such 
that $0\leq 2p\leq 2(j-k)$ and $f_1^\prime (2p) = h(p) =q_{k+p} \in i[\r 
P_0].$ [Note that 0 is considered to be an even number.] \pars
By *-transfer of the above statements (A)$,$ (B) and (C)$,$  $[g] \in Q_1$ iff there 
exists some $j,\ k \in \hypernat,\ k < j,$ and $f^\prime \in \Hyper (i[\r P_0])^{2(j-
k)}$ such that $[f^\prime]= [g]$ and $\Hyper h[[k,j]] \subset  \Hyper {\b d}.$ 
Moreover$,$ each $\Hyper h(r), \ r \in [k,j]$ is a distinct member of $\Hyper {\b 
d}.$ The conjunction ``codes'' for $i({\rm \sp and\sp}) \in A_1$ that are generated 
by each odd $2p+1$ are all the same and there are *-finitely many of them. 
Hence$,$ $Q_1$ is the *-finite $(\geq 1)$ conjunctions of distinct members of 
$d_1^\prime,$ no other *-propositions are utilized and since $\Hyper {\b G}(Q_1) = 
d_1^\prime,$ all members of $d_1^\prime $ are employed for these  
conjunctions. This completes the proof. \qed 
\vfil
\eject
{\bf Corollary 10.1.1.1} {\sl Let $ w \in \Hyper {\bf M_d} - \Hyper {\b d}$ be an 
ultraword for denumerable {\r d} such that $\b d \subset \Hyper {\b S}(\{w\}).$ 
Then $\Hyper {\b S}(\{w\}) \cap {\bf P_0} = \bf A \cup Q \cup d$ and $\rm  
A,\ Q,\ d$ are mutually disjoint. The set \r Q is composed of finite $\geq 1$ 
conjunctions of members of \r d and all of the members of \r d are employed to 
obtain these conjunctions.}\pars
Proof. Recall that due to the finitary character of our standard objects
$^\sigma {\b A} = \b A = \hyper {\b A} \cap {\bf P_0}$. In like manner$,$ since 
$\b d \subset d_1^\prime,$ $d_1^\prime \cap {\bf P_0} = \b d.$ Now ${\bf P_0} 
\cap Q_1$ are all of the standard members of $Q_1.$ For each $ k \in 
\hypernat,\ \Hyper h(k) = {\bf F_k} \in \Hyper {\b d}$ and conversely. 
Further$,$ $\bf F_k \in d$ iff $ k\in \nat.$ Restricting $k,\ j \in \nat$ in the 
above theorem yields standard finite $\geq 1$ conjunctions of standard members 
of $d_1^\prime;$ 
hence$,$ members of $\b d.$ Since ultraword $w \in Q_1,$ we know that there  
exists some $\eta \in \hypernat-\nat$ and $f_1^\prime \in \Hyper (i[\r P_
0])^{2\eta},$ where $f_1^\prime$ satisfies the *-transfer of the properties 
listed in the above theorem . Since finite conjunctions of standard members of 
$d_1^\prime$ are *-finite conjunctions of members of $d_1^\prime$ and $\b d = 
\b d \cap d_1^\prime,$ it follows that all possible finite conjunctions of 
members of \b d that are characterized by the function $f_1^\prime \in i[{\rm 
P_0}]^{2(j-k)}$ are members of $Q_1$ for each such $j,\ k < \eta.$ Also for 
such $j,\ k$ the values of $f_1^\prime$ are standard.                                                            
On the other hand$,$ any value of $f_1^\prime$ is nonstandard iff it 
corresponds to a 
member of $d_1^\prime - \b d.$ Thus $Q_1 \cap {\bf P_0} = \b Q$ and this 
completes the proof. \qed
If it is assumed that each member of d describes a Natural event (i.e. 
N-event) at times indicated by $\r X_i,$ dropping the $\r X_i$ may still yield a 
denumerable developmental paradigm without specifically generated symbols such 
as the ``$i.$'' Noting that $d_1^\prime$ is *-finite and internal leads to the 
conclusion that we can have little or no knowledge about the word-like 
construction of each member of $d_1^\prime-\b d.$ These pure nonstandard 
objects 
can be considered as describing pure NSP-world events$,$ as will soon be 
demonstrated.  
Therefore$,$ it is important to understand the following 
interpretation scheme$,$ where descriptions are corresponded to events. \parm
{\leftskip 0.5in \rightskip 0.5in \noindent {\it Standard or internal 
NSP-world events or sets of events are interpreted as directly or indirectly 
influencing N-world events. Certain external objects$,$ such as the standard part 
operator$,$ among others$,$ are also interpreted as directly or indirectly 
influencing N-world events.}\par}\parm 
Notice that standard events can directly or indirectly affect standard 
events. In the micro-world$,$ the term {{\it indirect evidence}} or verification 
is a different idea than indirect influences. You can have direct or indirect 
evidence of direct or indirect influences when considered within the N-world. 
An indirect influence occurs when there exists$,$ or there is assumed to 
exist$,$ a  mediating  ``something'' between two events. Of course$,$ indirect 
evidence refers to  behavior that can be observed by normally 
accepted human sensors as such  behavior is  assumed to be caused by unobserved 
events. 
However$,$ the evidence for pure NSP-world events that directly or indirectly 
influence N-world events must be indirect evidence under the above 
interpretation. \pars
In order to formally  consider NSP-world events for the formation of objective 
standard reality$,$ proceed as follows: let $\cal O$ be the subset of $\cal W$
that describes those Natural events that are used to obtain developmental or general 
paradigms and the like. Let ${\rm E_j} \in \cal O.$ Linguistically$,$ assume that 
each $\rm E_j$ has the spacing symbol $\sp$ 
immediately to the right. Thus 
within each $\r T_i,$ there is a finite symbol string ${\rm F_i = E_i} \in \cal O$ that can be 
joined by the justaposition (i.e. join) operation to other event descriptions. 
Assume that 
${\cal W}_1$ is the set of nonempty symbol strings (with repetitions) formed 
from members of $\cal O$ by the join operation. These finite strings of symbols 
generate the basic elements for our partial sequences. \pars
Obviously$,$ ${\cal W}_1 \subset {\cal W}.$ Consider $\r T_i^\prime = \{\r X\r 
W_i \mid \r X \in {\cal W}_1\}$ and note that in many applications the time 
indicator $\r W_i$ need not be of significance for a given $\rm E_j$
 in some of the strings. Obviously$,$ $\r T_i^\prime \subset \r T_i$ for each 
$i.$ For our isomorphism $i$  onto $A_1,$ the following hold.\parm
\line{ \hfil $ \forall y(y\in {\cal E}\to  (y\in \b T_i^\prime \iff \exists x\exists 
f\exists w((\emptyset \not= w \in F(i[{\cal O}])) \land (x \in \nat)\land 
 $\hfil}\pars
\line{\hfil $(f(0)=i[W_i])\land (f \in P)\land\forall z((z\in \nat)\land (0<z\leq x) \to $\hfil}\pars
\line{(10.1.4)\hfil $f(z) \in w)\land (f \in y  
)))).$\qquad\hfil }\pars
\line{\hfil$\forall x( x \in \nat \to \exists f\exists w((\emptyset \not= w 
\in F(i[{\cal O}])) \land (f \in P) \land$\hfil}\pars
\line{(10.1.5)\hfil $\forall z(z \in \nat \to (0< z\leq x \iff f(z) \in w)))).$\qquad\hfil }\pars
\line{\hfil$\forall w(\emptyset \not= w \in F(i[{\cal O}])\to \exists x\exists 
f\exists y((f \in P)\land (x \in \nat) \land(y \in \b T_i')\land$\hfil }\pars
\line{(10.1.6)\hfil $(f \in y)\land\forall z(z \in \nat \to (0< z\leq x \iff f(z) \in w)))).$ \qquad\hfil }\parm  
Since each finite segment of a developmental pardigm corresponds to a member of 
${\b T_i'},$ each nonfinite hyperfinite segment should correspond to a member of $\Hyper (\b T_i^\prime) - \b T_i^\prime$ and it should be certain individual segments of such members of $\Hyper (\b T_i^\prime) - \b T_i^\prime$ that correspond to the ultranatural events produced by an ultraword; UN-events that cannot be 
eliminated from an NSP-world developmental paradigms. [Note: For a scientific language, 10.1.4 - 10.1.6 and other such statements correspond to a ${\cal W}'$ as generated by, at least, a denumerable alphabet as used in 9.2, 9.3.]\par
\medskip
\leftline{\bf 10.2 Laws and Rules.}
\medskip 
One of the basic requirements of human mental activity is the ability to 
{recognize the symbolic differences} between finitely long strings of symbols as 
necessitated by our reading ability and to apply linguistic rules finitely 
many times. G\"odel numberings specifically utilize such recognitions and the 
rules for the generation of recursive functions must be comprehended with 
respect to finitely many applications. Observe that G\"odel number recognition 
is an   ``ordered'' process while some fixed intuitive order is not necessary 
for the application of the rules that generate recursive functions.\pars
In general$,$ the simplest  ``rule'' for ordered or unordered {finite human
choice}$,$ a rule that is assumed to be humanly comprehensible by finite 
recognition$,$ is to simply  {{\it list}} the results of our choice (assuming  
that they are symbolically representable in some fashion) as a partial finite 
sequence for ordered choice or as a finite set of finitely long symbol strings 
for an unordered choice. Hence$,$ the end result for a finite choice can itself 
be considered as an algorithm ``for that choice only.'' The next application 
of such a {{\it finite choice rule}} would yield the exact same partial 
sequence or choice set. Another more general rule would be a statement which 
would say that you should  ``choose a specific number of objects'' from a fixed 
set (of statements). Yet$,$ a more general rule would be that you simply are 
required to  ``choose a finite set of all such objects,'' where the term  
``finite'' is intuitively known. Of course$,$ there are numerous specifically 
described algorithms that will also yield finite choice sets.\pars
From the symbol string viewpoint$,$ there are trivial machine programmable  
algorithms that allow for the comparison of finitely long symbols with each 
member of a finite set of symbol strings B that will determine whether or not 
a specific symbol string is a member of B. These programs duplicate the 
results of human symbol recognition. As is well-known$,$ there has not been an 
algorithm described that allows us to determine whether or not a given finite 
symbol string is a member of the set of all theorems of such theories as formal {Peano 
Arithmetic}. If one accepts {Church's Thesis}$,$ then no such algorithm will ever 
be described. \pars
Define the general finite human choice relation on a set $A$ as $H_0(A) = 
\{(A,x) \mid x \in F_0(A)\},$ where $F_0$ is the finite power set operator 
(including the empty set = no choice is made). 
Obviously$,$ the inverse $H_0^{-1}$ is a function from $F(A)$ onto $\{A\}.$ 
There are choice operators that produce sets with  a specific number of elements 
that can be easily defined. Let $F_1(A)$  be the set of all singleton subsets 
of $A.$  The axioms of set theory state that such a set of singleton sets 
exists. Define $H_1(A) = \{(A,x) \mid x \in F_1(A) \},$ etc. Considering such functions as defined on sets $X$ that are members of 
a  superstructure$,$ then these relations are subsets of $\power X \times \power 
X$ and as such are also members of the superstructure.\pars
Let $\r A =\r P_0.$ Observe that $^\sigma {\bf H_0(A)} = \{ (\hyper {\b A}, x) \mid 
x \in F_0(\b A)\}$ and $\Hyper {\bf H_i(A)} = \{(\hyper {\b A}, x) \mid x \in \Hyper 
{(F_i(\b A))}\}\ (i\geq 0).$ Now $\Hyper {(F_0 (\b A))} = \Hyper {F_0}(\hyper 
{\b A})$ is the set of all *-finite subsets of $\hyper {\b A}.$ On the other 
hand$,$ for the $i > 0$ cardinal subsets$,$ $\Hyper {(F_i (\b A))} = F_i(\hyper {\b A})$ for each $ i \geq 1.$ 
With respect to an ultraword $w$ that generates the general and developmental 
paradigms$,$ we know that $ w \in \Hyper {\bf P_0} - {\bf P_0}$ and that 
$(\Hyper {\bf P_0}, \{x\}) \in \Hyper {\bf H_1}(\bf P_0).$ The actual finite 
choice operators are characterized by th set-theoretic second projector 
operator $P_2$ as it is defined on $H_i(\r A).$ This operator embedded by the 
injection $\theta$ is the same as $P_2$ as it is defined on $\bf H_i(A).$ Thus$,$ 
when $h = (\b A, x) \in \bf H_i(A),$ then we can define $ x = P_2(h) = C_i(h) 
= {\bf C_i}(h).$ The maps $C_i$ and $\bf C_i,$ formally defined below$,$ are the 
specific finite choice operators. For consistency$,$ we let $C_i$ and $\bf C_i$ 
denote the appropriate finite choice operators for $H_i(\r A)$ and $\bf H_i(A)$$,$ 
respectively. \pars
Since the $\Hyper P_2$ defined on say $\bf H_i(A)$ is the same as the 
set-theoretic second projection operator $P_2,$ it would be possible to denote 
$\Hyper {\bf C_i}$ as $\bf C_i$ on internal objects. For consistency$,$ the 
notation $\Hyper {\bf C_i}$ for these special finite choice operators is 
retained.  
Formally$,$ let ${\bf C_i\colon H_i(A)} \to F_i(\b A).$ Observe that $^\sigma {\bf 
C_i} = \{\Hyper (a,b) \mid (a,b) \in {\bf C_i}\} =\{((\hyper {\b A},b),b)\mid 
b\in F_i(\b A)\} \subset \Hyper {\bf C_i};$ and$,$ for $b \in F_i(\b A),\ {\bf 
C_i}((\b A,b)) = b$ implies that $^\sigma({\bf C_i}((\b A,b))) = \{\hyper 
a\mid a \in b\} = b$ from the construction of $\cal E.$ Thus in 
contradistinction to the consequence operator$,$ for each $(\hyper {\b A},b)\in 
{^\sigma {\bf H_i}},$ the image $(^\sigma {\bf C_i})((\hyper {\b A},b)) = 
{^\sigma({\bf C_i}((\b A, b)))} = (\Hyper {\bf C_i})((\hyper {\b A},b)) = b.$ 
Consequently$,$ the set map $^\sigma{\bf C_i}\colon {^\sigma{\bf H_i}} \to 
F_i(\b A) = {^\sigma(F_i(\b A))}$ and $\Hyper {\bf C_i} \mid {^\sigma {\bf 
H_i}} = {^\sigma {\bf C_i}}.$ Finally$,$ it is not difficult to extend these 
finite choice results to general internal sets.\pars
In the proofs of such theorems as 7.2.1$,$ finite and other choice sets are 
selected due to their set-theoretic existence. The finite choice operators 
$C_i$ are not specifically applied since these operators are only intended as 
a mathematical model for apparently effective human processes --- procedures 
that generate acceptable algorithms. As is well-known$,$  
there are other describable rules that also lead to finite or infinite 
collections of statements. Of course$,$ with respect to a G\"odel  encoding 
$i$ for the set of all words $\cal W$ the finite choice of readable sentences 
in $\cal E$ is one-to-one and effectively related to a finite and$,$ hence$,$ 
recursive subset of $\nat.$\pars
From this discussion$,$ the descriptions of the finite choice operators would 
determine a subset of the set of all algorithms (``rules'' written in the 
language $\cal W$) that allow for the selection of readable sentences. Notice 
that before algorithms are applied there may be yet another set of readable 
sentences that yields conditions that must exist prior to an application of 
such an algorithm and that these application rules can be modeled by members 
of $\cal E.$ \pars
In order to be as unbiased as possible$,$ it has been required for N-world 
applications that the set of all frozen segments be infinite. Thus$,$ within the 
proof of Theorem 7.2.1$,$ every N-world developmental$,$ as well as a general 
paradigm$,$ is a proper subset of a {*-finite NSP-world paradigm}$,$ and the *-finite 
paradigm is obtained by application of the *-finite choice operator $\Hyper 
{\bf C_0}.$ As has been shown$,$ such *-finite paradigms contain pure 
unreadable (subtle) sentences that may be interpreted for developmental 
paradigms as {pure refined NSP-world behavior} and for general paradigms as 
specific pure NSP-world ultranatural events or objects.\pars
Letting $\Gamma$ correspond to the formal theory of Peano Arithmetic$,$ then assuming 
Church's Thesis$,$ there would not exist a N-world algorithm (in any human 
language) that allows for the determination of whether or not a statement F in 
the formal language used to express $\Gamma$ is a member of $\Gamma.$ By 
application of the *-finite choice operator $\Hyper {\bf C_0},$ however$,$ there 
does exist a *-finite $\Gamma^\prime$ such that $^\sigma\Gamma = \Gamma 
\subset \Gamma^\prime$ and$,$ hence$,$ within the NSP-world a ``rule'' that allows 
the determination of whether or not $\r F \in \Gamma^\prime.$ If such internal 
processes mirror the only allowable procedures in the NSP-world for such a 
``rule,'' then it might be argued that we do not have an effective NSP-world 
process that determines whether or not F is a member of $\Gamma$ for $\Gamma$ 
is external. \pars 
As previously alluded to at the beginning of this section$,$ 
when a G\"odel encoding $i$ is utilized with the N-world$,$ the injection $i$ is 
not a surjection. When such G\"odel encodings are studied$,$ it is usually {\it 
assumed}$,$ without any further discussion$,$ that there is some human mental 
process that allows us to recognize that one natural number representation (whether in prime 
factored form or not) is or is not distinct from another such representation. 
It is not an unreasonable assumption to assume that the same {\it effective} 
(but external) process exists within the NSP-world. Thus within the {NSP-
world there is a  ``process'' that determines whether or not an object is a 
member of} $\hypernat -\nat= \nat_\infty$ or $\nat.$ Indeed$,$ from the ultraproduct 
construction of our nonstandard model$,$ a few differences can be detected by 
the human mathematician. Consequently$,$ this assumed NSP-world effective 
process would allow a determination of whether or not $\b F = [f_m]$ is a 
member of $\Gamma$ by recalling that $f_m \in P_m$ signifies that $[f_m] \in 
\Hyper {\Gamma} -\Gamma$ implies $m \in \nat_\infty\simeq
\Hyper {({i[{\cal W}]})} - {i[{\cal W}]} = \hyper A_1 - A_1.$\pars 
The above NSP-world recognition process is 
equivalent$,$ as defined in Theorem 7.2.1$,$ to various applications of a single 
(external) set-theoretic intersection. Therefore$,$ there are internal 
processes$,$ such as $\Hyper {\bf C_0},$ that yield pure NSP-world developmental 
paradigms and a second (external) but acceptable NSP-world effective process 
that produces specific N-world objects. Relative to our modeling procedures$,$ 
it can be concluded that both of these processes are {intrinsic ultranatural 
processes.} \pars 
With respect to Theorem 10.1.1$,$ the NSP-world developmental or general 
paradigm generated by an ultraword is *-finite and$,$ hence$,$ specifically NSP-
world obtainable prior to application of $\Hyper {\b S}$ through application 
of $\Hyper {\bf C_0}$ to $\Hyper {\b d}.$ However$,$ this composition can be 
reversed. The NSP-world (IUN) process $\Hyper {\bf C_1}$ can be applied to the 
appropriate $\Hyper {\bf M_d}$ type set and an appropriate ultraword $w \in 
\Hyper {\bf M_d}$ obtained. Composing $\Hyper {\bf C_1}$ with $\Hyper {\b S}$  
would yield $d_1^\prime$ in a slightly less conspicuous manner. Obviously$,$ 
different ultrawords generate different standard and nonstandard developmental 
or general paradigms.\pars 
To complete the actual mental-type processes that lead to the proper ordered 
event sequences$,$ the above discussion for the finite choice operators is 
extended to the human mental ability of ordering a finite set in terms of 
rational number subscripts. New choice operators are defined that model not 
just the selection of a specific set of elements that is of a fixed finite
cardinality but also choosing the elements in the required rational number 
ordering. 
The ultrawords $w$ that exist are *-finite in 
length. By application of the inverses of the $f$ and $\tau$ functions of 
section 7.1$,$ where they may be considered as extended standard functions 
$\hyper f$ and $\hyper \tau$$,$ there would be from analysis of extended 
theorem 7.3.2 a hyperfinite set composed of standard or nonstandard frozen 
segments contained in an ultraword. Further$,$ in theorem 7.3.2$,$ the chosen 
function $f$ does not specifically differentiate 
each standard or nonstandard frozen 
segment with respect to its ``time'' stamp subscript. There does exist$,$ 
however$,$ another function in the *-equivalence class $[g]=w$ that will make 
this differentiation. It should not be difficult to establish that after 
application of the ultralogic $\Hyper {\b S},$ there is applied an 
appropriate 
mental-like hyperfinite ordered choice operator (an IUN-selection 
process) and that this would yield that various types of event sequences. 
Please note that each event sequence has a beginning point of observation. 
This point of observation need not indicate the actual moment when a specific 
Natural system began its development.\pars 
Various subdevelopmental (or subgeneral) paradigms $\rm d_i$ are obtained by considering the 
actual descriptive content (i.e. events) of specific theories $\rm\Gamma_i$ that 
are deduced from hypotheses $\rm\eta_i,$ usually$,$ by finitary consequence operators 
$\rm S_i$ (the inner logics) that are compatible with $\r S.$ In this case$,$ 
$\rm d_i \subset S_i(\eta_i).$ It is also possible to include within $
\{\rm d_i\}$ and $\{\rm \eta_i\}$ the assumed descriptive {chaotic behavior} that 
seems to have no apparent set of hypotheses except for that particular 
developmental paradigm itself and no apparent deductive process 
except for the identity consequence operator. In this way$,$ such scientific 
nontheories can still be considered as a formal theory produced by a finitary 
consequence operator applied to an hypothesis. Many of these hypotheses 
$\rm \eta_i$ contain the so-called natural laws (or first-principles) peculiar to 
the formal theories $\rm \Gamma_i$  and the theories language$,$ where it is assume 
that such languages are at least closed under the informal conjunction and 
conditional.\pars
Consider each $\rm \eta_i$ to be a general paradigm. For the appropriate M type 
set constructed from the denumerable set $\r B = \{\bigcup\{\rm d_i\} 
\cup(\bigcup \{\eta_i\}),$ redefine $\rm M_B$ to be the smallest subset 
of $\rm P_0$ containing B and closed under finite $(\geq 0)$ conjunction. (The 
usual type of inductively defined $\rm M_B.$) Then there exist ultrawords $w_i 
\in \Hyper {\bf M_B} - \Hyper {\b B}$ such that ${\bf \eta_i} \subset \Hyper {\b 
S}(\{w_i\})$ (where due to parameters usually {\it ultranatural laws} exist in $\Hyper {\b S}(\{w_i\})- {\bf \eta_i}$) and ${\bf d_i} \subset \Hyper {\b S}(\{w_i\}).$ Using methods 
such as those in Theorem 7.3.4$,$ it follows that there exists some $w" \in 
\Hyper {\bf M_B} - {\bf M_B}$ such that $ w_i \in \Hyper {\b  S}(\{w"\})$ and$,$ 
consequently$,$ ${\bf \eta_i \cup d_i} \subset \Hyper {\b S}(\{w"\}).$  
Linguistically$,$ it is hard to describe the ultraword $w".$ Such a $w"$ might 
be called an {{\it ultimate  ultranatural hypothesis}} or {{\it the ultimate 
building plain}}. \pars
{\it Remark.} It is not required that the so-called Natural laws that appear in some 
of the $\rm \eta_i$ be either cosmic time or universally applicable. They could 
refer only to local first-principles. It is not 
assumed that those first-principles that display themselves in our local 
environment are universally space-time valid.\pars
Since the consequence operator $\r S$ is compatible with each $\rm S_i,$ 
it is useful to proceed in the following manner. First$,$ apply the IUN-process
$\Hyper {\b S}$ to $\{w"\}.$ Then ${\bf d_i \cup \eta_i} \subset \Hyper {\b 
S}(\{w"\}).$ It now follows that ${\bf d_i \cup \eta_i} \subset \Hyper {\b 
S}(\{w_i\}) \subset \Hyper {\bf S_i}(\Hyper {\b S}(\{w_i\})) \subset \Hyper 
{\bf S_i}(\Hyper {\bf S_i}(\{w_i\})) = \Hyper {\bf S_i}(\{w_i\}).$ Observe 
that for each $\rm a \in \Gamma_i$ there exists some finite $\rm F_i \subset 
\eta_i$ such that $\bf a \in S_i(F_i).$ However$,$ $\bf F_i \subset \eta_i$ for 
each member of $F(\rm \eta_i)$ implies that $\b a \in \Hyper {\bf S_i(F_i)} 
\subset \Hyper {\bf S_i}(\Hyper {\b S}(\{w_i\})).$ Consequently$,$ ${\bf 
\Gamma_i} \subset \Hyper {\bf S_i}(\Hyper {\b S}(\{w_i\})).$ The {ultimate 
ultraword} suffices for the descriptive content and inner logics associated 
with each theory $\rm \Gamma_i$.\pars
We now make the following observations relative to ``rules'' and deductive 
logic. It has been said that science is a combination of empirical data$,$ 
induction and deduction$,$ and that you can have the first two without the last. That this belief is totally false should be 
self-evident since the philosophy of science requires its own general rules 
for observation$,$ induction$,$ data collection$,$ proper experimentation and the 
like. All of these general rules require logical deduction for their 
application to specific cases --- the metalogic. Further$,$ there are specific 
rules for linguistics that also must be properly applied prior to scientific 
communication. Indeed$,$ we cannot even open the laboratory door --- or at least 
describe the process --- without application of deductive logic. The concept 
of deductive logic as being the patterns our ``minds'' follow and its use 
exterior to the inner logic of some theory should not be dismissed for even 
the (assumed?) mental methods of human choice that occur prior to 
communicating various scientific statements and descriptions. \pars 
Finally$,$ with respect to the {hypothesis rule} in [9]$,$ it 
might be argued that we can easily analyze the specific composition of all 
significant ultrawords$,$ as has been previously done$,$ and the composition of 
the nonstandard extension of the general paradigm. Using this assumed analysis 
and an additional alphabet$,$ one {\it might} obtain specific information about  
pure NSP-world ultranatural laws or refined behavior.  Such an argument would 
seem to invalidate the cautious hypothesis rule and lead to appropriate 
speculation. However$,$ such an argument would itself be invalid.\pars 
 Let ${\cal W}_1$ be an infinite set of meaningful readable sentences for some 
description and assume that ${\cal W}_1$ does not contain any infinite subset of 
readable sentences each one of which 
contains a mathematically interpreted entry such as a real number or the like. 
Since ${\cal W}_1 \subset {\cal W}$ and the totality $\r T_i = \{\r X
\r W_i \mid \r X \in {\cal W} \}$ is denumerable$,$ the subtotality 
${\r T}_i^\prime = \{\r X{\r W}_i \mid \r X \in {\cal W}_1 \}$ is also 
denumerable. Hence$,$ the external cardinality of $\Hyper {\bf T_i^\prime} \geq 
{\vert{\cal M}\vert}^+.$ \pars
Consider the following sentence \parm
\line{\hfil $ \forall z(z \in i[{\cal W}_1] \to 
\exists y\exists x((y \in 
A_1^{[0,1]}) \land (x \in {\bf T_i^\prime}) \land (y \in x) \land 
 $\hfil}\pars
\line{\hfil$((0,i({\r W_i})) \in y)\land
((1,z) \in y )))).$\hfil (10.2.1)}\parm
By *-transfer and letting  ``z'' be an element in $\Hyper (i[{\cal W}_1])-i[{\cal W}_1]$ it 
follows that we can have little knowledge about the remaining and what must be 
unreadable portions that take the ``X'' position. If one assumes 
that members of ${\cal W}_1$ are possible descriptions for possible NSP-world 
behavior at the time $t_i,$ then it may be assumed that at the time $t_i$ the 
members of$ \Hyper {\bf T_i^\prime} - {\bf T_i}$ describe NSP-world behavior at NSP-
world (and N-world) time $t_i.$ Now as i varies over $\hypernat,$ pure 
nonstandard subdevelopmental paradigms (with or without the time index 
statement $\r W_i$) exist with members in $\Hyper {\cal T}$ and may be 
considered as descriptions for time refined NSP-world behavior$,$ especially for 
a NSP-world time index $i \in \nat_\infty.$\pars   
\bigskip
\centerline{\bf CHAPTER 10 REFERENCES}
\medskip
\noindent {\bf 1} Beltrametti$,$ E. G.$,$ Enrico$,$ G. and G. Cassinelli$,$ {\it The 
Logic of Quantum Mechanics,} Encyclopedia of Mathematics and Its Applications$,$ 
Vol. 15$,$ Addison-Wesley$,$ Reading$,$ 1981.\pars
\noindent {\bf 2} d'Espagnat$,$ B.$,$ The quantum theory and realism$,$ Scientific 
America$,$ 241(5)(1979)$,$ 177. \pars
\noindent {\bf 3} {\it Ibid.} \pars
\noindent {\bf 4} {\it Ibid.}$,$ 181. \pars
\noindent {\bf 5} {\it Ibid.}$,$ \pars
\noindent {\bf 6} {\it Ibid.}$,$ 180.\pars
\noindent {\bf 7} Feinberg$,$ G.$,$ Possibility of faster-than-light particles$,$ 
Physical Review$,$ 159(5)(1976)$,$ 1089---1105.\pars
\noindent {\bf 8} Hanson$,$ W. C.$,$ The isomorphism property in nonstandard 
analysis and its use in the theory of Banach Spaces$,$ J. of Symbolic Logic$,$ 
39(4)(1974)$,$ 717---731.\pars
\noindent {\bf 9} Herrmann$,$ R. A.$,$ D-world evidence$,$ C.R.S. Quarterly$,$ 
23(2)(l1986)$,$ 47---54.\pars
\noindent {\bf 10} Herrmann$,$ R. A.$,$ The Q-topology$,$ Whyburn type filters and 
the cluster set map$,$ Proceedings Amer. Math. Soc.$,$  59(1)(1975)$,$ 161---
166.\pars
\noindent {\bf 11} Kleene$,$ S. C.$,$ {\it Introduction to  Metamathematics}$,$ D. 
Van Nostrand Co.$,$ Princeton$,$ 1950.\pars
\noindent {\bf 12} Prokhovnik$,$ S. J.$,$ {\it The Logic of Special Relativity}$,$ 
Cambridge University Press$,$ Cambridge$,$ 1967.\pars
\noindent {\bf 13} Stroyan$,$ K. D. and W. A. J. Luxemburg$,$ {\it Introduction to 
the Theory of Infinitesimals}$,$ Academic Press$,$ New York$,$ 1976. \pars
\noindent {\bf 14} Tarski$,$ A.$,$ {\it Logic$,$ Semantics$,$ Metamathematics}$,$ 
Clarendon Press$,$ Oxford$,$ 1969.\pars
\noindent {\bf 15} Thurber$,$ J. K. and J. Katz$,$ Applications of fractional 
powers of delta functions$,$ {\it Victoria Symposium on Nonstandard Analysis},
Springer-Verlag$,$ New York$,$ 1974.\pars
\noindent {\bf 16} Zakon$,$ E.$,$ Remarks on the nonstandard real axis$,$ {\it 
Applications of Model Theory to Algebra$,$ Analysis and Probability}$,$ Holt$,$ 
Rinehart and Winston$,$ New York$,$ 1969. \pars
\noindent {\bf 17} Note that the NSP-world model is not a local hidden 
variable theory.
\vfil
\eject
\centerline{\bf 11. "Things"}
\bigskip
\leftline{\bf 11.1 Propertons (Subparticles).}\par 
\medskip
What is a properton? Or$,$ what is an infant? Or$,$ better still$,$ what is a 
thing? I first used the name {infant} for these strange objects. I then 
coined the term subparticle and have even used the term "things." These three names do not convey the exact intuitively mean and are prone to incorrect mental images. As of 11 July 2012, the term "properton" is employed. As discussed, they carry physical or physical-like "properties" in a coded form. (As stated in [9]$,$ these objects are not to be described in terms 
of any geometric configuration. These {multifaceted things}$,$ these propertons$,$ 
are not to be construed as either particles nor waves nor quanta nor anything 
that can be represented by some fixed imagery. Propertons are to be viewed only 
operationally. Propertons are only to be considered as represented by a 
*-finite sequence $\{a_i\}_{i=1}^n,\ n \in \hypernat,$ of hyperreal numbers. Indeed$,$ the idea of the n-tuple 
$(a_1,a_2,\ldots,a_i,\ldots)$ notation is useful and we assume that 
$n$ is a fixed member of $\nat_\infty.$ The language of coordinates 
for this notation is used$,$ where the i'th coordinate means the i'th value of the 
sequence. Obviously$,$ 0 is not a domain member for our sequential representation. 
\pars
The first coordinate 
$a_1$ is a ``naming'' coordinate. The remaining coordinates are used to 
represent various real numbers$,$ complex numbers$,$ vectors$,$ and the like physical
qualities needed for different physical theories. For example$,$ $a_2 = 1$ might 
be a counting coordinate. Then $a_i,\ 3\leq i\leq 6$ are hyperreal numbers that
represent NSP-world 
coordinate locations of the properton named by $a_1$ --- $a_7,\ a_8$  
represent the positive or negative charges that can be assigned to every 
properton --- $a_9,\ a_{11},\ a_{13}$ {hyperreal representations} for the inertial$,$ 
gravitational and intrinsic (rest) mass. For vector quantities$,$ continue this 
coordinate assignment and assign specific coordinate locations for the vector 
components. So as not to be biased$,$ include as other coordinates hyperreal 
measures for qualities such as energy$,$ apparent momentum$,$ and all other physical  
qualities required within theories that must be combined in order to produce a 
reasonable description for N-world behavior. For the same reason$,$ we do not 
assume that such N-world properties as the uncertainty principle hold for 
the NSP-world. (See note (2) on page 116.)\pars
It is purposely assumed that the qualities represented by the coordinate 
$a_i, \ i \geq 3$ are not inner-related$,$ in  their basic construction$,$ by any 
mathematical relation since it is such {inner-relations} that are assumed to 
mirror the N-world laws that govern the development of not only our present 
universe but previous as well as future developmental alterations. The same 
remarks apply to any possible and distinctly \noindent different universes that may or 
not occur. Thus$,$ for these reasons$,$ we view the properton as being totally 
characterized by such a sequence $\{a_i\}$ and always proceed cautiously when 
any attempt is made to describe all but the most general properton 
behavior. Why have we chosen to presuppose that propertons are 
characterized by sequences$,$ where the coordinates are hyperreal 
numbers?\pars
For chapters 11$,$ 12 assume EGS. Let $r$ be a positive real number. The number $r$ can be represented by a 
decimal-styled number$,$ where for uniqueness$,$ the repeated 9s case is used for all 
terminating decimals. From this$,$ it is seen that there is\ a sequence $S_i$ of 
natural numbers such that $S_i/10^i \to r.$ Consequently$,$ for any $ \omega \in 
\nat_\infty=\hypernat-\nat,$ it follows that $\pm \hyper S_\omega/10^\omega \in  \monad {\pm r},$  
where $\hyper S_\omega \in \hypernat$ and $\monad {\pm r}$ is the {monad} 
about $\pm r.$ In [9]$,$ it is assumed that each coordinate $a_i,\ i \geq 3$ is 
characterized by the numerical quantity $\pm 10^{-\omega},\ \omega \in 
\nat_\infty.$ Obviously$,$ we need not confine ourselves to the number $10^{-
\omega}.$  \parm

{\bf Theorem 11.1.1} {\sl For each $0 < i \in \nat$, let $0 < m_i \in \nat$ and $m_i \to \infty.$ Let any $\omega,\ \lambda \in \nat_\infty.$ Then, for each $r \in \real,$ there exists a $b/\hyper m_\omega \in
\{x/\hyper {m_\omega} \mid ( x \in \Hyper {\b Z}) \land(\vert x\vert < \lambda\hyper{m_\omega}\},$ where $\hyper m_\omega \in \nat_\infty,$ and $b/\hyper m_\omega \approx r$ (i.e. $b/\hyper m_\omega \in \mu(r)$). If $r \not= 0,$ then $\vert b \vert \in \nat_\infty.$}\pars

Proof. For $r \in \real,$ there exists a unique integer $n  \in \b Z$ such 
that $n\leq r< n+1.$ Partition $[n,n+1)$ as follows: for each $0< i \in \nat,$ and $0< m_i \in \nat,$ consider $[n, n+1/m_i),\ldots,[n + (m_i -1)/m_i,n +1).$ Then there 
exists a unique $c_i \in \{0,1,\ldots, m_i-1\}$ such that $r \in [n 
+c_i/m_i,n+(c_i+1)/m_i).$ Let $S_i= (m_in +c_i)/m_i = f_i/m_i.$ Since $0\leq r - S_i < 1/m_i$ and $m_i \to \infty,$ then $S_i \to r.$ 
This yields two sequences $S\colon \nat \to \b Q$ and $ f\colon \nat \to \b 
Z,$ where, for each $\omega \in 
\nat_\infty, \ \hyper S_\omega =\hyper f_\omega/\hyper m_\omega \approx r$ and $\hyper f_\omega \in \Hyper {\b 
Z}.$ Observe that $\hyper f_\omega/\hyper m_\omega$ is a finite (i.e. limited) number and $\hyper m_\omega \in \nat_\infty.$   
Hence$,$ $ \vert \hyper f_\omega/\hyper m_\omega \vert < \lambda$ entails that $\vert 
\hyper f_\omega \vert < \lambda\hyper m_\omega.$ Therefore$,$ $\hyper f_\omega/\hyper m_\omega \in \{x/\hyper {m_\omega} \mid ( x \in \Hyper {\b Z}) \land(\vert x\vert < \lambda\hyper{m_\omega}\}$. If $\hyper f_\omega \in \b Z,$ then $\hyper f_\omega/\hyper m_\omega \approx 0.$ \qed
{\bf Corollary 11.1.1.1} {\sl For each $0 < i \in \nat$, let $0 < m_i \in \nat$ and $m_i \to \infty.$
Let any $\omega,\ \lambda \in \nat_\infty.$ Then, for each $r \in \real,$ there is a sequence  $f\colon \nat \to \b Z$ such that $\hyper f_\omega/\hyper m_\omega \in \mu(r).$ There are unique $n\in\b Z,\ c_\omega = 0$ or $c_\omega \in \nat_\infty$ such that $c_\omega \leq \hyper m_\omega -1$ and $\hyper f_\omega = \hyper m_\omega n + c_\omega.$}
\parm

For the {{\it ultra-properton}}$,$ each coordinate $a_i = 1/10^\omega \ i 
\geq 3$ and odd$,$ $a_i =- 1/10^\omega \ i 
\geq 4$ and even$,$ $\omega \in \nat_\infty.$ From the above theorem$,$ the 
choice of $10^{-\omega}$ as the basic numerical quantity is for 
convenience only and is not unique accept in its infinitesimal character. 
Of course$,$ the sequences chosen to represent the ultra-properton are pure 
internal objects and as such are considered to directly or indirectly affect 
the N-world. Why might the *-finite ``length''
 of such propertons (here is 
where we have replaced the NSP-world entity by its corresponding sequence)  
be of significance? \pars
First$,$ since our N-world languages are formed from a finite set of alphabets$,$ 
it is not unreasonable to assume that NSP-world ``languages'' are composed from 
a *-finite set of alphabets. Indeed$,$ since it should not be presupposed that there is 
an upper limit to the N-world alphabets$,$ it would follow that the basic 
NSP-world set of alphabets is an infinite *-finite set. Although the 
interpretation method that has been chosen does not require such a restriction 
to be placed upon NSP-world alphabets$,$ it is useful$,$ for consistency$,$ to assume 
that descriptions for substratum processes that affect$,$ in either a directly 
or indirectly detectable manner$,$ N-world events be so restricted. For the 
external NSP-world viewpoint$,$ all such infinite *-finite objects have a very 
significant common property. Note: in what follows ${\cal M}_1$ is the 
extended superstructure constructed on page 70. \parm
\vfil
\eject
{\bf Theorem 11.1.2} {\sl All infinite *-finite members of 
our (ultralimit) model $\Hyper {\cal M}_1$ have the same external cardinality which is
$\geq \vert {\cal M}_1\vert.$} \pars
Proof. Hanson [8] and Zakon [16] have done all of the difficult work for this 
result to hold. First$,$ one of the results shown by Henson is that all
infinite *-finite  members 
of our ultralimit model have the same external cardinality. Since our model 
is a comprehensive enlargement$,$ Zakon's theorem 3.8 in [16] applies. Zakon 
shows that there exists a *-finite set$,$ $A,$ such that $\vert A \vert \geq 
\vert {\cal M}_1 \vert= \vert {\cal R}\vert.$ 
Since $A$ is infinite$,$ Hanson's result now implies that all 
infinite *-finite members of 
our model satisfy this inequality.\qed
For an extended infinite standard set $\hyper A$ it is well-known that 
$\vert \hyper A \vert \geq {\vert{\cal M}_1\vert}^+.$  One may use these
various results and establish easily that there exist more than enough 
propertons to obtain all of the cardinality statements relative to the three 
substratum levels that appear in [9]  even if we assume that there are a 
continuum of finitely many properton qualities that are needed to create all 
of the N-world. \pars
Consider the following infinite set of statements expressed in an extended 
alphabet. \parm
\line{\hfil${\rm G_A = \{An\sp elementary\sp particle\sp}
\r k^\prime(\r i^\prime,\r j^\prime)\sp{\rm with\sp}$\hfil}\pars
\line{\hfil${\rm total\sp energy\sp}\r c^\prime{+}1/(\r n^\prime).\mid ((i,j,n) \in $\hfil}\pars
\line{\hfil \qquad$\nat^+ \times 
\nat^+ \times \nat^+)\land (1\leq k\leq m) \},$\hfil (11.1.1)}\parm
\noindent where $\nat^+$ is the set of all nonzero  natural numbers and $m \in \nat^+.$ 
Applying the same procedure that appears in the proof of Theorem 9.3.1 and 
with a NSP-world alphabet$,$ we obtain\par
\line{\hfil${\rm G_A^\prime = \{An\sp elementary\sp particle\sp} 
\r k^\prime(\r i^\prime,\r j^\prime)\sp{\rm with\sp}$\hfil}\pars
\line{\hfil${\rm total\sp energy\sp}\r c^\prime{+}1/(\r n^\prime).\mid ((i,j,n) \in $\hfil}\pars
\line{\hfil \qquad$\hypernat^+ \times 
\hypernat^+ \times \hypernat^+)\land (1\leq k\leq m )\},$\hfil (11.1.2)}\parm
Assume that there is at least one type of elementary particle with the properties 
stated in the set ${\rm G_A}.$ It will be shown in the next section 
that within the NSP-world there may be simple properties that lead to  
N-world energy being a manifestation of mass. For $c = 0,$ we have  another 
internal set of descriptions that forms a subset of $\rm G_A^\prime.$ \parm
\line{\hfil${\rm \{An\sp elementary\sp particle\sp} 
\r k^\prime(\r i^\prime,\r j^\prime)\sp{\rm with\sp}$\hfil}\pars
\line{\hfil${\rm total\sp energy\sp}\r c^\prime{+}1/(10^{\zeta^\prime}).\mid 
((i,j,\zeta) \in $\hfil}\pars
\line{\hfil \qquad$\hypernat^+ \times 
\hypernat^+ \times \hypernat^+)\land (1\leq k\leq m )\},$\hfil (11.1.3)}\parm
For our purposes$,$ (11.1.3) leads immediately to the not ad hoc concept of 
propertons with {infinitesimal proper mass}. As will be shown$,$ such 
infinitesimal proper mass can be assumed to characterize any possible zero 
proper mass N-world entity. The set $\rm G_A^\prime$ has meaning if there 
exists at least one natural entity that can possess the energy expressed by 
$\rm G_A,$ where this energy is measured in some private unit of measure.\pars
Human beings combine together finitely many sentences to produce 
comprehensible descriptions. Moreover$,$ all N-world human construction requires 
the composition of objectively real N-world objects. We model the idea of {{\it
finite composition}} or {{\it finite combination}} by an N-world process. 
This produces a corresponding 
NSP-world intrinsic ultranatural process {{\it ultrafinite composition}} or 
{{\it ultrafinite combination}} that 
can either directly or indirectly affect the N-world$,$ where its effect is 
{indirectly inferred}. \pars 
Let the index $j$ vary over a hyperfinite interval and fix the other indices. 
Then the set of sentences\parm
\line{\hfil${\rm G_A^{\prime\prime} = \{An\sp elementary\sp particle\sp} 
\r k^\prime(\r i^\prime,\r j^\prime)\sp{\rm with\sp}$\hfil}\pars
\line{\hfil${\rm total\sp energy\sp}\r c^\prime{+}1/(\r n^\prime).
\mid (j \in \hypernat^+)$\hfil}\pars
\line{(11.1.4)\hfil$\land (1\leq j\leq \lambda) \},$\qquad\hfil}\parm
\noindent where $\lambda \in \hypernat^+,\ 3\leq i \in \hypernat,\ n \in\nat_\infty$ 
and $1 \leq 
k\leq m,$ forms an internal linguistic object that can be assumed to describe a 
hyperfinite collection of ultranatural entities. Each member of $\rm 
G_A^{\prime\prime}$ has the $i$'th coordinate that  measures the proper 
mass and  is infinitesimal (with respect to NSP-world private units of 
measure). In the {N-world$,$ finite combinations yield an event.} Thus$,$  with 
respect to such sets as $\rm G_A^{\prime\prime},$ one can say that there are 
such N-world events iff there are ultrafinite combinations of NSP-world entities. And 
such ultrafinite combinations yield a  NSP-world event that 
is an ultranatural entity. \pars
Associated with such ultrafinite combinations for the entities described in 
$\rm G_A^{\prime\prime}$ there is a very significant procedure that yields 
the i'th coordinate value for the entity obtained by such ultrafinite combinations. 
Such entities are called {{\it intermediate propertons}}. Let $m_0\geq 0$ be the 
N-world proper mass for an assumed elementary particle denoted by $\rm 
k^\prime.$ If $m_0 = 0,$ then let $\lambda = 1.$ Otherwise$,$ from Theorem 11.1.1$,$ we know that there is a $\lambda \in\hypernat
$ such that $\lambda /(10^\omega) \in \monad {m_0},$  where $\omega \in 
\nat_\infty$ and since $m_0 \not=0,$ $\lambda \in \nat_\infty.$ 
Consequently$,$ for $b_n = 10^{-\omega},$ the *-finite sum
$$\sum_{n=1}^\lambda b_n=\sum_{n=1}^\lambda {{1}\over{10^\omega}} = {{\lambda}\over {10^\omega}}\leqno 
(11.1.5)$$
has the property that $\st {\sum_{n=1}^\lambda 1/(10^\omega)} = m_0.$ (Note the special summation notation for a constant summand.) The 
{standard part operator st is an important external operator that is a continuous 
[11] NSP-world process that yields N-world effects.} The appropriate 
interpretation is that\parm 
{\leftskip 0.5in \rightskip 0.5in \noindent {\it ultrafinite combinations of 
ultra-propertons yield an intermediate properton that$,$ after application 
of the standard part operator$,$ has the same effect as an elementary particle 
with proper mass $m_0.$}\par}\parm
An additional relevant idea deals with the interpretation that the *-finite 
set $\rm G_A^{\prime\prime}$ exists at$,$ say$,$ nonstandard time$,$ and that such a 
set is manifested at standard time when the operator st is applied. The 
standard part operator is one of those external operators that can be 
indirectly detected by the presence of elementary particles with proper mass 
$m_0.$ \pars
The above discussion of the creation of {intermediate propertons} yields a 
possible manner in which ultra-propertons are combined within the NSP-world
to yield appropriate energy or mass coordinates for the multifaceted 
propertons. But is there an indication that all standard world physical 
qualities that are denoted by qualitative measures begin as infinitesimals?   
\pars 
Consider the infinitesimal methods used to obtain such things as the charge on 
a sphere$,$ charge density and the like. In all such cases$,$ it is assumed that 
charge can be infinitesimalized. In 1972$,$ it was shown how a classical theory 
for the electron$,$ when infinitesimalized$,$ leads  to the point charge concept 
of quantum field theory and 
then how the *-finite many body problem produced the quasi-particle. [15] 
Although this method is not the same as the more general and less ad hoc 
properton approach$,$ it does present a procedure that leads to an 
infinitesimal charge density and then$,$  in a very ad hoc manner$,$ it is assumed 
that there are objects that when *-finitely combined together entail a real 
charge and charge density. Further$,$ it is the highly successful use of the 
modeling methods of infinitesimal calculus over hundreds of years that has
lead to our additional presumption that all coordinates of the basic 
sequential properton representation are a $\pm$ fixed infinitesimal. \pars
In order to retain the general independence of the coordinate representation$,$ 
{{\it independent *-finite coordinate summation}} is allowed$,$ recalling that 
such objects are to be utilized to construct many possible universes. [This 
is the same idea as *-finitely repeated simple affine or linear transformations.] Thus$,$ 
distinct from coordinatewise addition$,$ *-finitely many such sequences can be 
added together by means of a fixed coordinate operation in the following sense.
Let $\{a_i\}$ represent an ultra-properton. Fix the coordinate $j,$ then the 
sequence $\{c_i\},\ c_i = a_i ,\ i \not= j$ and $c_j = 2a_j$ forms an 
intermediate properton.  As will be shown$,$ it is only after the formation of 
such intermediate propertons that the customary coordinatewise 
addition is allowed and this yields$,$ after the standard part operator is 
applied$,$ representations for elementary particles. Hence$,$ from our previous 
example$,$ we have that ultrafinite combinations of ultra-propertons yield 
propertons with ``proper mass'' $\lambda/(10^\omega) \approx m_0$ while all 
other coordinates remain as $\pm 10^\omega.$ This
physical-like process is not a speculative ad hoc construct$,$ but$,$ rather$,$
it is modeled  after what occurs in our observable natural world.
Intuitively$,$ this  type of summation is modeled after the process of
inserting finitely  many pieces of information (mail) into a single
``postal box$,$'' where  these boxes are found in rectangular arrays in post
offices throughout the  world.\pars 
Now other ultra-propertons are ultrafinitely combined and yield for a specific 
coordinate the $\pm$ unite charge or$,$ if quarks exist$,$ other N-world charges$,$ 
while all other coordinates remain fixed as $\pm 1/(10^\omega)$$,$ etc. 
Rationally$,$ how can one conceive of a combination of these intermediate 
propertons$,$ a combination that will produce entities that can be 
characterized in a standard particle or wave language?\pars
Recall that a finite summation is a *-finite summation within the NSP-world. 
Therefore$,$ a finite combination of intermediate propertons is an allowed 
internal process. [Note that external processes are always allowed but with 
respect to our interpretation procedures we always have direct or indirect knowledge 
relative to application of internal processes. Only for very special 
and reasonable external processes do we have direct or indirect knowledge that 
they have been applied.]  Let $\zeta_i \in \monad 0,\ i= 1,\ldots,n.$ Then 
$\zeta_1 + \cdots \zeta_n \in \monad 0.$ The {final stage in properton 
formation for our universe} --- the final stage in particle or wave 
substratum formation --- would be {\it finite coordinatewise} summation of 
finitely many intermediate propertons. This presupposes that the N-world 
environment is characterized by but finitely many qualities that can be 
numerically characterized. This produces the following type of coordinate 
representation for a specific coordinate $j$ after $n$ summations with $n$ 
other intermediate propertons that have only infinitesimals in the $j$ 
coordinate 
position.
$$\sum_{i=1}^{\lambda}(1/(10^\omega) + \sum_{i=1}^n\zeta_i. \leqno (11.1.6)$$
Assuming $\lambda$ is one of those members of $\nat_\infty$ or equal to 1 as 
used in (11.1.5)$,$ 
then the standard part operator can now be applied to (11.1.6) and the result 
is the same as $\st {\sum_{i=1}^{\lambda}(1/(10^\omega)}.$\pars

The process outlined in (11.1.6) is then applied to finitely many distinct 
intermediate propertons -- those that characterize an elementary particle. 
The result is a properton each coordinate of which is infinitely close to 
the value of a numerical characterization or an infinitesimal. When the 
standard part operator is applied under the  usual coordinatewise procedure$,$ 
the coordinates are either the specific real coordinatewise characterizations 
or zero. Therefore$,$ N-world formation of particles$,$ the dense substratum 
field$,$ or even gross matter may be accomplished by a ultrafinite combination
of ultra-propertons that leads to the intermediate 
properton; followed by finite combinations of intermediate propertons that 
produce the N-world objects. Please note$,$ however$,$ that prior to application 
of the standard part operator such propertons retain infinitesimal nonzero
coordinate characterizations in other noncharacterizing positions. (See note (1) on 116.)\pars  
                          
We must always keep in mind the hypothesis law [9] and {avoid unwarranted 
speculation}. We do not speculate whether or not the formed particles have 
{point-like} or  ``{spread out}'' properties within our space-time environment. 
These additional concepts may be pure {catalyst} type statements within some 
standard N-world theory and could have no significance for either the N-world 
or NSP-world. \pars
With respect to field effects$,$ the cardinality of the set of all 
ultra-propertons clearly implies that there can be ultrafinite combinations   
of ultra-propertons  ``located'' at every   ``point'' of any 
finite dimensional continuum. Thus the field effects yielded by propertons 
may present a {completely dense continuum} type of pattern within the N-world 
environment although from the monadic viewpoint this is not necessarily how 
they   ``appear'' within the NSP-world.\pars
There are many scenarios for {quantum transitions} if such occur in objective 
reality. The simplest is a re-ultrafinite combination of the ultra-propertons 
present within the different objects. However$,$ it is also possible that this 
is not the case and$,$ depending upon the preparation or scenario$,$ the 
so-called   ``conservation'' laws do not hold in the N-world. \pars
As an example$,$ the {neutrino could be a complete fiction}$,$ only endorsed as a 
type of catalyst to force certain laws to hold under a particular scenario. 
Consider the set of sentences \parm
\line{\hfil${\rm G_B = \{An\sp elementary\sp particle\sp}
\r k^\prime(\r i^\prime,\r j^\prime)\sp{\rm with\sp}$\hfil}\pars
\line{\hfil${\rm total\sp energy\sp}\r c^\prime{+}\r n^\prime.\mid ((i,j,n) \in $\hfil}\pars
\line{\hfil \qquad$\nat^+ \times 
\nat^+ \times \nat^+)\land (1\leq k\leq m )\}.$\hfil (11.1.7)}\parm
It is claimed by many individuals that such objects as being described in $\rm 
G_B$ exist in objective reality. Indeed$,$ certain well-known scenarios for a 
possible cosmology require$,$ at least$,$ one  ``particle'' to be characterized by 
such a collection $ \rm G_B.$ By the usual method$,$ these statements are *-
transferred to \parm 
\line{\hfil${\rm G_B^\prime = \{An\sp elementary\sp particle\sp}
\r k^\prime(\r i^\prime,\r j^\prime)\sp{\rm with\sp}$\hfil}\pars
\line{\hfil${\rm total\sp energy\sp}\r c^\prime{+}\r n^\prime.\mid ((i,j,n) \in $\hfil}\pars
\line{\hfil \qquad$\hypernat^+ \times 
\hypernat^+ \times \hypernat^+)\land (1\leq k\leq m )\}.$\hfil (11.1.8)}\parm
\noindent Hence$,$ letting $n \in \nat_\infty$ then various   {``infinite'' 
NSP-world energies} emerge from our procedures. With respect to the total 
energy coordinate(s)$,$ ultra-propertons may also be ultrafinitely combined to 
produce such possibilities. Let $\lambda = 10^{2\omega}$ [ resp. $\lambda = 
\omega^2$] and $\omega \in \nat_\infty.$ Then 
$$\sum_{n=1}^\lambda {{1}\over{10^\omega}} =10^\omega \  {\rm 
[resp.\ }\sum_{n=1}^\lambda {{1}\over {\omega}} = \omega] \in 
\nat_\infty.\eqno (11.1.9)$$\par
Of course$,$ these numerical characterizations are external to the N-world. 
Various distinct {``infinite'' qualities can exist} rationally in the NSP-world 
without altering our interpretation techniques. The behavior of the infinite 
hypernatural numbers is very interesting when considered as a model for 
NSP-world behavior. A {transfer of finite energy}$,$ momentum and$,$ indeed$,$ all other 
N-world characterizing quantities$,$ back and forth$,$ between these two worlds is 
clearly possible without destroying  {NSP-world infinite conservation 
concepts}. 
\pars
Further$,$ observe that  various intermediate propertons carrying nearstandard 
coordinate values could be present at nearstandard space-time coordinates$,$ and 
application of the continuous and external standard part operator would 
produce an apparent not conserved N-world effect. These concepts will be 
considered anew when we discuss the Bell inequality.\pars

Previously$,$ ultrawords were obtained by application of certain concurrent 
relations. Actually$,$ basic ultrawords exist in any elementary nonstandard 
superstructure model$,$ as will now be established for the general paradigm. 
\pars
Referring back to $\rm G_A$ equation (11.1.1)$,$ for some fixed $k,\ 1\leq k\leq 
m,$ let\break
\vfil
\eject 
\noindent $\r h_k \colon \nat^+ \times \nat^+ \times \nat^+ \to  \rm G_A$ be 
defined as follows: $\r h_k(i,j,n) = $\break ${\rm An\sp elementary\sp particle\sp} \r k^\prime 
(\r i^\prime,\r j^\prime){\rm \sp with\sp total\sp energy\sp} 
\r c^\prime{+} 1/(\r n^\prime).$ Since the set $F(\nat^+ \times \nat^+ \times 
\nat^+)$ is denumerable$,$ there exists a bijection $\r H\colon \nat \to$\break$ 
F(\nat^+ \times \nat^+ \times \nat^+).$ For each $1\leq \lambda \in \nat$ and fixed 
$i,\ n \in \nat^+,$ 
let 
${\rm G_A(\lambda)} = \{{\rm An\sp elementary\sp particle\sp} \r k^\prime 
(\r i^\prime,\r j^\prime) {\rm \sp with\sp total\sp energy\sp} 
\r c^\prime{+}1/(\r n^\prime).\mid (1\leq j'\leq \lambda) \land (j' \in \nat^+) \}.$ 
Let $p \in \nat.$ If $\vert \r H(p) \vert \geq 2,$ define finite $\r M(\r h_k[\r H(p)]) = 
\{{\rm A_1\sp and\sp A_2\sp and\sp\cdots \sp and\sp} \r A_m\},$ where ${\r A_j} \in 
\r h_k[\r H(p)], m = \vert\r H(p)\vert.$ If $\vert \r H(p) \vert \leq 1,$ then define
$\r M(\r h_k[\r H(p)]) = \emptyset.$ Let $\r M^0 = \bigcup\{\r M(\r h_k[\r H(p)])\mid p \in \nat 
\}.$ Please note that the $\r k^\prime$ represents the  ``type'' or name of the 
elementary particle$,$ assuming that only finitely many different types exist,
$\r i^\prime$ is reserved for other purposes$,$ and the $\r j ^\prime$ the number of such elementary particles of 
type $\r k^\prime.$ \parm
{\bf Theorem 11.1.3} {\sl For any $i,\ n,\ \lambda \in \hypernat^+,$ such that $2 
\leq \lambda,$ there exists $w \in \Hyper {\bf M^0 - G_A},\ \Hyper{\bf 
G_A(\lambda)} \subset \Hyper{\b S(\{w\})}$ and if $A \in \Hyper {\bf G_A} - 
\Hyper{\bf G_A(\lambda)},$ 
then $A \notin \Hyper{\b S(\{w\})}.$} \pars
Proof. Let $i,\ j,\ \lambda \in \nat^+$ and $ 2 \leq \lambda.$ Then there 
exists some $ r\in \nat$ such that $\r h_k[\r H(r)] = {\rm G_A(\lambda)}.$ From 
the construction of $\r M^0,$ there exists some $r^\prime \in \nat$ such that 
$\r w(r^\prime) = {\rm An\sp elementary\sp particle\sp k^\prime(i^\prime,1^\prime) \sp with\sp total\sp energy \sp}\break{\rm 
c^\prime{+}1/(n^\prime). 
\sp and\sp An\sp elementary\sp 
particle\sp k^\prime(i^\prime,2^\prime) \sp with\sp total\sp energy 
\sp} \break {\rm c^\prime{+}1/(n^\prime). \sp and\sp \cdots\sp and\sp An\sp elementary\sp 
particle\sp k^\prime(i^\prime,\lambda^\prime)\sp\break with\sp total\sp}\break{\rm 
energy \sp c^\prime{+}1/(n^\prime). } \in \r M[\r h_k[\r H(r^\prime)].$ Note that 
$\r w(r^\prime) \notin {\rm G_A},\ \r h_k[\r H(r^\prime)] \subset \r S(\{\r 
w(r^\prime)\})$ and if ${\rm A} \in {\rm G_A} - \r h_k[\r H(r)],$ then $\r A \notin \r 
S(\{\r w(r^\prime)\}).$ The result follows by our embedding and *-transfer. \qed
The ultrawords utilized to generate various propertons$,$ whether obtained as 
in Theorem 11.1.3 or by concurrent relations$,$ are called {{\it ultramixtures}} 
due to their applications. The ultrafinite choice operator $\bf C_1$ can 
select them$,$ prior to application of $\Hyper {\b S}.$ Moreover$,$ application of 
the ultrafinite combination operator entails a specific intermediate 
properton with the appropriate nearstandard coordinate characterizations. 
Please notice that the same type of sentence collections may be employed to 
infinitesimalize all other quantities$,$ although the sentences need not have 
meaning for certain popular N-world theories. Simply because substitution of 
the word ``charge'' for  ``energy'' in the above sentences  $\rm G_A$ does not 
yield a particular modern theory description$,$ it does yield the infinitesimal 
charge concept prevalent in many older classical theories. \pars
Using such altered $\rm G_A$ statements$,$ one shows that there does exist 
ultramixtures $w_i$ for each intermediate properton and$,$ thus$,$ a single
ultimate ultramixture $w$ such that $\Hyper {\b S}(\{w_i\}) \subset 
\Hyper {\b S}(\{w\}).$ Each elementary particle may$,$ thus$,$ be assumed to 
originate from $w$ through application of the ultralogic $\Hyper {\b S}.$\pars

Recall that if a standard $A \subset \real$ is infinite$,$ then it is external$,$ 
and if $B$ is internal$,$ $A\subset B,$ then $B \not= A.$ Therefore$,$ there 
exists some $\eta \in B$ such that $\eta \notin A.$ This simple fact yields 
many significant nonstandard results. For example$,$ as the next theorem shows$,$ 
if $\eta \in \nat_\infty,$ then there exists some $\lambda \in \nat_\infty$ 
such that $10^{2\lambda} < \eta.$\parm
\vfil
\eject
{\bf Theorem 11.1.4} {\sl Let $f\colon \nat \to \nat$ and $f[\nat]$ be 
infinite. If $\eta \in \nat_\infty,$ then there exists some $\lambda \in 
\nat_\infty$ such that $\hyper f(\lambda) < \eta.$}\pars
Proof. For $\eta \in \nat_\infty,$ consider the nonempty internal set $B = \{\hyper 
f(x)\mid (\hyper f(x) < \eta) \land (x 
\in \hypernat) \}\subset \hyper f[\hypernat].$ 
Let $n \in \nat.$ Then $f(n) \in \nat$ and $f(n) < \eta$ imply that $(^\sigma 
f)[\nat] = {^\sigma (f[\nat])} = f[\nat]  \subset B.$ Since 
$f[\nat]$ is infinite$,$ it is external. Thus $\hyper f\colon \hypernat \to \hypernat$ implies 
that there exists some $\lambda \in \hypernat$ such that $\hyper f(\lambda) \in B - f[\nat].$ 
However, $\hyper f$ is a function. Hence, $\lambda \in \nat_\infty$ and $\hyper f(\lambda) < \eta.$ 
\qed
Theorem 11.1.4 has many applications and can be extended to other functions not 
just those with domain and codomain $\nat,$ and other $B$ type relations.\par
\bigskip
\leftline{\bf 11.2 Ultraenergetic Propertons (Subparticles).} 
\medskip
There is a possibility that propertons can have additional and unusual 
properties when they are generated by statements such as $\rm G_B.$ With 
respect to the translated $\rm G_B^\prime$ statements$,$ we have coined the term
{{\it ultraenergetic}} to discuss propertons that have various infinite 
energies. I again note that such ultraenergetic propertons may be considered 
as under the control of our previously discussed ultralogics

and {{\it 
ultranatural choice operators}} (i.e. hyperfinite choice). Further$,$ it is 
possible to place these ultraenergetic 
propertons into pools of infinite energy that are mathematically termed as
``{galaxies}'' and that have interesting mathematical properties. However$,$ these 
properties will not be discussed in this present book.\pars
{\bf If} our universe or any portion of it began or exists at this present 
epoch 
in a {state of ``infinite'' energy}$,$ then the ultraenergetic propertons 
could play a critical role. [I point out that general developmental paradigms 
indicate$,$ as will be shown$,$ that the actual state of affairs for any 
beginnings of our universe cannot be known by the present methods of the 
scientific method.] Now$,$ any such {singularity} that might exist 
cosmologically$,$ even in our local environment$,$ may owe its existence to 
various ultralogically generated ultraenergetic propertons. Of course$,$ this 
is pure speculation$,$ but these NSP-world alternative explanations for {\it 
assumed} quantum physical phenomena yield indirect evidence for the acceptance 
of the NSP-world model. \pars
Quantum mechanics has now become highly {positivistic} in character although 
certain previous states of affairs have been partially accepted. This 
important possibility was stated by Bernard d'Espagnet with respect to one of 
our 
 preliminary investigations --- the experimental disproof of the {Bell 
inequality} and the local variable concept [17] --- that ``seems to imply that 
in some sense all of the objects [particles or aggregates] constitute an 
indivisible whole'' [2]. One aspect of the MA-model$,$ (11.4.5) of section 11.4$,$ 
can be used as an aid to model this statement.\pars
The ideas developed within our theory of developmental paradigms do not 
contradict d'Espagnet's {(weak) definition of realism}. He simply requires that 
if we can describe a relation between physical entities produced by some 
experimental process$,$ a relation that is not observed and$,$ thus$,$ not described 
prior to the experiment$,$ then their must be a cause that has produced this new 
relation. I don't believe that it is necessary$,$ under his definition$,$ 
that this cause be describable. \pars
In May 1984$,$ this author became aware of the {Bell inequality} and d'Espagnat's 
discussion of local realism [3]. In particular$,$ we discovered that d'Espagnat 
may$,$ to some degree$,$ embraced our the statements that appear in (11.4.5) 
section 11.4 as 
an explanation for this experimental disproof. ``Perhaps in such a world the 
concept of an independent existing reality can retain  some meaning$,$ but it 
will be altered and one remote from everyday experience'' [4]. But is there a 
NSP-world cause$,$ indeed$,$ a mechanism that for such behavior?
\pars
There are many scenarios as to how  {``instantaneous'' informational signals} 
may 
be transmitted within the NSP-world without violating {Einstein separability} 
in the N-world (i.e. no influence of any kind within the N-world can propagate 
faster than the speed of light [5].) The basic N-world interpretation for any 
verified effect of the Special Theory would imply that the only reason for 
Einstein separability is a relation between those propertons that create the 
N-world and the NSEM field propertons. But ultraenergetic propertons are 
not of either of these types and need not interact with the NSEM field for 
many reasons. The most obvious is that the {NSEM field is not dense from the 
NSP-world viewpoint but is scattered}. Obviously ultraenergetic propertons 
may be used for this purpose. Recall that we should be very careful when 
speculating about the NSP-world due to the difficulty of describing refined 
behavior. However$,$ this should not completely restrain us$,$ especially when the 
general paradigm method states that such things as these  exist logically. 
\pars
It is possible to describe a mechanism and a possible new type of properton 
that can send N-world {instantaneous informational signals} between all standard 
material particles$,$ field objects or aggregates and {not violate N-world 
Einstein separability}. One possibility is that these influences would be 
imparted by means of {independent coordinate summation} to the propertons 
that 
comprise these objects and yet in doing so these new entities could not be 
humanly detected$,$ {{not detectable}} except as far as the instantaneous state change indicates,
 since the total energy (in this case classical kinetic) 
utilized by this NSP-world mechanism would be infinitesimal. If it is an 
instantaneous energy change$,$ then$,$ as will be shown$,$ only that specific energy 
change would appear in the N-world. \pars  
From the methods employed to construct propertons$,$ it is immediately clear 
that there exists a very ``large'' quantity of propertons that are not used for 
standard particle and field effect construction. This can be seen by allowing 
the i'th symbol is such statements as $\rm G_A^\prime$ to vary from 1 to some 
value in $\nat_\infty.$ The cardinality of such collections of statements 
would be great than or equal to ${\vert {\cal M}\vert}^+.$ We simply pass 
this external cardinality statement to the propertons being described. \pars
 It 
is a basic tenet of infinitesimal reasoning that without further justification 
the only properties that we should associate with such unutilized objects are 
of the simplest classical type. The logic of particle physics allows us to 
logically accept the existence of such propertons without any additional 
justification. Let $\lambda \in \nat_\infty.$ Then $(1/\lambda)^4 \in \monad 
0.$ 
Let a pure NSP-world properton$,$ not one used to construct a universe$,$ have 
mass coordinate of the value $m =(1/\lambda)^4.$ Call this properton P and 
have P increase its velocity over a finite NSP-time interval from zero to 
$\lambda.$ Propertons the attain such velocities are called 
{{\it ultrafast propertons}}.\pars
Extending the classical idea of kinetic energy to the NSP-world it follows 
that the (kinetic) energy attained by this properton$,$ when it reaches its 
final velocity$,$ is $(1/2)(1/\lambda)^2 \in \monad 0.$ Suppose that there is a 
continuum or less of positions within our universe. Since $c < \lambda,$ the 
kinetic energy used to accelerate enough of these propertons to the $\lambda$ 
velocity so that they could effect every position in our universe would be 
less than  
$(1/2)(1/\lambda) \in \monad 0,$ assuming the energy is additive in the 
NSP-world. \pars 
Each of these ultrafast propertons besides altering some 
other specific coordinate would also add its total kinetic energy  
to the intermediate properton since the new intermediate properton simply 
includes this new one. 
But then suppose that our universe has existed for less than or equal to a 
continuum of time. Then since $2c < 
\lambda,$
once again the amount of energy 
that would be added to our universe over such a time period$,$ if each of 
these informational propertons combined with one member of an 
intermediate$,$ would be infinitesimal. All the state alterations give 
the N-world appearance of being {instantaneously obtained} although the 
existence of the 
$G$ function of Theorem 7.5.1 clearly states that in the NSP-world such 
alterations are actually hypercontinuous and hypersmooth. \pars
What if the state change itself depends upon the velocity of such an ultrafast 
properton? We use kinetic energy as an example. Say the change is in the 
kinetic energy coordinate in the standard amount of $h.$ Then all one needs to 
consider is an ultrafast properton with infinitesimal mass 
$m = 2h(1/\lambda)^2.$ Moving with a velocity of $\lambda,$ such an ultrafast 
properton has the requisite kinetic energy. Things can clearly be arranged 
so that all other coordinates of such ultrafast propertons are 
infinitesimal. Independent coordinate summation for any finite number of 
alterations will leave all other nonaltered coordinates of the  
intermediate properton infinitely close to the original values for the 
alteration is but obtained by the addition of a finite number of new 
propertons to the collection. \pars
We acknowledge that the N-world inner coordinate relations have been used to 
obtain these alterations. This need not be the way it could be done. 
Can we describe the method of capture and other sorts of behavior? 
Probably too much 
has already been described in the language of this book. One should not forget 
that descriptions may exist for such NSP-world behavior but not in a readable 
language. \pars
The state of affairs described above lends credence to d'Espagnate's 
explanation of why the Bell inequality is violated and gives further evidence 
for the acceptance of the NSP-world model.  
``The basic law that signals 
cannot travel faster than light is demoted from a property of external 
[N-world] reality to a feature of mere communicatable human experience. 
....the concept of an independent or external reality can still be retained as 
a possible explanation of {\it observed} regularities in experiments. It is 
necessary$,$ however$,$ that the violation of Einstein separability be included as 
a property$,$ albeit a well-hidden and counterintuitive property....''[6]  \par                 
\bigskip
\leftline{{\bf 11.3 More on Propertons (Subparticles),}}\par
\medskip
The general process for construction of all N-world fundamental entities from 
propertons can be improved upon or achieved in an alternate fashion. One of 
the basic assumptions of subatomic physics is that in the Natural-world two 
fundamental subatomic objects$,$ such as two electrons$,$ cannot be differentiated 
one from another by any of its Natural-world properties. One of the 
conclusions of what comes next is that in the NSP-world this need not be the 
case. Of course$,$ this can also be considered but an auxiliary result and need 
have no applications.
At a 
particular instant of (universal) time$,$ it is  possible to associate with 
each entity a distinct ``name'' or identifier through properton 
construction. This is done through application of the properton naming 
coordinate  $a_1.$ As will be shown$,$ the concept of independent *-finite 
coordinate summation followed by n-tuple vector addition can be accomplished 
by means of a simple linear transformation. However$,$ by doing so$,$ the concept 
of the *-finite combinations or the gathering together of propertons  
as a NSP-world physical-like process is suppressed. Further$,$ a simple method 
to identify each N-world entity or Natural system is not apparent. Thus$,$ we 
first keep the above two processes so as to adjoin to each entity constructed 
an appropriate identifier. \pars

A standard properton is modeled by a finite collection of numerical or coded descriptive physical characteristics. These characteristics are represented by coordinates within n-tuples. Other identifiers can also be included as specific coordinates. Included within these propertons are those of the following special type. \pars 

Informally, consider the denumerable set all prime numbers P, a bijection $\r h \colon \nat' \to \r P,\ \nat' = \nat - \{0\}$, and the  sequence $\r g \colon \nat  \to \b Q,$ the set of rational numbers, where $g(n) = 1/10^n.$ Let $\omega \in \hypernat' - \nat' = \hypernat - \nat = \nat_\infty.$ Since $g \to 0,\ n \to \infty,$ then $\hyper g(\omega) = 1/10^\omega \in \monad {0},$ the set of all infinitesimals. \parm

{\bf Definition 11.3.1, Ultra-propertons.} Let even $K >2, K \in \nat$ and $f\colon [1,K] \to \{1/10^\omega\}$. Then ${\cal C} = \{(\Hyper h(i), 1, -f(1), f(2),\ldots,-f(K-1), f(K))\mid i \in \hyper \nat' \}.$ Each member of the set $C$ represents an ultra-properton. \parm

For Definition 11.3.1, it is assumed that there is no more than $K$ physical or physical-like numerical or coded descriptive characteristics for the any elementary entity. \parm

\noindent {\bf Theorem 11.3.1.} {\sl Consider any nonempty internal $D$ and $A \subset D$ such that $|A| < |{\cal M}_1|^+$. Then there exists a hyperfinite $B_A$ such that $A \subset B_A \subset D.$ }\pars
Proof. Let $\cal F$ be the finite power set operator. That is for any set $X, \ {\cal F}(X)$ is the set of all finite subsets of 
$X$, where a set $Y$ is finite if it is empty or there exists an $n \in \nat'= \nat - \{0\}$ and a bijection $f'\colon [1,n] \to Y.$ In our structure, there is a least $n \in \nat',$ such that internal $D \in \hyper X_n.$ If $x \in A,$ then $x\in D, \ x \in \hyper X_{n-1},$  $\Hyper {\cal F}(D) \in \Hyper X_{n+2}$. If $y \in \Hyper {\cal F}(D),$ then $y \in \hyper X_{n+1}.$  \pars

Consider the internal binary relation $C = \{(x,y) \mid (x \in y)\land (y \in \hyper X_{n+1})\land ( x \in D)\land(x \in \hyper X_{n-1}) \land (y\in \Hyper {\cal F}(D))\}.$ Let $\{(x_1,y_1),\ldots,(x_m,y_m)\} \subset C.$ Then  $y' = y_1 \cup \cdots \cup y_m$ is an internal subset of $\Hyper {\cal F}(D).$ Since the domain of $C$ is  $D$, $A \subset D$  and $|A| < |{\cal M}_1|^+,$ then by saturation there exists a $B_A \in \Hyper {\cal F}(D)$ such that $A \subset B_A \in \Hyper {\cal F}(D)$ and $B_A \subset D$. This complete the proof. \qed

\noindent {\bf Corollary 11.3.1.1} {\sl Consider any $\Hyper E$ and $A \subset \Hyper E$ such that $|A| < |{\cal M}_1|^+$. Then there exists a hyperfinite $B_A$ such that $A \subset B_A \subset \Hyper E.$ }\parm

\noindent {\bf Theorem 11.3.2} {\sl Consider any nonempty hyperfinite $A \subset \hypernat'$. Then there exists a $\gamma \in \hypernat$ such that $A \subset [1, \gamma]$.} \pars
Proof. Every nonempty finite subset $F$ of $\nat'$ has a greatest member $M_F \in \nat$. That is if $x \in F,$ then $x \in [1,M_F].$ By *transfer, $A$ has a *greatest member $\gamma \in \hypernat'$ such that if $x \in A,$ then $x \in [1,\gamma]$ \qed
Let $r_1 \in \real.$ By Theorem 11.1.1 in Herrmann (1979-93), there is a $\lambda_1 \in \nat_\infty$ such that $\lambda_1/10^\omega \in \monad {|r|}.$ Hence, $\st {(\lambda_1/10^\omega)} = |r|.$ Then there are $K$, $\lambda_i, \ i \in [1,K]$ that yield the $K$ characteristics. For an elementary entity $\r e_j$, some characteristics can be 0, meaning that the measure has value 0. Throughout the combining processes, if a coordinate retains its infinitesimal value $\pm 1/10^\omega,$ this indicates that the characteristic has no meaning for $\r e_j.$ In order to indicate these differences, any characteristic that has measure 0 is obtained from a combination of two ultra-propertons. The standard part physical realization operator $\tt St$ is only applied to coordinates of the intermediate properton representations with the form $\pm \lambda/10^\omega,$ where $\lambda \geq 2.$\pars

There are other characteristics such as spin, where the 0 takes on a different meaning. However, such coding is rather arbitrary and can be replaced with non-zero numbers or non-zero codings for the characteristics so as to not confuse them with a 0 measurement. For the needed intermediate properton $e_1$, with a third coordinate  charcteristic  under independent coordinate addition, the set of ultra-propertons $\{(\Hyper h(i),1,-1/10^\omega, \ldots, 1/10^\omega)\mid i \in [1,\lambda_1]\}$ is employed. Hence, the first intermediate properton is $(\Pi_{1}^{\lambda_1}, \lambda_1, -\lambda_1/10^\omega,1/10^\omega, \ldots, 1/10^\omega).$ For a forth coordinate intermediate properton for value $r_2,$ consider $\{(\Hyper h(i),  1,-1/10^\omega, \lambda_2/10^\omega,\ldots, 1/10^\omega)\mid i \in [\lambda_1 +1,\lambda_1 +\lambda_2]\}.$ Continue these definition for each member of $[1,K].$ Thus the entire collection of ultra-propertons used to obtain one of the $e_1$ entities is $\lambda_1 + \cdots +\lambda_K = \delta_1 \in \nat_\infty.$ \pars

It is assumed that there are a nonempty countable (i.e non-zero finite or denumerable) collection of $\{e_i\}$ needed. Thus there is a non-zero finite or denumerable set $\{\delta_i\}$ and in the finite case, consider $\sum \delta_i \in \nat_\infty.$ Next consider $\{\delta_i\mid i \in \nat'\}.$ Then $\{\delta_i\mid i \in \nat'\} \subset \hypernat.$ The $|\{\delta_i\mid i \in \nat'\}| < |{\cal M}_1|^+$. Hence, there is a $\gamma_1 \in \nat_\infty$ such that $\{\delta_i\mid i \in \nat'\} \subset [1,\gamma_1]$ by application of Corollary 11.3.1.1 and Theorem 11.3.2. Thus, in both cases, there is a $\Gamma_1 \in \nat_\infty$ such that $\{\delta_i\} \subset [1,\Gamma_1]$. This shows that there are ``enough'' ultra-propertons to produce the set $\{e_i\}.$ For another type of elementary particle, simply repeat this for the identifiers $h(i),\ i > \Gamma_1].$ Then continue by induction. \pars

For this application, it appears unnecessary to consider more than $H$, where  $1 \leq H \in \nat$, different types of elementary entities. The set of ultra-propertons $\{(\Hyper h(i), 1,-1/10^\omega, \ldots, 1/10^\omega)\mid i \in \hypernat'\}= {\cal C}$ is an internal set and as such the hyperfinite operator $\Hyper {\cal F}$ is defined for it. For properton generation, a universe can be considered as a collection of physical-systems. Hence application of a finite iteration $\hyper {\cal F}^i$ to $\cal C$ yields $\bigcup \{\hyper {\cal F}^i({\cal C})\mid (0\leq i \leq n)\land (i \in \nat)\},$ an internal collection that is sufficient to generate the physical-systems for any of the presently considered cosmologies. To accommodate the formation of the physical-like systems, internal $X$ that is disjoint from $\bigcup \{\hyper {\cal F}^i({\cal C})\mid (0\leq i \leq n)\land (i \in \nat)\}$ is adjoined to $\bigcup \{\hyper {\cal F}^i({\cal C})\mid (0\leq i \leq n)\land (i \in \nat)\}$. \pars

Relative to the GGU-model and generation of a universe, a hyperfinite $\Hyper {\b I}^q(i,j)$ yields a universe-wide frozen-frame. (See the appendix prior to the symbols page.) Each instruction $x \in \Hyper {\b I}^q(i,j)$, yields a physical or physical-like system. The physical-systems are disjoint. Each collection of ultra-propertons that yields a specififc physical-system is distinct from the set of ultra-propertons that yields any other physical-system. Hence, each physical-system within a universe-wide frozen-frame has a distinct identifier via the collection of all of the identifiers for the ultra-propertons or the intermediate propertons employed to produce the physical-system.\par

\bigskip
\leftline{\bf 11.4 MA-Model.} 
\medskip
In this section$,$ we look back and gather together various observations 
relative to formal theorems that yield the concept I have described as the 
Metaphoric-Anamorphosis (i.e. MA) model. The different types of developmental paradigms 
that can be selected by ultrafinite (i.e. ultranatural) choice and a few of our previous results 
leads immediately to the following logically acceptable possibilities. [Of 
course$,$ as is the case with all mathematical modeling$,$ simply because a 
possibility exists it need not be utilized to describe an actual scenario and 
if it is used$,$ then it need not be an objectively real description.]  
\parm
{\leftskip=0.5in \rightskip =0.5in \noindent (11.4.1) {\it Entire microscopic$,$ 
macroscopic or large scale natural systems can apparently appear or disappear or 
be physically altered suddenly.}\par}\pars
{\leftskip=0.5in \rightskip =0.5in \noindent (11.4.2) {\it Theorem 7.3.1 
shows clearly that the suddenly concept in (11.4.1) is justifiable. 
Ultralogics$,$ ultrawords$,$ the 
intermediate properton  
and the 
ultrafast properton concept are possible mechanisms that can yield the 
behavior described in (11.4.1).}\par}\pars
{\leftskip=0.5in \rightskip =0.5in \noindent (11.4.3) {\it All such 
alterations may occur in an ultracontinuous manner.}\par}\pars
{\leftskip=0.5in \rightskip =0.5in \noindent (11.4.4) {\it None of these 
NSP-world concepts are related to the notion of hidden variables.}\par}\pars
{\leftskip=0.5in \rightskip =0.5in \noindent (11.4.5) {\it Any numerical 
quantity associated with any elementary particle$,$ field effect or aggregate is 
associable with every numerical quantity associated with every other 
elementary particle or aggregate by means of hypercontinuous$,$ hyperuniform and 
hypersmooth pure NSP-world functions. These functions may be interpreted as 
representing the IUN-altering process of utilizing ultrafinite composition 
(i.e. {ultranatural composition}) in order to  ``change'' any elementary 
particle$,$ field effect or aggregate into any type of elementary particle$,$ 
field effect or aggregate.}\par}\pars
Statement (11.4.5) is particular significant in that it may be coupled with 
ultralogics and 
ultrafinite choice operators and entails an additional manifestation for the 
possibility that there is no N-world independent existing objective reality. 
{\bf Further$,$ notice that depending upon the space-time neighborhood$,$ statements 
such as (11.4.1) need not be humanly verifiable (i.e. they may be 
undetectable).}\par
\bigskip
\centerline{\bf CHAPTER 11 REFERENCES}
\medskip
\noindent {\bf 1} Beltrametti$,$ E. G.$,$ Enrico$,$ G. and G. Cassinelli$,$ {\it The 
Logic of Quantum Mechanics,} Encyclopedia of Mathematics and Its Applications$,$ 
Vol. 15$,$ Addison-Wesley$,$ Reading$,$ 1981.\pars
\noindent {\bf 2} d'Espagnat$,$ B.$,$ The quantum theory and realism$,$ Scientific 
America$,$ 241(5)(1979)$,$ 177. \pars
\noindent {\bf 3} {\it Ibid.} \pars
\noindent {\bf 4} {\it Ibid.}$,$ 181. \pars
\noindent {\bf 5} {\it Ibid.}$,$ \pars
\noindent {\bf 6} {\it Ibid.}$,$ 180.\pars
\noindent {\bf 7} Feinberg$,$ G.$,$ Possibility of faster-than-light particles$,$ 
Physical Review$,$ 159(5)(1976)$,$ 1089---1105.\pars
\noindent {\bf 8} Hanson$,$ W. C.$,$ The isomorphism property in nonstandard 
analysis and its use in the theory of Banach Spaces$,$ J. of Symbolic Logic$,$ 
39(4)(1974)$,$ 717---731.\pars
\noindent {\bf 9} Herrmann$,$ R. A.$,$ D-world evidence$,$ C.R.S. Quarterly$,$ 
23(2)(l1986)$,$ 47---54.\pars
\noindent {\bf 10} Herrmann$,$ R. A.$,$ The Q-topology$,$ Whyburn type filters and 
the cluster set map$,$ Proceedings Amer. Math. Soc.$,$  59(1)(1975)$,$ 161---
166.\pars
\noindent {\bf 11} Kleene$,$ S. C.$,$ {\it Introduction to  Metamathematics}$,$ D. 
Van Nostrand Co.$,$ Princeton$,$ 1950.\pars
\noindent {\bf 12} Prokhovnik$,$ S. J.$,$ {\it The Logic of Special Relativity}$,$ 
Cambridge University Press$,$ Cambridge$,$ 1967.\pars
\noindent {\bf 13} Stroyan$,$ K. D. and W. A. J. Luxemburg$,$ {\it Introduction to 
the Theory of Infinitesimals}$,$ Academic Press$,$ New York$,$ 1976. \pars
\noindent {\bf 14} Tarski$,$ A.$,$ {\it Logic$,$ Semantics$,$ Metamathematics}$,$ 
Clarendon Press$,$ Oxford$,$ 1969.\pars
\noindent {\bf 15} Thurber$,$ J. K. and J. Katz$,$ Applications of fractional 
powers of delta functions$,$ {\it Victoria Symposium on Nonstandard Analysis},
Springer-Verlag$,$ New York$,$ 1974.\pars
\noindent {\bf 16} Zakon$,$ E.$,$ Remarks on the nonstandard real axis$,$ {\it 
Applications of Model Theory to Algebra$,$ Analysis and Probability}$,$ Holt$,$ 
Rinehart and Winston$,$ New York$,$ 1969. \pars
\noindent {\bf 17} Note that the NSP-world model is not a local hidden 
variable theory.\pars
\noindent {\bf 18} Herrmann$,$ R. A.$,$ Fractals and ultrasmooth microeffects,
J. Math. Physics$,$ 30(4)$,$ April 1989$,$ 805---808.\pars
\noindent {\bf 19} Davis$,$ M.$,$ {\it Applied Nonstandard Analysis,} John Wiley 
\& Sons$,$ New York$,$ 1977. \parm
(1) For propertons, only two possible intrinsic properties for elementary particle formation are here considered. Assuming that there are such things as particles or elementary particles, then they would be differentiated one from the other by their intrinsic properties that are encoded within properton coordinates. When there are particle interactions, these intrinsic properties can be altered or even changed to extrinsic properties. How the alteration from intrinsic to extrinsic occurs probable cannot be known since it most likely is an ultranatural event. For further results on this subject, see http://arxiv.org/abs/quant-ph/9909078\pars
(2) To conceive of propertons properly, quantum theory is viewed as an approximation. Moreover, in terms of physically determined units, the numerical characteristics produced by applications of the standard part operator are  considered as exact. \vfil\eject

\centerline{\bf Ultra-logic-systems}\parm

\noindent{\bf 1. Logic-System Generation for Instructions}\parm

As is customary, the nonstandard model used in all of the articles on the GGU-model is a polysaturated polyenlargements (Lobe and Wolff, 2000; Stroyan and Bayod, 1986). In this paper, $q =1,2,3,4$. These numbers denote the four primitive-time intervals (Herrmann 2006) employed for the GGU-model.  The ultraword approach to generate a universe is replaced with an ultra-logic-system. This is a hyperfinite logic-system where, after application of the extended logic-system algorithm, generates each member of the hyperfinite instruction paradigm $d^q_x$ in the proper $\leq_{d^q_x}$ order such that $\b d_q \subset d^q_x \subset \hyper {\b d_q},$ where $q = 1,2,3,4$ and $x = \lambda, \nu\lambda, \mu\lambda,\nu\gamma\lambda$ respectively. Finally, in this article, the term "subparticle" was previously used. To prevent incorrect mental images as to models for subparticles, the term "properton" replaces the term "subparticle." Without visualizing, a properton is an entity characterized only by a list of properties.\pars

The primitive entity that yields physical reality for any GGU-model generated universe is dense collection of ultra-propertons. 
 When first conceived this author had not investigated quantum field theory and did not base proprtins upon any quantum theoretic approach. All of the GGU-model entities and processes can be considered as existing in a {\bf background universe} or {\bf substratum world}. This world can be considered as a physical-like world, where the rules that govern universe formation are distinct from those processes and rules that govern the development of \underbar{any} physical universe. They are simple rules that only refer to counting. This substratum world is also interpreted philosophically in other ways. \pars

If necessary for a specific physical theory, any continuity requirement is satisfied by the properton field (Herrmann, 1983, 1989). For our universe, a collection of propertons has been shown to be closely associated with relativistic effects (Herrmann, 2003). No other known primitive entities, such as superstrings, will have any effect upon the application of propertons as the primitive entities that generate a universe. The processes used to obtain particles and all other physical entities from ultra-propertons need not correspond to the rules of quantum field theory or any additional rules like how quarks combine to form particles. \pars
For our universe, quantum field theory contains descriptions (rules or instructions) that produce such particles from  immaterial fields. Such fields are quantum  mechanical systems and, when represented, have various degrees of freedom. These are but parameters that contribute to the overall state of the system. For various particles, parameters for physical measures or states are the characterizing features of propertons. The physical appearance  and disappearance of particles are trivial applications of properton processes. For quantum field theory, one has the ``creation'' and ``annihilation'' operators that mathematically yield the same results.\pars

For the GGU-model, quantum theory does not produce steps in a development since the method of production must be universe and physical law independent. For our universe, the development ``satisfies'' the predictions of accepted physical theories. I personally consider quantum theory as mostly a product of human imagination that predicts behavior, behavior that we cannot otherwise comprehend. That is, it is a model that mimics. \pars

The GGU-model can be based upon observable human behavior and the mathematics predicts, for our universe, behavior that satisfies the behavior predicted by accepted physical theories. There is a vast amount of evidence for the predicted GGU-model processes. Whether such processes exist in some sort of reality is a philosophic choice. One can make this choice based upon various factors. One can choose to accept properton existence based upon the same philosophy expressed by those that accept that entities postulated in quantum field and particle theory exist.\pars  

The concept of instructions or rules is generalized to instructions that yield a physical reality from combinations of propertons. They are substratum laws. (So as not to confuse these with physical laws, they are called instructions. Further, in what follows, the events that correspond to each $\r f^q(i,j)$ are denoted by $\r E^q(i,j). $) This does \underbar{not} mean that the rules used in quantum theory (QT) actually yield each $\r E^q(i,j).$ As mentioned, what this signifies is that the QT rules are verified via the production of event sequences that yield our universe. For the GGU-model, the physical realization of each $\r f^q(i,j)$ is not the result of any of these physical theories. These theories are but verified by each realized $\r f^q(i,j)$ and they allow us to predict what behavior occurred in or will occur within other realized $\r f^q(p,k).$ For the GGU-model, the ``instructions'' are rather simple ones that lead to all the characteristics that allow one to identify any material entity for any of the presently known cosmologies. \pars 

 Rather than the $\r f^q(i,j)$ being a general description, one considers instructions or rules $\r I^q(i,j)$ - a nonempty finite subset of L, which is equivalent to a single word in L. These sets of instructions - instruction-sets - (also called instruction-information) are also indexed in the same way as the general descriptions and determine the {\it instruction paradigm} ${\cal I}_q$. Indeed, there is an injection $\r H$ on  $\r d_q$  onto ${\cal I}_q$, where $\r H(\r f^q(i,j)) = \r I^q(i,j)$ and $(i,j)$ varies over the same set of integers and natural numbers. There is one instruction paradigm for each pre-designed universe and there can be a vast collection of such universes. Rather than simply applying this bijection as a means to reproduce each of the instruction paradigm results from the developmental paradigm results, what follows is a duplicate of these results and how they are obtained in terms of instruction paradigm notation. \pars

Relative to the GGU-model and generation of a universe, a hyperfinite $\Hyper {\b I}^q(i,j)$ yields a universe-wide frozen-frame. Each instruction $x \in \Hyper {\b I}^q(i,j)$, yields a physical or physical-like system. The physical-systems are disjoint. Each collection of ultra-propertons that yields a specific physical-system is distinct from the set of ultra-propertons that yields any other physical-system. Hence, each physical-system within a universe-wide frozen-frame has a distinct identifier via the collection of all of the identifiers for the ultra-propertons or the intermediate propertons employed to produce the physical-system.\parm

\noindent{\bf 2.  Logic-System Generation for the Type-1 Interval.}\parm

The notation in all that follows is from Herrmann (2006). Notice that there are two different $t$ sequence notations. One $\r t $ is in the informal world, while another $t$ is in the formal standard superstructure. These two sequence are, of course, consider as equivalent since the set of objects that informally yield the informal $\r t$ are also formally present within the standard superstructure. The informal composition $\r f^q =\r I^q \circ \r t^q,$ when embedded relative to $\cal E$  is denoted by $\b f^q = \b I^q \circ t^q$ since the $\r t^q$ is not embedded relative to $\cal E$ and it merely generates a rational number sequence for the embedded informal paradigm.  These different notations are eliminated and only the math-italics font is employed. This is the customary practice throughout Herrmann (1979 - 1993). Notation for  informal natural, rational and real numbers, if applicable, is usually the same for the informal and more formal superstructure objects. Each $t^q(i,j)$ is a rational number. Each $\b f(i,j)$ is a nonempty instruction-set.\pars

Each member of ${\cal I}_q$ is now considered as determined by a function defined on a set $R_q$ of rational numbers, \b Q. The members of $R_q$ carry the restricted rational simple order and the order $\leq_{{\cal I}^q}$ for the members of ${\cal I}_q$ (the lexicographic order) is order isomorphic  to $R_q$ in the obvious way. Each interval partition is of the form $[c_i, c_{i+1})$ (with a closed interval in two cases), where $i \in \b Z$ and $\b Z$ is the set of integers, and $t^q(i,0) = i,\ t^q(i+1,0) =i +1.$  Then each member of $(c_i,c_{i+1})$ is a defined rational number $t^q(i,j),$ where $i < j < i+1.$ For example, consider $[c_2,c_3)$. Then $t^q(2,1) = 3- 1/2,\ t^q(2,2) = 3 - 1/4,\ t^q(2,3) =  3 - 1/8,$ then, in general, $t^q(2,j) = 3 -1/2^j.$ Hence, $\r f^q(2,0) <_{{\cal I}^q} \r f^q(2,1) <_{{\cal I}^q} \r f^q(2,2) <_{{\cal I}^q} \cdots  <_{{\cal I}^q} \r f(3,0).$ (The order $\leq_{\underline{\cal I}^q}$ is  lexicographic and is isomorphic to the rational number order for a specific set of rational numbers.)\pars

Let ${\cal I}_1$ be the standard instruction paradigm.  An instruction paradigm is defined mathematically in the exact same manner as that of the developmental paradigm in Herrmann (2006) and is equivalent to the range of a sequence $\rm g' \colon \nat\to \power {\rm L},$ where L is our denumerable general language. The first case illustrated for the GGU-model is for a developing universe starting with a frozen segments (frame) instruction-set $\r g'(0)$. For the other three GGU-model cases, this sequence is appropriately modified. In all cases, the $(\r f^q(i,j),\r f^q(p,k))$ is equivalent to ``If $\r f^q(i,j),$ then $\r f^q(p,k)).$ This notation will be simplified later. \pars

For the type-1 case $[0,b],\ b >0$, as indicated above, a denumerable instruction paradigm displays a refined form. 
For $1 < m \in \nat,\  {\cal I}_1 = \{\r f^1(i,j)\mid (0 \leq i \leq m)\land(i \in \b Z)\land (j \in \nat).$  Using ${\cal I}_1,$ consider the following logic-system.\pars

Due to the simplicity and special nature of the logic-systems used, a simplified algorithm is employed. 
The basic logic-system algorithm is re-defined for sets of two distinct objects $\{ \r A, \r B \}$. If a deduction yields $\r C$ and $\r C$ is a member of $\{ \r A, \r B \}$, then the  ``other'' member is a deduction. Hence, if A is deduced, then from $\{ \r A, \r B\}$, B is deduced. This can be written as $\{\r A, \r B \} - \{\r A\}$ is deduced. In general, this approach is only valid for these special collections of two element sets. This process mimics the proposition-logic modus ponens rule of inference $\{(\r X\to \r Y, \r X, \r Y) \mid \r X,\ \r Y \ {\rm are\ propositions}\}.$ However, for both logic-systems only one member of any two element set is deducible. \parm

{\bf Definition 2.1} Let $i \in \b Z.$ For each $n \in \nat,$ let $\r k^1_i(n) = \{ \{ \r f^1(i,j), \r f^1(i,j+1)\}\mid (0\leq j \leq n-1)\land (j\in \nat)\},\ \r K^1(n) =\bigcup \{\r k^1_i(n)\mid (0\leq i < m)\land(i \in \b Z)\}.$ Finally, let finite $\Lambda^1(n) = \{\r f^1(0,0)\} \cup \r K^1(n) \cup \{ \{ \r f^1(p-1,n),\r f^1(p,0) \}\mid (0< p \leq m)\land (p \in \b Z)\}$ and ${\cal L}^1 =\{\Lambda^1(x)\mid x \in \nat\}.$  The set $\{ \{\r f^1(p-1,n),\r f^1(p,0)\}\mid (0< p \leq m)\land (p \in \b Z)\}$ is called the ``jump elements." Also, each $\Lambda^1(n)$ is a finite set. \parm 
In general, members in ${\cal L}^q$ can be characterized by a first-order sentence.  When the deduction algorithm is applied to $\Lambda^1(n)$ the result is an ordered set of words from L -  the ordered instruction paradigm. In accordance with the juxtaposition join operator that yields words in L, this ordered instruction paradigm is a word in L. It can be obtained using the spacing symbol where each member of this paradigm is considered a sentence. For a  multi-universe theory, each such universe is a portion of each of the original members of the instruction paradigm.  \pars

In order to make the notation as simple as possible for the next construction, notice that ${\cal L}^1$ is denumerable.  Let $\nat -\{0\} = \nat'.$ Thus, there is a bijection $\r D^1 \colon \nat' \to {\cal L}^1.$ We use the subscript notation for this bijection. Thus, consider ${\cal L}^1 = \{\r D^1_i\mid i \in \nat'\}$. For each $n \in \nat'$, define $\r M^1_n = \{\{\r D^1_1,\ldots,\r D^1_n\}\}.$ Let ${\cal M}^1 =\{ M^1_n\mid n \in \nat'\}.$ The set $\r M^1_n = \{\{\r D^1_1,\ldots,\r D^1_n\}\}$, as before, can be considered as a single word-like object. \pars

(There are a few typographic errors in Herrmann (2006) and (2006a). For example, in Theorem 4.1, $m >0$ should read $m >1,$ and $\Hyper {\b D}$, should read $\Hyper {\b D}_1.$ In Herrmann (2006a), page 12, in the first (4), the $\nu \in \Hyper {\b Z}^{\geq 0} - \b Z$ should be replaced with $\nu \in \Hyper {\b Z}^{\leq 0} - \b Z$, $\gamma \in \Hyper {\b Z}^{\leq 0} - \b Z$ should be replaced with $\gamma \in \Hyper {\b Z}^{\geq 0} - \b Z$.)\pars

A finite consequence  operator S is defined in Herrmann (1979 - 1993, p. 65).  However, a new simplified logic-system ${\cal S}^q,\ q = 1,2,3,4$ is defined. When a logic-system is applied, it generates a specific finite consequence operator. It is the logic-system algorithm that does this.  In this article, this algorithm is explicitly noted since only logic-systems are used.  In general, logic-systems are stated in terms of metamathematics n-tuples. If a set  $\{\r A,\r B,\r C,\ldots,\r D\}$ is used as an hypothesis, then it is word-like since the objects the logical deduction models via the algorithm yields words or word-like objects. \pars

Define ${\cal M}^q, \ q = 2,3,4,$ in the same manner as ${\cal M}^1,$ from members of  ${\cal L}^q$. For each $\r G^q \in {\cal M}^q$, there exists a unique $n \in \nat'$ such that $\r G^q \in \r M^q_n.$  
This $\r G^q = \{\r D^q_1,\ldots,\r D^q_n\},\ \r D^q_i \in {\cal L}^q, \ 1\leq i\leq n.$  \pars 

Define the logic-system that generates $\r S^q$ as
 ${\cal S}^q = \{\{x,y\}\mid (\exists n (n \in \nat')) \land (x \in \r M^q_n) \land (y \in {\cal L}^q )\land (y \in x)\}.$ (This definition can be further described in order to  characterize the doubleton set notion and can include all necessary bounds for the quantifiers.) Further, under the simplification used here, each member of ${\cal S}^q$ is a propositional tautology. Notice that $\r M^q$ is a function with values a singleton set containing an n-set (i.e. a set of ``n'' members).\pars 
 Usually, such a logic-system would use ordered pairs  to model the rules of inference. Within these rules, finite conjunctions are displayed as first coordinates via n-sets. Again the simplified doubleton-set approach is  used here, where one of these sets is $\{\{D_1\}, D_1\}$.   \pars

Hypotheses are considered as members of a set (a 1-ary relation), when part of a logic-system. They are, usually, considered as a list of the members of this set.   In general, a logic-system, when considered as an operator,  is defined on subsets of the language employed. \pars

From the definitions employed for the logic-systems used here, the properties of the logic-system algorithm ${\cal A}$ 
can be explicitly described in set-theoretic notation. For these applications, $\cal A$ is a function defined on various defined logic-systems and a set of hypotheses. For example, the entire set of deductions or the order in which the deductions are made, among a few other characteristics. In our application to a logic-system, the notation used signifies all of the ``deduced'' results the algorithm produces when the logic-system is applied to a set of hypotheses. This yields the same results as a corresponding finite consequence operator. What the notation indicates is that the finite consequence operator is being displayed in a more refined and explicit manner. Hence, the algorithm and its relation to the logic-system can be embedded into the formal structure via formalizable characteristics.\pars

When the application characteristics are *-transferred, then the notation $\hyper {\underline{\cal A}}$ is employed. The process of applying the algorithm to the logic-system ${\cal S}^q,$ that is applied it to a set of hypotheses Y, is denoted by ${\cal A}(({\cal S}^q, \r Y)).$ Hence, $\cal A$ is defined upon a set of ordered pairs. The result of ${\cal A}(({\cal S}^q,\r Y))$ is a set. An additional step can be included for this specific algorithm, where $\r Y$ is removed. When this is done the algorithm is denoted by ${\cal A}'$. The necessary informally and, hence, formally described properties are specifically displayed. In general, the $q$ notion is not included as part of the $\cal A$ notation unless confusion would result. \pars

For the denumerable set ${\cal L}^1$, notice that for any $\Lambda^1(k),\ k \in \nat$  there exists an $k' \in \nat$ and $\r X^1_{k'} \in {\r M}^1_{k'},$ such that $\Lambda^1(k) \in {\cal A}'(({\cal S}^1,\{\r X^1_{k'}\}))$ and, in this case, finite choice yields the $\Lambda^1(k)$ logic-system. Notice that the logic-system $\Lambda^1(k)$ is considered as a set-theoretic set. 
Then the logic-system algorithm $\cal A$ is applied to  $(\Lambda^1(k),\{\r f^1(0,0)\}),$ where $\r f^1(0,0)$ is the only hypothesis contained in the logic-system. This yields $\r f^1(i,j)\in {\cal I}_1$ as a deduction from $\r f^1(0,0)$. Conversely, if $\r f^1(i,j) \in {\cal I}_1,$ then there is an $\r X^1_{k'} \in {\r M}^1_{k'}$ and a logic-system $\Lambda(k) \in {\cal A}'({\cal S}^1,\{\r X^1_{k'}\})$ such that application of the logic-system algorithm $\cal A$ to  $(\Lambda^1(k),\{\r f^1(0,0)\})$ yields $\r f^1(i,j)$ as a deduction from $\r f^1(0,0)$.  \pars

The informal algorithm $\cal A$ is defined on any logic-system that contains an hypothesis and, in this paper, such a logic-system is $\Lambda^q (x)$ and application is on $(\Lambda^q(x), \r Y)$ where Y is an hypothesis contained in the logic-system and containing but one member. Due to the construction of the $\Lambda^q (x),$ this yields a partial sequence of members of ${\cal I}_q$. This sequence is denoted by ${\cal A}[(\Lambda^q,\r Y)].$ This sequence represents the steps in the deduction and satisfies the $\leq_{{\cal I}^q_x}$ order. Also, for this case, ${\cal A}((\Lambda^q(x), \r Y)) = {\cal I}^q_{x} \subset {\cal I}_q.$ Significantly, for $n,\ k \in \nat,\ n \leq k, {\cal A}((\Lambda^1(n), \r Y)) \subset {\cal A}((\Lambda^1(k), \r Y))$ and ${\cal A}[(\Lambda^1(k),\r Y)]|[1,n] = {\cal A}[(\Lambda^1(n),\r Y)].$ \pars

In the usual way, all of the above informally defined objects are embedded relative to  $\cal E.$ When the informal set-theoretic expresses are considered as embedded into the standard superstructure, all of the bold font conventions defined in Herrmann (1979-1993) are observed. All other embedded symbols retain their math-italics form. Where script notation is used, an underline is used in place of the bold face font. All the following results are relative to our nonstandard model $\hyper {\cal M} = \langle \Hyper {\b Q}, \in , = \rangle$ or $\hyper {\cal M} = \langle \Hyper {\real }, \in ,= \rangle$ (Herrmann, (1979 -  1993)).   \parm 

\noindent{\bf Theorem 2.1} {\it Consider primitive time interval $1 = [0,b],  b>0.$ It can always be assumed that interval 1 is partitioned into two or more intervals $[c_0,c_1),\ldots$ $ [c_{m-1}, c_m], \ c_m = b, \ m >1,\ m \in \b Z.$ Let $\underline{\cal I}_1$ be an instruction paradigm order isomorphic to the rational numbers $R_1 \subset [0,b].$ For any $\lambda \in \nat_\infty,$ there exists a unique hyperfinite $\hyper {\b \Lambda^1}(\lambda) \in \Hyper {\underline{\cal L}}^1$ and a $\lambda' \in \hypernat$ such that the ultra-word-like $X^1_{\lambda'} \in  \Hyper {\b M}^1_{\lambda'}$ and ultra-logic-system $\hyper {\b \Lambda^1}(\lambda) \in  
\hyper { {\underline{\cal A}}'} ((\Hyper {\underline{\cal S}^1},\{X^1_{\lambda'}\}))$ and $\sig{\underline{\cal I}_1} \subset 
\hyper { {\underline{\cal A}}}(({\hyper {\b \Lambda^1}}(\lambda),\{\Hyper {\b f}^1(0,0)\})) = {\cal I}^1_\lambda\subset \Hyper {\underline{\cal I}}_1$. 
 Also the $\hyper { {\underline{\cal A}}}[(\hyper{\b \Lambda^1}(\lambda),\{\Hyper {\b f}^1(0,0)\})]$ *steps satisfy the $\leq_{{\cal I}^1_\lambda}$ order  and  $(\Hyper {\underline{\cal I}_1} - \sig{\underline{\cal I}_1})\cap \hyper { {\underline{\cal A}}}((\hyper {\b \Lambda^1}(\lambda),\{\Hyper {\b f}^1(0,0)\}))= $ an infinite set.} 
\pars 

Proof. This follows in the same manner as Theorem 4.1 in Herrmann (2006) by *-transfer of the appropriate first-order statements that precede this theorem statement. Also note that since for every $n \in \nat'$, the $\Lambda(n)$ is finite, then, via the identification process, $\sig {\b \Lambda(n)} =\b \Lambda(n).$ It also follows that $\Hyper {\b \Lambda(n)} = \b \Lambda(n)$ under the customary conventions. Since for any $n,\ k \in \nat',\ n \leq k$, $ {\cal A}((\b \Lambda(n),\{{\b f}^1(0,0)\})) \subset  {\cal A}((\b \Lambda(k),\{ \b f^1(0,0)\})),$  from the  above and, via *-transfer, it follows that  $^\sig{\underline{\cal I}_1} \subset  \hyper {\underline{\cal A}}((\hyper {\b \Lambda}^1(\lambda).\{\Hyper {\b f}^1(0,0)\})) = {\cal I}^1_\lambda \subset \Hyper {\underline{\cal I}}_1.$  From the definition of $\Lambda^1(n),$ these steps numbers are order isomorphic the set of rational numbers $R_1.$ Hence, $\hyper {\underline{\cal A}}((\hyper {\b \Lambda}^1(\lambda),\{\Hyper {\b f}^1(0,0)\}))$ is *order isomorphic to a hyperfinite subset of $\Hyper {\b Q}.$  
Since there are infinitely many $i < \lambda$ and $i \in \nat_\infty$,  there are infinitely many $\Hyper {\b f}(i,j) \in \hyper {\underline{\cal A}}((\hyper{\b \Lambda^1}(\lambda),\{\Hyper {\b f}^1(0,0)\}))\subset \Hyper {\underline{\cal I}_1},$ where $ \Hyper {\b f}(i,j) \in \Hyper {\underline{\cal I}_1} - \sig{\underline{\cal I}_1}.$ These are interpreted as ultranatural events but in some cases may differ from physical events only in their primitive time identifications. This completes the proof. \qed 

By considering the definition of ${\cal L}^1$, it follows that the given $1 < m \in \nat,\ \hyper {\b \Lambda^1}(\lambda)$ is precisely 
$\{\Hyper {\b f}^1(0,0)\} \cup \{\bigcup\{\Hyper{\b k}^1_i\mid 0 \leq i <m\}\} \cup \{\{\Hyper {\b f}^1(p-1,\lambda),\Hyper {\b f}^1(p,0)\}\mid (0< p\leq m) \land (p \in \Hyper {\b Z})\}.$ Of significance is the fact that the steps in the *-deduction $\hyper { {\underline{\cal A}}}((\hyper{\b \Lambda^1}(\lambda),\{\Hyper {\b f}^1(0,0)\}))$ preserve the order $\leq_{\Hyper {\underline{\cal I}}_1}$. Notice that $\hyper {\b \Lambda^1}(\lambda)$ is obtained by hyperfinite choice.  Further, any $\Hyper {\b f}^1(i,j) \in \{\Hyper {\b f}^1(x,y)\mid (0\leq x <m)\land(0\leq y \leq \lambda)\land(x \in \Hyper {\b Z})\land( y \in \hypernat)\} \cup \{\Hyper {\b f}^1(m,0)\}$ is a hyperfinite *-deduction from $\b f^1(0,0) = \Hyper {\b f}^1(0,0).$ And, it also follows that the set of all such *deductions yields a hyperfinite set ${\cal I}^1_\lambda$ such that $\sig{\underline{\cal I}_1} \subset {\cal I}^1_\lambda \subset \Hyper {\underline{\cal I}_1}.$\pars

For the GGU-model and each of the four cases, an internal nonempty set X disjoint from $\bigcup \{\hyper {\cal F}^i({\cal C})\mid (0\leq i \leq n)\land (i \in \nat)\}$ is adjoined when info-fields are employed. It is assumed that a physical universe is a collection of many ($> 1$) physical-systems. By the way each universe-wide frozen frame is constructed, each $\Hyper {\b f}^q(i,j)$ is an infinite hyperfinite set that contains for standard $(i,j)$, the instruction-set $\b f^q(i,j),$ and the set $\Hyper {\b f}^q(i,j) - \b f^q(i,j) \not= \emptyset.$  Members of $\Hyper {\b f}^q(i,j) - \b b^q(i,j)$ can yield the same physical-systems or members can yield physical-like systems. Each physical-system corresponds to one member of $\Hyper {\b f}^q(i,j).$ \pars

Consider each $\Hyper {\b f}^q(i,j) \in \sig{\underline{\cal I}_q}.$ There is a function $G^q_x(i,j) \colon \Hyper {\b f}^q(i,j) \to \power{ \bigcup \{\hyper {\cal F}^i({\cal C})\mid (0\leq i \leq n)\land (i \in \nat)\} \cup X}.$ The image of $G^q_x(i,j)$ is an info-field. Each  $\hyper y \in \Hyper {\b f}(i,j) $ determines the properton composition for each of the physical-systems. For a $z \in \Hyper {\b f}(i,j)$ that is not equivalent to a $\hyper y,$ $G^q_x(i,j)(z)$ determines a physical-like system relative to X. \parm

\noindent{\bf 3. Logic-System Generation for the Type-2 Interval}\parm

For the type-2 case $[0,+\infty)$, a denumerable instruction paradigm displays a refined form. 
For this case, ${\cal I}_2 = \{\r f^2(i,j)\mid (0 \leq i )\land(i \in \b Z)\land (j \in \nat).$  Using ${\cal I}_2,$ consider the following logic-system.\parm

{\bf Definition 3.1} Let $0 \leq i \in \b Z.$ For each $n \in \nat,$ let $\r k^2_i(n) = \{\{\r f^2(i,j), \r f^1(i,j+1)\}\mid (0\leq j \leq n-1)\land (j\in \nat)\}.$ For $0< m \in \b Z,$ let $\r K^2(m,n) =\bigcup \{\r k^2_i(n)\mid (0\leq i < m)\land(i \in \b Z)\}.$ Finally, let $\Lambda^2(m,n) = \{\r f^2(0,0)\} \cup \r K^2(m,n) \cup \{\{\r f^2(p-1,n),\r f^2(p,0)\}\mid (0< p\leq m)\land (p \in \b Z)\}\cup \{\{\r f^2(m,j), \r f^2(m,j+1)\}\mid (0\leq j < n)\land (j \in \nat)\},$ and ${\cal L}^2 = \{\Lambda^2(x,y)\mid (0 \leq x \in \b Z)\land (y \in \nat)\}.$ Notice that if $0 \leq i < k, \ i,\ k \in \b Z$, then ${\cal A}((\Lambda^2(i,j), \{\r f^2(0,0)\})) \subset {\cal A}((\Lambda^2(k,n),\{\r f^2(0,0)\}))$ for any $j,\ n \in \nat.$ 
Also, each $\Lambda^2(m,n)$ is a finite set. (Notice that members in ${\cal L}^2$ can be characterized by a first-order sentence.)\parm

Consider any $\Lambda^2(q,k).$ Then there exists an $q'k'\ \in \nat'$ ($q'k'$ is a natural number in $\nat'$) and the 
$q'k'$-set $\r X^2_{q'k'} \in \r M^2_{q'k'},$ such that $\Lambda^2(q,k) \in {\cal S}^2(\{\r X^2_{q'k'}\})$ and, in this case, finite choice yields the $\Lambda^2(q,k)$ logic-system. Then the logic-system algorithm $\cal A$ applied to $(\Lambda^2(q,k),\{\r f^2(0,0)\})$ yields $\r f^2(q,k)$ as a deduction from $\r f^2(0,0).$ Further, $\r f^2(q,k) \in {\cal I}_2.$ Conversely, if $\r f^2(q,k) \in {\cal I}_2,$ then there exists an $q'k' \in \nat'$ and an $\r X^2_{q'k'} \in {\r M}^2_{q'k'}$ and a logic-system $\Lambda(q,k) \in { {\cal A}'}(({\cal S}^2,\{\r X^2_{q'k'}\}))$ such that application of the logic-system algorithm $\cal A$ to  $(\Lambda^2(q,k),\{\r f^2(0,0)\})$ yields a deduction of $\r f^2(q,k)$ from $\r f^2(0,0)$. \parm

\noindent{\bf Theorem 3.1} {\it Consider primitive time interval $2 = [0,+\infty).$ It can always be assumed that interval 2 is partitioned  into intervals $[c_0,c_1),\ldots [c_{m-1}, c_m), \ m >1,\ m \in \b Z.$ Let $\b d_2$ be an instruction paradigm order isomorphic to the rational numbers $R_2 \subset [0,+\infty).$ For any $\lambda \in \nat_\infty$ and $\nu \in 
\Hyper {\b Z} - \b Z, \ \nu >0,$ there exists a unique hyperfinite $\hyper {\b \Lambda^2}(\nu,\lambda) \in \Hyper {\underline{\cal L}^2}$ and $\nu', \lambda'  \in \hypernat$ such that the ultra-word-like $X^2_{\nu'\lambda'} \in \Hyper {\b M}^2_{\nu'\lambda'}$ and ultra-logic-system $\hyper {\b \Lambda^2}(\nu,\lambda) \in 
\hyper { {\underline{\cal A}}'}((\Hyper {{\cal S}^2,\{X^2_{\nu'\lambda'}\}}))$ and $\sig{\underline{\cal I}_2} \subset \hyper { {\underline{\cal A}}}((\hyper {\b \Lambda^2}(\nu,\lambda),\{\Hyper {\b f}^2(0,0)\})) = 
{\cal I}^2_{\nu\lambda}\subset \Hyper {\underline{\cal I}}_2$.
Also the $\hyper { {\underline{\cal A}}}[(\hyper {\b \Lambda^2}(\nu,\lambda),\{\Hyper {\b f}^2(0,0)\}))]$ *steps satisfy the $\leq_{{\cal I}^2_{\nu\lambda}}$ order  and  $(\Hyper {\underline{\cal I}_2} -  \sig{\underline{\cal I}_2})\cap \hyper { {\underline{\cal A}}}(\hyper {\b \Lambda^2}(\nu,\lambda),\{\Hyper {\b f}^2(0,0)\})) = $ an infinite set.} 
\pars
Proof. As in Theorem 2.1, the proof follows by *-transfer of the appropriate formally presented material that appears above in this section 3.\qed

By considering the definition of ${\cal L}^2$, it follows that the $\hyper {\b \Lambda^2}(\nu,\lambda)$ is precisely 
$ \{\Hyper {\b f}^2(0,0)\} \cup \{\bigcup\{\Hyper {\b k}^2_i\mid 0 \leq i< \nu\}\} \cup \{(\{\Hyper {\b f}^2(p-1,\lambda),\Hyper {\b f}^2(p,0)\}\mid (0 < p\leq\nu)\land (p \in \Hyper {\b Z})\} \cup \{\{\Hyper {\b f}^2(\nu,j),\Hyper {\b f}^2(\nu,j+1)\}\mid (0\leq j < \lambda)\land(j \in \hypernat)\}$. Of significance is the fact that the steps in the *-deduction $\hyper { {\underline{\cal A}}}[\hyper {\b \Lambda^2}(\nu,\lambda),\{\Hyper {\b f}^2(0,0)\}))]$ preserve the order $\leq_{\Hyper {\underline{\cal I}}_2}$.  Notice that $\hyper {\b \Lambda^2}(\nu,\lambda)$ is obtained by hyperfinite choice. Further, any $\Hyper {\b f}^2(i,j) \in \{\Hyper {\b f}^2(x,y)\mid (0\leq x \leq \nu)\land(0\leq y \leq \lambda)\land(x \in \Hyper {\b Z})\land (y \in \hypernat)\}$ is a hyperfinite *-deduction from ${\b f}^2(0,0).$ And, it also follows that the set of all such *deductions yield a hyperfinite set ${\cal I}^2_{\nu\lambda}$ such that $\sig{\underline{\cal I}_2} \subset {\cal I}^2_{\nu\lambda} \subset \Hyper {\underline{\cal I}_2}.$\parm

\noindent{\bf 4. Logic-System Generation for the Type-3 Interval}\parm

For the type-3 case $(-\infty,0]$, a denumerable instruction paradigm displays a refined form. 
For this case, $\r d_3 = \{\r f^3(i,j)\mid (i\leq 0 )\land(i \in \b Z)\land (j \in \nat).$  Using $\r d_3,$ consider the following logic-system.\parm

{\bf Definition 4.1} Let $i \in \b Z,\ i\leq 0.$ For each $n \in \nat,$ let $\r k^3_i(n) = \{\{\r f^2(i,j), \r f^1(i,j+1)\}\mid (0\leq j \leq n-1)\land (j\in \nat)\}.$ For $m \in \b Z\, m<0,$ let $\r K^3(m,n) =\bigcup \{\r k^3_i(n)\mid (m \leq i < 0)\land(i \in \b Z)\}.$ Finally, let $\Lambda^3(m,n) = \{\r f^3(m,0)\} \cup \r K^3(m,n) \cup \{\{ \r f^3(p-1,n),\r f^3(p,0)\}\mid (m < p \leq 0)\land (p \in \b Z)\},$ and ${\cal L}^3 =\{\Lambda^2(x,y)\mid (0 \leq x \in \b Z)\land (y \in \nat)\}.$ Notice that if $ i < k \leq 0, \ i,\ k \in \b Z$, then ${\cal A}((\Lambda^3(i,j)),\{\r f^3(m,0)\})) \subset {\cal A}((\Lambda^3(k,n), \{\r F^3(m,0)\}))$ for any $j,\ n \in \nat.$ 
 Also, each $\Lambda^3(m,n)$ is a finite set. (Notice that members in ${\cal L}^3$ can be characterized by a first-order sentence.)\parm

Consider any $\Lambda^3(q,k).$ Then there exists an $q'k' \in \nat$ and $\r X^3_{q'k'} \in \r M^3_{q'k'},$ such that $\Lambda^3(q,k) \in {\cal A}'(({\cal S}^3,\{\r X^3_{q'k'}\}))$ and, in this case, finite choice yields the $\Lambda^3(q,k)$ logic-system. Then the logic-system algorithm $\cal A$ applied to $(\Lambda^3(q,k),\{\r f^3(q,0)\}))$ yields $\r f^3(q,k)$ as a deduction from $\r f^3(q,0).$ Further, $\r f^3(q,k) \in {\cal I}_3.$ Conversely, if $\r f^3(q,k) \in \r {\cal I}_3,$ then there is an $\r X^3_{q'k'} \in \r {\cal M}^3_{q'k'}$ and a logic-system $\Lambda(q,k) \in \r S^3(\{\r X^3_{q'k'}\})$ such that application of the logic-system algorithm $\cal A$ to  $(\Lambda^3(q,k),\{\Hyper {\b f}^3(q,0)\}))$ yields $\r f^3(q,k)$ as a deduction from $\r f^3(q,0)$. \parm

\noindent{\bf Theorem 4.1} {\it Consider primitive time interval $3 = (-\infty, 0].$ It can always be assumed that interval 3 is partitioned  into intervals $ \ldots, [c_{-2},c_{-1}), [c_{-1},c_0]$. Let $\b d_3$ be an instruction  paradigm order isomorphic to the rational numbers $R_3 \subset (-\infty,0].$ 
For any $\lambda \in \nat_\infty,\ \mu \in \Hyper {\b Z} - \b Z,\ \mu < 0,$
 there exists a unique hyperfinite $\hyper {\b \Lambda^3}(\mu,\lambda) \in \Hyper {\underline{\cal L}^3}$ and $\mu', \lambda'  \in \hypernat$ such that the ultra-word-like $X^3_{\mu'\lambda'} \in \Hyper {\b M}^3_{\mu'\lambda'}$ and ultra-logic-system 
$\hyper {\b \Lambda^3}(\mu,\lambda) \in 
\hyper { {\underline{\cal A}}'}((\Hyper {{\cal S}^3,\{X^3_{\mu'\lambda'}\}}))$ and $\sig{\underline{\cal I}_3} \subset \hyper { {\underline{\cal A}}}(\hyper {\b \Lambda^3}(\mu,\lambda),\{\Hyper {\b f}^3(\mu,0)\})) = {\cal I}^3_{\mu\lambda}\subset \Hyper {\underline{\cal I}}_3$.
Also the $\hyper { {\underline{\cal A}}}[(\hyper {\Lambda^3}(\mu,\lambda),\{\Hyper {\b f}^3(\mu,0)\}))]$ *steps 
satisfy the $\leq_{{\cal I}^3_{\mu\lambda}}$ order  and  $(\Hyper {\underline{\cal I}_3} - 
^{\underline{\cal I}_3})\cap \hyper { {\underline{\cal A}}}((\hyper {\b \Lambda^3}(\mu,\lambda)) ,\{\Hyper {\b f}^3(\mu,0)\}))= $ an infinite set.} 
\pars

Proof. As in Theorem 3.1, the proof follows by *-transfer of the appropriate formally presented material that appears above in this section 3.\qed

By considering the definition of ${\cal L}^3$, it follows that the $\hyper {\b \Lambda^3}(\mu,\lambda)$ is precisely 
$ \{\Hyper {\b f}^3(\mu,0)\} \cup \{\bigcup\{\Hyper {\b k}^3_i\mid \mu \leq i < 0\}\} \cup \{\{\Hyper {\b f}^3(p-1,\lambda),\Hyper {\b f}^2(p,0)\}\mid (\mu < p \leq 0)\land (p \in \Hyper {\b Z})\} $. Of significance is the fact that the steps in the *-deduction $\hyper {\underline{\cal A}}[(\hyper {\b \Lambda^3}(\mu,\lambda), \{\Hyper {\b f}^3(\mu,0)\}]$ preserve the order $\leq_{\Hyper {\underline{\cal I}}_3}$. Notice that $\hyper {\b \Lambda^3}(\mu,\lambda)$ is obtained by hyperfinite choice.  Further, any $\Hyper {\b f}^3(i,j) \in \{\Hyper {\b f}^3(x,y)\mid (\mu \leq x < 0)\land(0\leq y \leq \lambda)\} \cup \{{\b f}^3(0,0\}$ is a hyperfinite *-deduction from $\b f^3(\mu,0).$ And, it also follows that the set of all such *deductions is a hyperfinite set ${\cal I}^3_{\nu\lambda}$ such that $\sig{\underline{\cal I}_3} \subset {\cal I}^3_{\nu\lambda} \subset \Hyper {\underline{\cal I}_3}.$\parm

\noindent{\bf 5. Logic-System Generation for the Type-4 Interval}\parm

\noindent{\bf Theorem 5.1} {\it Consider primitive time interval $4 = (-\infty, +\infty).$ It can always be assumed that interval 4 is partitioned  into intervals $ \ldots, [c_{-2},c_{-1}), [c_{-1},c_0), \ldots$. Let $\b d_4$ be a instruction paradigm order isomorphic to the rational numbers $R_4 \subset (-\infty,+\infty).$ For any $\lambda \in \nat_\infty,\ \nu, \gamma \in \Hyper {\b Z} - \b Z,$ such that $\nu \leq 0, \ \gamma\geq 0,$ 
 there exists a unique hyperfinite $\hyper {\b \Lambda^4}(\nu,\gamma,\lambda) \in \Hyper {\underline{\cal L}^4}$ and $\nu', \gamma', \lambda'  \in \hypernat$ such that the ultra-word-like $X^4_{\nu'\gamma'\lambda'} \in \Hyper {\b M}^4_{\nu'\gamma'\lambda'}$ and ultra-logic-system $\hyper {\b \Lambda^4}(\nu,\gamma,\lambda) \in 
\hyper { {\underline{\cal A}}'}((\Hyper {{\cal S}^4,\{X^4_{\nu'\gamma'\lambda'}\}}))$ and $\sig{\underline{\cal I}_4} \subset \hyper { {\underline{\cal A}}}((\hyper {\b \Lambda^4}(\nu,\gamma,\lambda),\{\Hyper {\b f}^4(\nu,0)\})) = {\cal I}^4_{\nu\gamma\lambda}\subset \Hyper {\underline{\cal I}}_4$. Also the $\hyper {\underline{\cal A}}[(\hyper {\b \Lambda^4}(\nu,\gamma,\lambda),\{\Hyper {\b f}^4(\nu,0)\})]$ *steps  
satisfy the $\leq_{{\cal I}^4_{\nu\gamma\lambda}}$ order  and $(\Hyper {\underline{\cal I}_4} - \sig{\underline{\cal I}_4})\cap \hyper {\underline{\cal A}}((\hyper {\b \Lambda^4}(\nu,\gamma,\lambda),\{\Hyper {\b f}^4(\nu,0)\}))= $ an infinite set.} \parm 

By considering the definition of ${\cal L}^4$, it follows that the $\hyper {\b \Lambda^4}(\nu,\gamma,\lambda)$ is precisely 
$ \{\Hyper {\b f}^4(\nu,0)\} \cup  \{\bigcup\{\Hyper {\b k}^4_i\mid (\nu \leq i < \gamma) \land (i \in \Hyper {\b Z})\}\} \cup  \{\{\Hyper {\b f}^4(p-1,\lambda),\Hyper {\b f}^4(p,0)\}\mid (\nu<p\leq \gamma)\land(p \in \Hyper {\b Z})\} \cup\{\{\Hyper {\b f}^4(\gamma,j),\Hyper {\b f}^4(\gamma,j+1)\}\mid (0\leq j < \lambda)\land (j \in \hypernat)\}.$ Of significance is the fact that the steps in the *-deduction $\hyper {\underline{\cal A}}[(\hyper {\b \Lambda^4}(\nu,\gamma,\lambda),\{\Hyper {\b f}^4(\nu,0)\})]$ preserve the order $\leq_{\Hyper {\underline{\cal I}}_4}$. Notice that $\hyper {\b \Lambda^4}(\nu,\gamma,\lambda)$ is obtained by hyperfinite choice.  Further, any $\Hyper {\b f}^4(i,j) \in \{\Hyper {\b f}^4(x,y)\mid (\nu \leq x \leq \gamma)\land(0\leq y \leq \lambda)\}$ is a hyperfinite *-deduction from $\b f^4(\nu,0).$ And, it also follows that the set of all such *deductions is a hyperfinite set 
${\cal I}^4_{\nu\gamma\lambda}$ such that $\sig{\underline{\cal I}_4} \subset {\cal I}^4_{\nu\gamma\lambda} \subset \Hyper {\underline{\cal I}_4}.$\pars
The above established theorems and appropriate definitions all hold for the development paradigms and yield ultra-logic-systems that can replace the ultralogic notion.\parm

\noindent{\bf 6. The Complete GGU-model Scheme}\parm

For the $({\tt St}G_q)$ is defined in Herrmann (2006a). The following scheme is not in composition notational form due to one application of choice and the step-by-step application of $({\tt St}G_q).$ It represents an ordered application of the GGU-model operators. For $q = 1,2,3,4,$ $x = \lambda, \nu\lambda,\mu\lambda, \nu\gamma\lambda$ with or without commas, respectively. The $a,b,c$ take the appropriate value for a specific $q.$

$$({\tt St}G_q)(\hyper {\underline{\cal A}}((\hyper {\b \Lambda^q}(a),\{\Hyper {\b f}^q(b,c)\})))(\hyper {\underline{\cal A}'}((\Hyper {{\cal S}^q,\{X^q_{a'}\}}))).$$

The operators $\underline{\cal A}$ and $\underline{\cal A}'$ have characterizing first-order statements. These statements need not capture all of  of the intuitive statements that describe the algorithms. T he results of application of $\underline{\cal A}'$ as formalized can show major aspects of the algorithm's selection process. For example,
$$\forall x \forall y\forall z\forall w((w \in \underline{\cal M}^q)\land(y \in {\cal F}(\underline{\cal L}^q)) \land(y\in w)\land (x \in  \underline{\cal A}'((\underline{\cal S}^q,y)) \to$$ 
$$\exists p((p \in \underline{\cal S}^q)\land (y \in p)\land (x \in p)\land (y \not= x)).$$
Of course, the natural numbers and the embedded $\b D^q$ can also be employed.\parm

\noindent{\bf 7. The Participator Universe.}\parm

For the GGU-model, one of the most difficult requirements is to include the concept of the ``participator'' universe. As stated at the May 1974 Oxford Symposium in Quantum Gravity, Patton and Wheeler describe how existence of human beings alter the universe to various degrees. ``To that degree the future of the universe is changed. We change it. We have to cross out that old term `observed' and replace it with the new term `participator.' In some strange sense the quantum principle tells us that we are dealing with a participator universe." (Patton and Wheeler (1975, p. 562).) This aspect of the GGU-model is  only descriptively displayed in section 4.8 in Herrmann (2002). It is now possible to obtain formally the collection of pre-designed universes that satisfies this participator requirement. \pars

The previous notation is modified for finitely many ($> 0$) instruction paradigms as previously denoted by ${\cal I}_q.$ From the construction of each instruction paradigm using  \r L, it follows that there is, at least, a sequence of possible alterations. An instruction paradigm is a nonempty subset of the subsets of \r L. Hence, the collection of all such instruction paradigms is a member of $\power{\power {\r L}}$. There can be infinitely many basic universes. These are universes prior to participator alterations. For each of these, there is an collection of ultra-word-like objects of the appropriate type. What follows next is for an arbitrary member of  this collection of ultra-world-like objects and, hence, an arbitrary basic universe. \pars

So as to include the type of universe being considered, let $q \colon  \nat' \times [1,4] \to \power{\power {\r L}}.$ Then for a specific $p \in [1,4]$ an expression $\{x \mid (x = \r q(n,p)) \land (n \in \nat')\} = {\cal I}_{p}$ represents this denumerable set of instruction paradigms for type-p universes. (If $n' \not= m ,$ then $q(n',p) \not= q(m,p).$)  Let  $\{x \mid \exists p (p \in [1,4] \land(x =  {\cal I}_{p})\} =   {\cal I}.$  Then, a  specific $\r f^1(0,0)$ is further identified relative to the sequences. As an example, $\r f^{\r q(3,1)}(0,0) \subset \r L$ represents a specific $n =3$ type-1  instruction-set. Thus, as embedded into the formal structure this last expression reads $\b f^{\b q(3,1)}(0,0) \subset \b L$.\pars

An original alteration can be miniscule and made in one or more of the necessary parameters that are satisfied by a specific cosmology. This can be done in such a way that only miniscule alterations in physical-system satisfy  the alterations. On the other hand, a highly altered cosmology can also occur. An alteration is local prior to it being propagated during a universe's development. Although not specifically included, and indeed the definition would need to be altered slightly, other sequences $\r q'$ can be used where various values can be empty or repeated. Also, for a universe with infinitely many local alterations at the same moment, one can include the obvious change in the $\r q$ sequence, where the sequence $q$ is a sequence of type-4 except that each image is a ``universe.'' \pars

For the GGU-model, the various members of ${\cal I}_p$ satisfy the ``participator" requirements, when participators exist, for each of the known suggested cosmologies. When embedded into the formal structure, properties of $q$ can be easily characterized, using various forms, in a first order language. For example, 
$$\forall x \forall p((x \in \nat') \land(p \in [1,4]) \to $$ $$\exists y\exists z((z \in \underline{{\cal I}})\land (y \in z) \land (\b q(x,p) = y)).$$ Consider a specific $p' \in [1,4]$. Then 
$$\forall y((y \in \underline{{\cal I}}_{p'})) \to \exists x((x \in \nat') \land (\b q(x,p') =y)).$$ 
In *-transfer form, these two sentences read
$$\forall x \forall p((x \in \hypernat') \land(p \in [1,4]) \to $$ $$\exists y\exists z((z \in \Hyper{\underline{{\cal I}}})\land (y \in z) \land (\hyper{\b q(x,p)} = y)).$$
$$\forall y((y \in \Hyper {\underline{{\cal I}}}_{p'})) \to \exists x((x \in \hypernat') \land (\hyper{\b q}(x,p') =y)).$$ 

The previous four theorems are all relative to a specific instruction paradigm and each holds for a collection of these instruction paradigms. Thus, the notion can be added and the additional statement that the results hold for each $n \in \hypernat'$ and each $p \in [1,4].$ The special processes noted in the scheme in section 6 are applied to each set of instruction paradigms.  \pars 

Each member of In Herrmann (2002), hyperfast propertons are mentioned as mediators for the automatic selection of a particlar member of ${\cal I}_{p}.$ Each not realized member of $\{{\cal I}\}_p$ is termed as covirtual. Notice that from the transformed formal statements above, there exist a member of  ${\underline{{\cal I}}}_{p'}$ for each $\gamma \in \hypernat' - \nat'.$ These can be used for various interpretations using either a $\r q,\ \r q',\ q$ type of sequence. Further, GGU-model predicted processes and entities can aid in comprehening the notion of the non-temporal and its relation to the temporal.  \parm

\centerline{\bf References}\par\bigskip
\noindent{H}errmann, R. A., (2006), The GGU-model and Generation of Developmental Paradigms, http://arxiv.org/abs/math/0605120 \pars
\noindent{H}errmann, R. A., (2006a), General Logic-systems that Determine Significant Collections of Consequence Operators, \hfil\break http://arxiv.org/abs/math/0603573 \pars
\noindent{H}errmann, R. A., (2003), The Theory of Infinitesimal Light-Clocks, \hfil \break http://arxiv.org/abs/math/0312189 \pars
\noindent{H}errmann, R. A., (2002), Science Declairs that Our Universe IS Intelligently Designed, Xulon Press, Longwood, FL. (and elsewhere)\pars
\noindent{H}errmann, R. A., (1979 - 1993). The Theory of Ultralogics, \hfil\break  
http://raherrmann.com/books.htm \pars
http://arxiv.org/abs/math/9903082 http://arxiv.org/abs/math/9903081 
\noindent{H}errmann, R. A., (1989), Fractals and ultrasmooth microeffects, J. Math. Physics, 30(4), :805-808. (Note that there are typographical errors in this paper. In the proof of Theorem 4.1, in equations $h(x,c,d), \ G_j(x),$ the ) + )) should be )) + ). In $G_j(x),$ the second $c$ should be replaced with $a_j.$ On page 808, the second column, second paragraph, line six, $\st {D}$ should read $\st {\hyper D}$ and, trivially, $x \in \mu (p)$, should read $x \in \mu(p) \cap \hyper {D}.$ In the proof of Theorem 3.1, first line $\Hyper {\bf R}^m$ and should read ${\b R}^m.$\pars
\noindent{H}errmann, R. A., (1983). Mathematical philosophy and developmental processes, Nature and System 5(1/2):17-36.\pars
\noindent{L}oeb, P. A. and M. Wolff, (2000). Nonstandard Analysis for the Working Mathematician, Kluwer Academic Publishers, Boston.\pars

\noindent{P}aton, C. M. and J. A. Wheeler, (1975). Is physics legislated by a cosmogony? In {\it Quantum Gravity} ed. Isham, Penrose and Sciama, Clarendon Press, Oxford. \pars
\noindent{S}troyan, K. D. and J . M. Bayod, (1986). Foundations of Infinitesimal Stochastic Analysis, North Holland, N.Y.\pars
\vfil\eject
\newdimen\fullhsize
\fullhsize=4.5in \hsize=2.0in 
\def\fullline{\hbox to \fullhsize}
\def\makeheadline{{\vbox to 0pt{\vskip-22.5pt
  \fullline{\vbox to 8.5pt{}\the\headline}\vss}\nointerlineskip}}
\headline={\ifnum\pageno=100\hfil\tenbf\folio\else\ifodd\pageno\rightheadline 
\else\leftheadline\fi\fi}
\def\rightheadline{\sl\hfil SYMBOLS \hfil\tenbf\folio}
\def\leftheadline{\tenbf\folio\hfil Robert A. Herrmann \hfil}
\voffset=1\baselineskip

\output={\if L\lr
         \global\setbox\leftcolumn=\columnbox \global\let\lr=R
\else \doubleformat \global\let\lr=L\fi
\ifnum\outputpenalty>-20000 \else\dosupereject\fi}
\def\doubleformat{\shipout\vbox{\makeheadline
\fullline{\box\leftcolumn\hfil\columnbox}
\makefootline}
\advancepageno}
\def\columnbox{\leftline{\pagebody}} 
\let\lr=L \newbox\leftcolumn

\noindent Symbol\dotfill Page no.\pars
\noindent$A_h$\dotfill 7.\pars 
\noindent$H_t$\dotfill 7.\pars
\noindent${\cal A}$\dotfill 7.\pars
\noindent$\cal W$\dotfill 7.\pars
\noindent{\bf ZFC}\dotfill 8.\pars 
\noindent{\bf ZFH}\dotfill 8.\pars
\noindent$D$\dotfill 9.\pars
\noindent$W_Y$\dotfill 9.\pars
\noindent$B_Y$\dotfill 9.\pars
\noindent${\cal P}(B_Y)$\dotfill 9.\pars
\noindent$H^n= H^{[0,n]}$\dotfill 9.\pars 
\noindent$P_H$\dotfill 10.\pars
\noindent$P_{i[{\cal W}]} = P$\dotfill 10.\pars
\noindent$f(k) \leq_f f(j)$\dotfill 10.\pars
\noindent$\sim$\dotfill 10.\pars
\noindent$[f]$\dotfill 11.\pars
\noindent$\vert \cdot\vert$\dotfill 11.\pars
\noindent${\cal E}$\dotfill 11\pars
\noindent$C$\dotfill 11.\pars
\noindent$F(A)$\dotfill 11.\pars
\noindent$C\colon \power A \to \power A$\dotfill 11.\pars
\noindent$k\subset F(A) \times A$\dotfill 12.\pars
\noindent$C_k$\dotfill 12.\pars
\noindent$k_c$\dotfill 13.\pars
\noindent$\Theta \colon {\cal W} \to {\cal E}$\dotfill 13.\pars 
\noindent$\cal X$\dotfill 17.\pars
\noindent$\theta \colon i[{\cal W}] \to {\cal E}$\dotfill 20.\pars
\noindent${\bf R_A}$\dotfill 21.\pars
\noindent${\cal M}$\dotfill 22.\pars
\noindent$\hyper {{\cal M}} = \langle \hyper {{\cal N}},\in, = 
\rangle$\dotfill 23.\pars
\noindent$^\sigma B$\dotfill 23.\pars
\noindent${\cal Y}$\dotfill 23.\pars
\noindent${\rm L}_0$\dotfill 24.\pars
\noindent$\{{\rm P}_i\mid i \in \omega \}$\dotfill 25.\pars
\noindent${\rm B}^\prime$\dotfill 26.\pars
\noindent${\rm B}$\dotfill 26.\pars
\noindent${\rm C}_1$\dotfill 26.\pars
\noindent${\rm BP}$\dotfill 26.\pars
\noindent$j \colon {\rm B} \to \{{\rm P}_i\mid i \in \omega \}$\dotfill 26.\pars
\noindent${\rm BP}_0$\dotfill 26.\pars
\noindent${\rm E}_0$\dotfill 26.\pars
\noindent$\vdash_a$\dotfill 29.\pars
\noindent$T = i[{\cal W]}$\dotfill 30.\pars
\noindent${\rm S}$\dotfill 32, 65.\pars
\noindent$\rm S_0$\dotfill 33.\pars
\noindent${\rm BP}_{0(}$\dotfill 34.\pars
\noindent(MP)\dotfill 34.\pars
\noindent${\rm Rn}(f^\prime)\ D(f^\prime)$\dotfill 34.\pars
\noindent$\hyper {P_{L_1}}$\dotfill 35.\pars
\noindent${\rm BPC|}$\dotfill 37.\pars
\noindent$\Pi_W$\dotfill 37.\pars
\noindent$f_w$\dotfill 37.\pars
\noindent${\rm C}$\dotfill 38.\pars
\noindent$F_\nu,\ f^b$\dotfill 38.\pars
\noindent$\vdash_\pi$\dotfill 38.\pars
\noindent$\Pi$\dotfill 38.\pars
\noindent${\rm BPC}_0$\dotfill 38.\pars
\noindent$g^{b}$\dotfill 39.\pars
\noindent$G_n$\dotfill 39.\pars
\noindent$\#(\rm A)$\dotfill 40.\pars
\noindent${\cal G}$\dotfill 41.\pars
\noindent$\rm B\leq_\# \rm D$\dotfill 43.\pars
\noindent$[f] \leq_B [g]$\dotfill 44.\pars
\noindent$A \leq_B D$\dotfill 44.\pars
\noindent${\cal C}^\prime$\dotfill 45.\pars
\noindent${\cal C}_f^\prime$\dotfill 45.\pars
\noindent$C_1 \leq C_2$\dotfill 45.\pars
\noindent$C_1 \wedge C_2)$\dotfill 45.\pars
\noindent${\cal C} = {\cal C}^\prime - \{U\}$\dotfill 46.\pars
\noindent${\cal C}_f = {\cal C}_f^\prime - \{U\}$\dotfill 46.\pars
\noindent$K$\dotfill 47.\pars
\noindent$K_\infty$\dotfill 47.\pars
\noindent${\rm MP}_n$\dotfill 49.\pars
\noindent${\cal H} \vdash_n X$\dotfill 49.\pars
\noindent$\lceil \ \rceil $\dotfill 53.\pars
\noindent$(A)_R$\dotfill 54.\pars
\noindent$RR(E,A)$\dotfill 54.\pars
\noindent$TC$\dotfill 54.\pars
\noindent$SS$\dotfill 54.\pars
\noindent${\cal M}_1$\dotfill 54.\pars
\noindent${\cal M}_2$\dotfill 54.\pars
\noindent$\hyper p$\dotfill 55.\pars
\noindent$\hyper{\hyper {\cal N}}$\dotfill 55.\pars
\noindent$\hyper {\underline{C}}$\dotfill 57.\pars
\noindent$\rm T_i$\dotfill 61.\pars
\noindent$\cal T$\dotfill  63.\pars 
\noindent$\widetilde{\cal T}$\dotfill 63.\pars
\noindent${\cal S}$\dotfill 63.\pars
\noindent$B$\dotfill 65.\pars
\noindent$\rm P_0$\dotfill 65.\pars
\noindent$M$\dotfill 66.\pars
\noindent$\rm d$\dotfill 66.\pars
\noindent$\rm P_m$\dotfill 67.\pars
\noindent$M_B$\dotfill 68.\pars
\noindent$N(G)$\dotfill 69.\pars
\noindent${\cal B}$\dotfill 69.\pars
\noindent EGS\dotfill 70.\pars
\noindent${\cal F}(A,B)$\dotfill 70.\pars
\noindent$D(a,b)$\dotfill 72.\pars
\noindent$\cal Q$\dotfill 77.\pars 
\noindent$\monad 0^+$\dotfill 77.\pars
\noindent$M(\cal Q)$\dotfill 78.\pars
\noindent${\rm G_A}$\dotfill 82.\pars
\noindent$\cal O$\dotfill 91.\pars
\noindent$H_0(A)$\dotfill 93.\pars
\noindent$H_1(A)$\dotfill 93.\pars
\noindent$C_i$\dotfill 93.\pars
\noindent$\{a_i\}$\dotfill 99.\pars

\end

\end

\end